\documentclass[10pt]{article}
\usepackage{}
\usepackage{amssymb}
\usepackage{mathrsfs}
\usepackage{amsmath,amssymb,theorem}
\usepackage{graphicx}
\usepackage{color}
\usepackage{ifpdf}
\usepackage{subfigure}
\usepackage{psfrag}
\usepackage{color}
\usepackage{enumerate}
\usepackage[numbers,sort&compress]{natbib}
\usepackage{epstopdf}
\usepackage{appendix}
\newtheorem{thm}{Theorem}[section]
\newtheorem{conj}[thm]{Conjecture}
\newtheorem{cor}[thm]{Corollary}%[section]
\newtheorem{lem}[thm]{Lemma}%[section]
\theorembodyfont{\rmfamily}
\newtheorem{defn}[thm]{Definition}%[section]
\def\pf{\bigskip\noindent {\bf Proof.}~~}

\def\dfn#1{{\sl #1}}

\def\less{\backslash}

\def\pf{\bigskip\noindent {\emph{Proof.}}~~}

\def\qed{ \hfill $\square$}

\newtheorem{ob}[thm]{Observation}

\title{The saturation number of $C_6$}
\author{
\small Yongxin Lan$^1$, Yongtang Shi$^2$,  Yiqiao Wang$^3$ and Junxue Zhang$^2$\\
\small $^1$ School of Science, Hebei University of Technology, Tianjin 300401, P.R. China\\
\small $^2$ Center for Combinatorics and LPMC, Nankai University, Tianjin 300071, P.R. China\\
\small $^3$ School of Management, Beijing University of Chinese Medicine, Beijing 100029, P.R. China\\
\small Emails: yxlan@hebut.edu.cn, shi@nankai.edu.cn, yqwang@bucm.edu.cn, jxuezhang@163.com
}
\date{}
\begin{document}\maketitle
\begin{abstract}
A graph $G$ is called $C_k$-saturated if $G$ is $C_k$-free but $G+e$ is not for any $e\in E(\overline{G})$.
The saturation number of $C_k$, denoted $sat(n,C_k)$,
is the minimum number of edges in a $C_k$-saturated graph on $n$ vertices.
Finding the exact values of $sat(n,C_k)$ has been one of the most intriguing open problems in extremal graph theory.
In this paper, we study the saturation number of $C_6$. We prove that ${4n}/{3}-2 \le sat(n,C_6) \le {(4n+1)}/{3}$ for all $n\ge9$, which significantly improves the existing lower and upper bounds for $sat(n,C_6)$.\\

\noindent\textbf{Keywords:} saturation number; cycles; saturated graph\\
\end{abstract}

\section{Introduction}
\baselineskip 16pt

All graphs considered in this paper are finite and simple. Throughout the paper we use the terminology and notation of \cite{W2001}.
For a graph $G$, we use $e(G)$ to denote  the number of edges, $|G|$ the number of vertices, $\delta(G)$  the minimum degree and $\overline{G}$ the complement of $G$. For  $A,B\subseteq V(G)$, let $B\setminus A:=B-A$
%Let $\overline{G}$ denote the complement graph of $G$.
and $N_G(A)$ denote the subset  of $V(G)\less A$  in which each vertex is adjacent to some vertex of $A$ in $G$. If $A=\{v\}$, then we write as $B\setminus v$ and $N_G(v)$. The degree of vertex $v$ in $G$, denoted $d_G(v)$, is the size of $N_G(v)$. Let $N_G[v]=N_G(v)\cup\{v\}$. For a graph $G$ and $u,v\in V(G)$, let $d_G(u,v)$ denote the length of a shortest path between $u$ and $v$.
Without confusion, we abbreviate  as $N(A)$, $N(v)$, $N[v]$, $d(v)$ and $d(u,v)$, respectively. We use $k$-path (resp. $k$-cycle) to denote the path (resp. $k$-cycle) of length $k$.
For positive integer $k$, let $[k]:=\{1,2,\ldots,k\}$.
%We denote a path, a cycle, a star, and a complete graph with $n$ vertices by $P_n$, $C_n$, $S_n$, and $K_n$, respectively.

Given a family of graphs $\mathcal{F}$, a graph is $\mathcal{F}$-free if it contains no any member in $\mathcal{F}$ as a subgraph.
A graph $G$ is called \dfn{$\mathcal{F}$-saturated} if $G$ is $\mathcal{F}$-free but $G+e$ is not for any $e\in E(\overline{G})$.
%is a subgraph of $G$, but for any $e\in E(\overline{G})$, some member of $\mathcal{F}$ is a subgraph of $G+e$.
The \dfn{saturation number} of $\mathcal{F}$, denoted $sat(n,\mathcal{F})$,
is the minimum number of edges in an $\mathcal{F}$-saturated graph on $n$ vertices.
If $\mathcal{F}=\{F\}$, then we write as $sat(n,F)$.
%A graph $G$ is called $F$-saturated if $G$ does not contain $F$ as a subgraph (not
%necessarily induced) but the addition of any missing edge to $G$ creates a copy of $F$. The saturation number $sat(n,F)$
%is defined as the minimum number of edges in an $F$-saturated graph on $n$ vertices.
Erd\H{o}s, Hajnal and Moon~\cite{EHM1964} initiated the study of the saturation numbers of graphs and proved that $sat(n,K_{k+1})=(k-1)n-{{k}\choose{2}}$ with $K_{k-1}\vee\overline{K_{n-k+1}}$ as the unique extremal graph.
%Furthermore, they proved that equality holds
%only for the graph $K_{k-1}\vee \overline{K_{n-k+1}}$,
%where $\vee$ denotes the standard graph joining operation.
Later, K\'{a}szonyi and Tuza \cite{KT1986} extended this to the case when $K_k$ is replaced by any family $\mathcal{F}$ of graphs, showing that $sat(n,\mathcal{F})=O(n)$.
Since then, a large quantity of work in this area has been carried out in determining the saturation numbers of  complete multipartite graphs, cycles, trees, forests and hypergraphs, see e.g.~\cite{P1999,BFP2010,CFG2008,FJMTW12,GS2007,HL2020,MHHG21,O1972,T1989,C2014}.
Surveys on the saturation problem of graphs and hypergraphs could be found in~\cite{FFS2011,P2004}.
%We briefly survey some results in two main lines of research in this area, namely the saturation numbers for cycles and complete multipartite graphs.
% We  first mention some results for complete multipartite graphs.
%When all but at most one parts have size 1, Pikhurko \cite{P1999} and  Chen,  Faudree, and Gould \cite{CFG2008} independently
%determined the saturation number of complete multipartite graphs with sufficiently large order.
%When there are at least two parts of size at least 2, the exact value was only known for $K_{2,2}$, $K_{2,3}$ and $K_{3,3}$.
%The exact value for $K_{2,2}$ was first determined by Ollmann \cite{O1972}. Later on, a shorter proof was given by Tuza \cite{T1989}.
%For $K_{2,3}$, there have been several papers aiming to determine $sat(n,K_{2,3})$ over the years.
%This was finally achieved by Chen \cite{C2014} who confirmed a conjecture of Bohman,  Fonoberova, and Pikhurko  \cite{BFP2010}
%that $sat(n, K_{2,3})= 2n-3$ for all $n\geq 5$.
%The exact value for  $K_{3,3}$ was determined by Shenwei  Huang et al. \cite{HL2020} that $sat(n, K_{3,3}) = 3n-9 $ for $n \ge 9$.
%For the case where the graph has $r$ parts and all parts have size 2,
%Gould and Schmitt \cite{GS2007} conjectured that $sat(n, K_{2,\ldots,2})=\lceil((4r-5)n-4r^2+6r-1)/2\rceil$, and
%they proved the conjecture when the minimum degree of the $K_{2,\ldots,2}$-saturated graphs is $2r-3$.
It is worth noting that the proofs of almost all results are quite technically involved and require significant efforts.
%Other than these graphs, no result on the exact value is known.

Finding the exact values of $sat(n,C_k)$ has been one of the most intriguing open problems in extremal graph theory. Until now, the exact values of $sat(n,C_k)$ are known for a few $k$. The result of Erd\H{o}s, Hajnal and Moon \cite{EHM1964} pointed out the star $S_n$ is the unique extremal graph for $C_3$. Ollmann \cite{O1972} determined all extremal graphs for $C_4$. Later, Tuza \cite{T1989} gave a shorter proof for the exact value of $sat(n, C_4)$. Fisher, Fraughnaugh and Langley \cite{FFL1997} derived an upper bound for $sat(n, C_5)$ by constructing a class of $C_5$-saturated graphs. Chen finally \cite{C2009,C2011} confirmed the upper bound obtained in \cite{FFL1997} is the exact value for $sat(n, C_5)$ and also characterized all extremal graphs. For the case of $C_n$, Bondy \cite{B1972} first showed that $sat(n,C_n)\geq \left\lceil\frac{3n}{2}\right\rceil$ for all $n>6$, which was later improved to $sat(n,C_n)\geq \left\lfloor\frac{3n+1}{2}\right\rfloor$ for all $n\ge20$ in \cite{CCES1986, CE1983, CES1982, LJZY1997}.
We summarize some of these results as follows.

\begin{thm}
Let $n$  be a positive integer.
\begin{enumerate}[(a)]
\item {\bf \cite{O1972}} $sat(n,C_3)=n-1$ for any $n\ge3$.
\item {\bf \cite{O1972,T1989}} $sat(n,C_4)=\left\lfloor\frac{3n-5}{2}\right\rfloor$ for any $n\ge5$.
\item {\bf \cite{C2009,C2011,FFL1997}} $sat(n, C_5)= \left\lceil\frac{10(n-1)}{7}\right\rceil$ for any $n\geq 21$.
\end{enumerate}
\end{thm}

For the general case $k\ge6$, many researchers devoted to the study of the saturation numbers of $C_k$, see \cite{B1972,BCEPST1996,CES1982,CE1983,CCES1986,FK2011,LJZY1997,ZLS2015}. We present the previous best results as follows.
%The following  is a result of Barefoot, Clark, Entringer, Porter, Sz\'ekely and Tuza~\cite{BCEPST1996}.
%
%\begin{thm}[\cite{BCEPST1996}]\label{BCEPST1996}
% For $n\ge k\ge6$, $$\left(1+\frac{1}{2k+8}\right)n\leq sat(n,C_k)\leq c_kn+O(1),$$
% where $c_k=2-\frac{1}{2^{k/2}-2}$ when $k$ is even and $c_k=2-\frac{1}{3\cdot2^{k-3}/2-2}$ when $k$ is odd.
%\end{thm}
% who gave the problem of estimating $sat(n,C_k)$ for $3 \leq k \leq n$.
\begin{thm}[\cite{GLS2006,ZLS2015}]\label{ZLS2015}
 For any $n\ge9$, we have $\left\lceil\frac{7n}{6}\right\rceil-2\le sat(n,C_6)\le \left\lfloor\frac{3n-3}{2}\right\rfloor$.
\end{thm}

\begin{thm}[\cite{FK2011}]
For all $k\geq 7$ and $n\geq 2k-5$, we have
$\left(1+\frac 1 {k+2}\right)n-1\leq sat(n,C_k)\leq \left(1+\frac 1 {k-4}\right)n+{{k-4}\choose 2}.$
\end{thm}

Based on the constructions that yield the upper bounds, F\"{u}redi and Kim believed the upper bounds are approximately optimal and proposed the following conjecture, which is widely open.

\begin{conj}[\cite{FK2011}]\label{conj}
There exists a constant $k_0$ such that
$sat(n,C_k)= \left(1+\frac 1 {k-4}\right)n+O(k^2)$ for all $k\geq k_0$.
\end{conj}

%Finding exact values of $sat(n, C_k)$ is far from trivial. Until now, only for $k=3,4,5$ the exact values are known. Ollmann \cite{O1972} first investigated the star $S_n$ is unique extremal graph for $C_3$, and then determined all extremal graph for $C_4$. Later, Tuza \cite{T1989} gave a shorter proof for the exact value of $sat(n, C_4)$. An upper bound for $sat(n, C_5)$ was derived by Fisher, Fraughnaugh and Langley \cite{FFL1997}. Chen \cite{C2009,C2011} confirmed the above upper bound is also the lower bound and determined all extremal graphs. Notice that these extremal graphs for $C_4$ and $C_5$ are not unique. We summarize their results as follows.
%
%\begin{thm}
%\begin{enumerate}[(a)]
%\item {\bf \cite{O1972}} $sat(n,C_3)=n-1$ for any $n\ge3$.
%\item {\bf \cite{O1972,T1989}} $sat(n,C_4)=\lfloor\frac{3n-5}{2}\rfloor$ for any $n\ge5$.
%\item {\bf \cite{C2009,C2011,FFL1997}} $sat(n, C_5)= \lceil\frac{10(n-1)}{7}\rceil$ for any $n\geq 21$.
%\end{enumerate}
%\end{thm}

It seems quite difficult to determine the exact values of $sat(n,C_k)$ for any $k\ge6$ as mentioned
by Faudree, Faudree and  Schmitt \cite{FFS2011}.
%In \cite{BCEPST1996}, the authors also gave the better upper bound for $C_6$, that is, $sat(n,C_6)\le 3n/2$ for any $n\ge11$. Zhang, Luo and Shigeno \cite{ZLS2015} recently improved the previously bounds in \cite{BCEPST1996} and further proved the following result.
In this paper, we prove the following result, which improves Theorem~\ref{ZLS2015}.

\begin{thm}\label{boundC6}
Let $n\ge9$ be a positive integer and $n \equiv \varepsilon \ (\text{mod } 3)$ with $\varepsilon\in\{0,1,2\}$. Then
\begin{enumerate}[(a)]
\item $sat(n,C_6)\le\frac{4n}{3}+\frac{\varepsilon^2}{2}-\frac{5\varepsilon}{6}$;
\item $sat(n,C_6)\ge \frac{4n}{3}-2$.
\end{enumerate}
%	For $n \geq 9$,	$ \frac{4n}{3}-2 \le sat(n,C_6)\le\frac{4n-\varepsilon}{3}+{\varepsilon \choose 2}$, where $\varepsilon\equiv(n \text{ mod } 3)$.
\end{thm}
%We proof Theorem \ref{boundC6} by proofing Theorem \ref{upperC6}, Theorem \ref{d1lower} and Theorem \ref{d2lower}.
%
%\begin{thm}\label{upperC6}
%Then $sat(n,C_6)\le\frac{4n-\varepsilon}{3}+{\varepsilon \choose 2}$, where $\varepsilon\equiv(n \text{ mod } 3)$.
%\end{thm}

%\begin{thm}\label{upperC6}
%	For $k \ge 1$,
% $sat(n,C_6) \leq \begin{cases}
% \frac{4n}{3}     &\text{if $n=3k+6$}, \\
% \frac{4n-1}{3}   & \text{if $n=3k+7$},\\
% \frac{4n+1}{3}   & \text{if $n=3k+8$}.\\
%
% \end{cases}$\\
%\end{thm}
From Theorem~\ref{boundC6}, one can see $sat(n,C_6)=4n/3+O(1)$ for any $n\ge9$, which confirms that the constant $k_0$ in Conjecture~\ref{conj} should be at least $7$.

We organize our paper as follows. In Section 2, we prove Theorem~\ref{boundC6}($a$) by giving a new construction. In Section 3, we  prove Theorem~\ref{boundC6}($b$) by assuming  Theorem~\ref{main} holds. To complete the proof of Theorem~\ref{main}, we investigate that  we just need to prove Theorems~\ref{-320}, \ref{-3little} and \ref{f7vge0} in Section 4. Several lemmas are provided in Section 5 but their proofs are contained in Appendix.  We then prove Theorems~\ref{-320}, \ref{-3little} and \ref{f7vge0} in Sections 6-8.
We need to introduce more notation.  Given a graph $G$, let $\mathcal{F}_G$ denote the family of all graphs $G+e$ containing a $6$-cycle as a subgraph, where $e\in E(\overline{G})$, and $C_6(e)$ denote the family of $6$-cycles containing $e$ in $G+e$. By abusing notation, we also use $C_6(e)$ to denote one member in $C_6(e)$. Denote by $P(uv)$ the path with ends $u$ and $v$. Given a graph $G$ with vertex partition $V(G)=V_1\cup \cdots\cup V_s$ and $x\in V(G)$, let $N_i(x)=N(x)\cap V_i$ and $n_i(x)=|N_i(x)|$.\qed
%Let $\Theta_5$ denote the graph obtained from $5$-cycle by joining two nonadjacent vertices of $5$-cycle.

\section{Proof of Theorem~\ref{boundC6}($a$)}
\begin{figure}[htbp]
\centering
\includegraphics[scale=0.3]{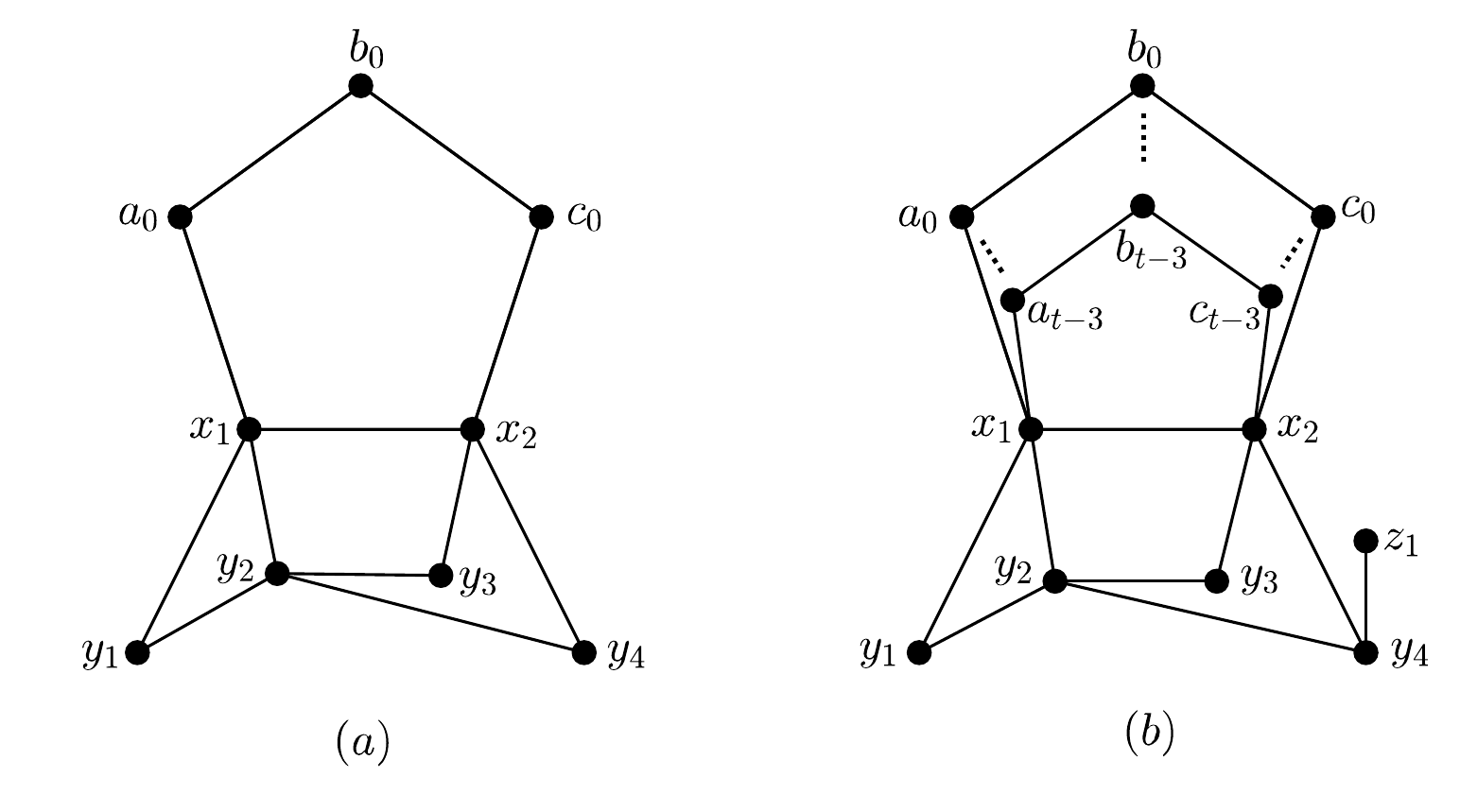}
\caption{($a$) ${H}$ and ($b$) $G_t^1$.}\label{fig GH}
\end{figure}

Let $n$ and $\varepsilon$ be given as in the statement. Let $H$ be the graph depicted in Figure~\ref{fig GH}($a$). Let $P^i$ be a path with vertices $a_i,b_i,c_i$ in order. For $n=3t$ with $t\ge3$, let $G_t^{0}$ be a graph on $n=3t$ vertices obtained from $H\cup P^1\cup \cdots\cup P^{t-3}$ by adding edges $x_1a_i$ and $x_2c_i$ for any $i\in[t-3]$.
Then $G_t^0$ is $C_6$-free and $$e(G_t^0)=e(H)+e(P^1)+\cdots+e(P^{t-3})+2(t-3)=12+2(t-3)+2(t-3)=4t=\frac{4n}{3}.$$
For $n=3t+\varepsilon$ with $t\ge3$ and $\varepsilon\in[2]$, let $G_t^{\varepsilon}$ be a graph on $n=3t+\varepsilon$ vertices obtained from $G_t^0\cup K_{\varepsilon}$ by joining $y_4$ to all vertices of $K_{\varepsilon}$. The graph $G_t^1$ is depicted in Figure~\ref{fig GH}($b$). Then for any $\varepsilon\in[2]$, $G_t^{\varepsilon}$ is $C_6$-free and $$e(G_t^{\varepsilon})=e(G_t^0)+e(K_\varepsilon)+\varepsilon=4t+\frac{\varepsilon(\varepsilon-1)}{2}
+\varepsilon=\frac{4n-\varepsilon}{3}+\frac{\varepsilon(\varepsilon-1)}{2}=\frac{4n}{3}+\frac{\varepsilon^2}{2}-\frac{5\varepsilon}{6}.$$

We next prove that  $G_t^{\varepsilon}$ is $C_6$-saturated for any $t\ge3$ and $\varepsilon\in\{0,1,2\}$. Observe that  if $G_t^1$ is $C_6$-saturated, then both $G_t^0$ and $G_t^2$ are $C_6$-saturated.  Thus, we shall prove that  $G_t^1$ is $C_6$-saturated for any $t\ge3$. It has been proved in \cite{ZLS2015} that $H$ is $C_6$-saturated. Thus we just need to show that
$G_t^1+z_1v\in \mathcal{F}_{G_t^1}$ for any $v\in V(H)\setminus y_4$, and $\{G_t^1+a_ia_j, G_t^1+a_ib_j,G_t^1+a_ic_j,G_t^1+b_ib_j,G_t^1+b_ic_j,G_t^1+c_ic_j\}\subseteq \mathcal{F}_{G_t^1}$  for any $0\le i< j\le t-3$. It is clear that the latter holds since $C_6(a_ia_j)=a_ib_ic_ix_2x_1a_j$,  $C_6(a_ib_j)=a_ib_ic_ix_2c_jb_j$, $C_6(a_ic_j)=a_ix_1y_2y_3x_2c_j$, $C_6(b_ib_j)=b_ia_ix_1x_2c_jb_j$, $C_6(b_ic_j)=b_ia_ix_1a_jb_jc_j$ and $C_6(c_ic_j)=c_ib_ia_ix_1x_2c_j$. For $v\in \{x_1,y_1,c_0,y_3\}$,  $C_6(z_1v)=z_1P^1v$ where $P^1=y_4x_2y_3y_2$ or $P^1=y_4y_2x_1x_2$. For $v\in \{x_2,a_0\}$, $C_6(z_1v)=z_1P^2v$, where $P^2=y_4y_2y_1x_1$. For $v\in\{y_2,b_0\}$, $C_6(z_1v)=z_1P^3v'v$, where $P^3=y_4x_2x_1$ and $v'\in N(x_1)\cap N(v)$. Hence, for any $t\ge3$ and $\varepsilon\in\{0,1,2\}$, $G_t^{\varepsilon}$ is $C_6$-saturated. Thus, $sat(n,C_6)\le e(G_t^{\varepsilon})\le 4n/3+\varepsilon^2/2-5\varepsilon/6$, as desired.\qed

%\begin{figure}
%	\begin{center} \subfigure[]{\includegraphics[width=7.5cm]{9.pdf}}
%		\subfigure[]{\includegraphics[width=7.5cm]{3k67.pdf}}
%		\caption{(a)${G_0}$ and (b)$G_t^1$}\label{fig GH}
%	\end{center}
%\end{figure}

\section{Proof of Theorem~\ref{boundC6}($b$)}
We call a vertex of degree two in a graph a {\bf root}. A root is {\bf good} if it does not lie in any triangle, otherwise it is {\bf bad}.
For any graph $G$, let $B(G)$ denote the set of all bad roots in $G$, $B_1(G)=\{v\in B(G)| d(v')\ge3 \text{ for every } v'\in N(v)\}$ and $B_2(G)=B(G)\less B_1(G)$.
%For any $v\in V(G)$ with $d(v)=2$, say $N(v)=\{v_1, v_2\}$, let  $B_1(G)=\{v|d(v)=2,~v_1v_2\in E(G),$ $d(v_1)\ge 2$, and $d(v_2)\ge 3
%\}$ and $B_2(G)=\{v|d(v)=2, ~v_1v_2\in E(G), ~d(v_1)=2$ or $d(v_2)=2$\}.
%vertices $v\in V(G)$ such that $d(v)=2$ and $v$ is in a triangle,
%$B_1(G)$ the set of vertices $v\in B(G)$ such that $v$ has no neighbor with degree two, and
%$B_2(G)$ denote the set of vertices $v\in B(G)$ such that $v$  has neighbor with degree two.
%Clearly, $B(G)=B_1(G)\cup B_2(G)$.
 %It follows from the following result that Theorem~\ref{boundC6}($b$) holds.

%\begin{lem}\label{path joint}
%Let $G$ be a $C_6$-saturated graph and $e\in E(\overline{G})$. Let $P=P(xy)$ be $p$-path with $e\in E(P)$ in $G+e$ and $Q=P(xy)$ be $(6-p)$-path in $G$.
%Then $Q$ and $P'$ have at least three common vertices for any $p$-path $P'=P(xy)$ of $G$.
%
%\vspace{-12pt}
%
%\end{lem}
%\pf The statement obviously holds because $Q$ is a subgraph of $G$ and $G$ is $C_6$-free.\qed

\begin{lem}[\cite{ZLS2015}]\label{2intriang}
	Let $G$ be a $C_6$-saturated graph with $\delta(G)=2$.
%Let $A=\{v \in V(G)|d(v)=2, v\text{  is in a triangle} \\ \text{ and has no neighbor with degree two}\}$.
Then at least one of the following is satisfied:\vspace{-8pt}

\begin{enumerate}[(1)]
\item There exists a root $\alpha \in B_1(G)$ with $N(\alpha)=\{\alpha_1, \alpha_2\}$ such that $N(\alpha_i)=\{\alpha,\alpha_{3-i},x_i\}$ for any $i\in[2]$ and $x_1\ne x_2$.\vspace{-8pt}
%We can choose  a root $\alpha \in B_1(G)$ such that $N(\alpha)=\{\alpha_1, \alpha_2\}$, $N(\alpha_1)=\{\alpha,\alpha_2,x_1\}$ and $N(\alpha_2)=\{\alpha,\alpha_1,x_2\}$, where $x_1\ne x_2$;\vspace{-8pt}

%$V_1=\{\alpha_1, \alpha_2, \alpha\}$ and $V_2=\{x_1, x_2\}$ satisfy $N(\alpha)=\{\alpha_1, \alpha_2\}$,  $N_2(\alpha_1)=\{x_1\}$ and $N_2(\alpha_2)=\{x_2\}$, where $x_1\neq x_2$;\vspace{-8pt}

\item $D(X) \ge 3|X|$ holds for $X=V(G) \setminus (\{v \in V(G)|d(v)=2\}\less B_1(G))$, where $D(X)=\sum_{v \in X}d(v)$.
\end{enumerate}
\end{lem}

In fact, we can derive a stronger result than Lemma~\ref{2intriang}.

\begin{lem}\label{new2intriang}
	Let $G$ be a $C_6$-saturated graph with $\delta(G)=2$. Then we must have $D(X) \ge 3|X|$ for $X=V(G) \setminus (\{v \in V(G)|d(v)=2\}\less B_1(G))$, where $D(X)=\sum_{v \in X}d(v)$.\vspace{-12pt}
\end{lem}

\pf Suppose not. By Lemma~\ref{2intriang}, we see Lemma~\ref{2intriang}(1) holds. Let $\alpha$ be such a root. Let $N(\alpha)$, $N(\alpha_1)$ and $N(\alpha_2)$ be defined as in the statement of the Lemma~\ref{2intriang}.
Then $G$ has a $5$-path $P(\alpha x_1)$ because $G+\alpha x_1\in\mathcal{F}_G$. Let $P(\alpha x_1)=\alpha a_1a_2a_3a_4x_1$. Note that $d(\alpha_i)=3$ for any $i\in [2]$. If $a_1=\alpha_2$, then $\alpha_1=a_i$ for some $i\in\{2,3,4\}$ because $G$ is $C_6$-free.
But then $d(\alpha_1)\ge4$, a contradiction. If $a_1=\alpha_1$, then $a_2=\alpha_2$ and $a_3=x_2$. But then $G$ has a copy of $C_6$ with vertices $x_1,\alpha_1,\alpha,\alpha_2,x_2,a_4$ in order, a contradiction.\qed\medskip

This immediately leads to the following theorem.

\begin{thm}\label{goodroot}
Let $n\ge6$ be a positive integer. Let $G$ be a $C_6$-saturated graph on $n$ vertices with $B_2(G)=\emptyset$. If $\delta(G)=2$ and $G$ has no good root, then $e(G)\ge4n/3-2$.\vspace{-12pt}
\end{thm}
\pf Since $G$ has no good root, we see $\{v \in V(G)|d(v)=2\}=B_1(G)\cup B_2(G)$. Hence, $\{v \in V(G)|d(v)=2\}\less B_1(G)=B_2(G)=\emptyset$.
By Lemma \ref{new2intriang}, $D(V(G))=\sum_{v\in V(G)}d(v)\ge3|V(G)|$. Then $e(G)=(\sum_{v\in V(G)}d(v))/2\ge3n/2\ge4n/3-2$.\qed\medskip

To complete the proof of Theorem~\ref{boundC6}($b$), we just need to prove the following Theorem.

\begin{thm}\label{main}
Let $n\ge6$ be a positive integer. Let $G$ be a $C_6$-saturated graph on $n$ vertices with $B_2(G)=\emptyset$. If $\delta(G)=1$, or $\delta(G)=2$ and $G$ has at least one bad root, then $e(G)\ge 4n/3-2$.
\end{thm}

\noindent{\bf Proof of Theorem~\ref{boundC6}($b$):}  Let $G$ be a $C_6$-saturated graph on $n\ge9$ vertices with $e(G)=sat(n,C_6)$. If $B_2(G)=\emptyset$, then by Theorems~\ref{goodroot} and \ref{main}, we have $e(G)\ge 4n/3-2$. So we may assume $B_2(G)\ne\emptyset$. Let $G_1$ be a graph obtained from $G$ by deleting every vertex of $B_2(G)$.  Clearly, $G_1$ is a $C_6$-saturated graph and $B_2(G_1)=\emptyset$. If $|G_1|\le5$, then $G_1$ is a complete graph and so $e(G_1)={|G_1|(|G_1|-1)}/{ 2}>4|G_1|/3-2$. If $|G_1|\ge6$, then  $e(G_1)\ge4|G_1|/3-2$ by Theorem~\ref{main}. It is easy to see that $e(G[B_2(G)])+e(B_2(G), V(G)\setminus B_2(G))=3|B_2(G)|/2$. Thus, $e(G)=e(G_1)+e(B_2(G), V(G)\setminus B_2(G))+e(G[B_2(G)])\ge 4|G_1|/3-2+3|B_2(G)|/2>4n/3-2$.\qed
\medskip

\section{Proof of Theorem~\ref{main}}
Let $G$ be a $C_6$-saturated graph on $n$ vertices. If $\delta(G)\ge3$, then $e(G)\ge3n/2>4n/3-2$ for any $n\ge6$.
So we assume that $\delta(G)\le2$.
Let $\alpha\in V(G)$ such that $d(\alpha)=\delta(G)$. Set $N[\alpha]=\{\alpha,\alpha_1,\ldots,\alpha_{\delta(G)}\}$.
Denote by $\Theta_5$ the graph obtained from a $5$-cycle by joining two non-adjacent vertices. Let $S=\{v\in V(G)|d(v)=2\}$, $\mathcal{C}_4=\{v\in S| v \text{ is contained in a 4-cycle}\}$, $\mathcal{C}_5=\{v\in S| v\text{ is contained in a 5-cycle}\}$, $S_5=\{v\in S| v\text{ is contained in a }\Theta_5\}$, $S_4=\{v\in S| v\text{ is not contained in a }\Theta_5\}$, $S_3=\{v| v\in \mathcal{C}_4-\mathcal{C}_5\}$,  $S_2=\{v| v\in \mathcal{C}_5-\mathcal{C}_4\}$ and $S_1=S-\mathcal{C}_4-\mathcal{C}_5$.
 We choose such $\alpha$ satisfying:  when $\delta(G)=1$, the number of $4$-cycles containing $\alpha_1$ is minimum; when $\delta(G)=2$, $\alpha$ is a good root and $\alpha \in S_i$  for some $i \in [5]$ such that $i$ is as small as possible.
%We choose such a vertex $\alpha$ satisfying:  when $\delta(G)=1$, the number of $4$-cycles containing $\alpha_1$ is minimum; when $\delta(G)=2$, $\alpha$ is a good root and \textcolor{blue}{the number of $4$-cycles and $\Theta_5$ containing $\alpha$ is minimum. 这里应该是Θ图的个数最小吧}
Let $V_1=N[\alpha]$ and $V_{i}=\{x\in V(G)| d(x,\alpha)=i\}$ for any $i\ge2$. Clearly, $V_i=\emptyset$ for any $i\ge6$ because $G$ is $C_6$-saturated, that is, any two non-adjacent vertices are the ends of some $5$-path in $G$. Thus, $V_1,\ldots,V_5$ form a partition of $V(G)$. We define a function as follows.

%Then $V(G)=V_1\cup \cdots \cup V_5$ because any two nonadjacent vertices are the ends of some $5$-path in $G$. We define a function as follows.
%For $x\in V(G)$, let $N_i(x)=N(x)\cap V_i$ and $n_i(x)=|N_i(x)|$.

\begin{defn}\label{gfunc}
For any $x\in V_i$ and $i\ge1$, let
\begin{equation*}
g(x)=
\begin{cases}
\frac{1}{2}n_i(x)-\frac{4}{3} &\;\text{ if } \; i=1;\\[2mm]
n_{i-1}(x)+\frac{1}{2}n_i(x)-\frac{4}{3} &\;\text{ if } \;i\ge2.
\end{cases}
\end{equation*}
\end{defn}

%For any $x\in V_i$ and $i\ge1$, we define a function as follows:
%\begin{equation*}
%g(x)=
%\begin{cases}
%\frac{1}{2}n_i(x)-\frac{4}{3} &\;\text{ if } \; i=1;\\
%n_{i-1}+\frac{1}{2}n_i(x)-\frac{4}{3} &\;\text{ if } \;i\ge2.
%\end{cases}
%\end{equation*}

\noindent{\bf Remark 1.} For any $x \in V_i$ with $i\ge2$,  $g(x)\ge n_{i-1}(x)/3+n_i(x)/6\ge\frac{2}{3}$ if $n_{i-1}(x)\ge 2$ or $n_i(x)\ge2$; $g(x)=\frac{1}{6}$ if $n_{i-1}(x)=1$ and $n_i(x)=1$; and $g(x)=-\frac{1}{3}$ if $n_{i-1}(x)=1$ and $n_i(x)=0$.\medskip

 We define $N_{i}^{\star}(x):=N_{i}(x)\cap V_i^{\star}$ and $n_i^{\star}(x):=|N_i^{\star}(x)|$, where $\star\in\{-,+,1,2,-1,-2\}$. By Remark 1, for $i\ge2$, $V_i$ can be partitioned into $V_i^+,V_i^-$ such that
 $$V_i^{+}:=\{x\in V_i | g(x)\ge0\} \text{ and } V_i^-=V_i\less V_i^+=\{x\in V_i | g(x)=-1/3\}.$$
Moreover, we partition $V_i^+$ into $V_i^1,V_i^2$ and $V_i^-$ into $V_i^{-1},V_i^{-2}$ such that $$V_i^{2}:=\{x\in V_i^+ | g(x)\ge2/3\} \text{ and }  V_i^{1}:=V_i^+\less V_i^2=\{x\in V_i^+ | g(x)=1/6\};$$
$$V_i^{-1}:=\{x\in V_i^-|n_{i+1}^2(x)\ge2\} \text{ and } V_i^{-2}:=V_i^{-}\less V_i^{-1}.$$\vspace{-12pt}
%$V_i^{+}:=\{x\in V_i | g(x)\ge0\} \text{ and } V_i^-=V_i\less V_i^+=\{x\in V_i | g(x)=-1/3\}.$
%Moreover, we partition $V_i^+$ into $V_i^1,V_i^2$ such that $V_i^{2}:=\{x\in V_i^+ | g(x)\ge2/3\} \text{ and }  V_i^{1}:=V_i^+\less V_i^2=\{x\in V_i^+ | g(x)=1/6\}$;
%$V_i^-$ into $V_i^{-1},V_i^{-2}$ such that
%$V_i^{-1}:=\{x\in V_i^-|n_{i+1}^2(x)\ge2\} \text{ and } V_i^{-2}:=V_i^{-}\less V_i^{-1}$.

%Clearly, $\sum_{x\in V_1}g(x)=e(G[V_1])-\frac{4|V_1|}{3}\ge \min\{1-\frac{8}{3},2-4\}=-2$.

\begin{figure}[htbp]
	\centering \includegraphics[scale=0.3]{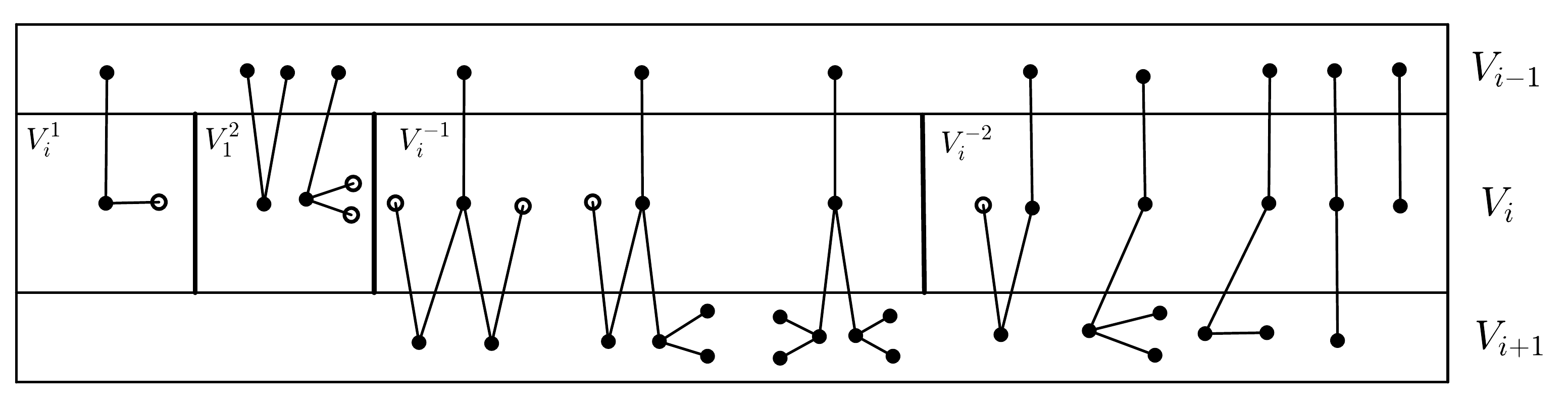}\vspace{-5mm}
	\caption{The partition of $V_i$.}
	\label{6}
\end{figure}
For any $i\ge2$, the partition of $V_i$ is depicted in Figure~\ref{6}, where the hollow vertices do not necessarily belong to the set shown in the Figure. By Definition \ref{gfunc},
\begin{align}
e(G)=&\sum_{i=1}^5e(G[V_i])+\sum_{i=2}^5e(V_i,V_{i-1})
= \sum_{i=1}^5\frac{1}{2}\sum_{x\in V_i}n_i(x) +\sum_{i=2}^5\sum_{x\in V_i}n_{i-1}(x)\notag\\
=&\sum_{i=1}^5\sum_{x\in V_i}\left(g(x)+\frac{4}{3}\right)=\sum_{x\in V_1}g(x)+\sum_{x\in V(G)\setminus V_1}g(x)+\frac{4}{3}n\notag\\
=&{\sum_{x\in V(G)\setminus V_1}g(x)+\frac{4}{3}n+e(G[V_1])-\frac{4}{3}|V_1|}\notag\\
=&{\sum_{x\in V(G)\setminus V_1}g(x)+\frac{4}{3}n+(|V_1|-1)-\frac{4}{3}|V_1|}\notag\\
=&{\sum_{x\in V(G)\setminus V_1}g(x)+\frac{4}{3}n-\frac{1}{3}|V_1|-1}\notag\\
\ge& \sum_{x\in V(G)\setminus V_1}g(x)+\frac{4}{3}n-2,
\end{align}
where $e(V_i,V_{i-1})$ is the number of edges with one end in $V_i$ and the other end in $V_{i-1}$.

To obtain the desired lower bound in Theorem~\ref{main}, we just need to show that $\sum_{x\in V(G)\setminus V_1}g(x)\ge0$.
%By inequalities~(2) and (4), it is enough to prove Theorem~\ref{f7vge0}. Thus, in the rest of the paper, we will prove Theorem~\ref{f7vge0} by a series of lemmas.
{The basic idea of our proof is as follows: we first allocate a charge to each vertex in $G$ such that the initial charge function is the function $g$ defined in Definition~\ref{gfunc}, and then reallocate the charge to each vertex in $G$ such that the final charge at each vertex of $V(G)\setminus V_1$ is at least zero. To reallocate the charge to each vertex in $G$, we need to define a series of charge functions, which can be regarded as two stages: at the first stage, the charge function $g_i$ for any $i\in[5]$ is defined;   at the second stage, the  charge function $f_i$ for any $i\in [7]$ is defined, such that  $\sum_{x\in V(G)}g_{i-1}(x)=\sum_{x\in V(G)}g_{i}(x)$ for any $i\in[5]$ and $\sum_{x\in V(G)}f_{j-1}(x)=\sum_{x\in V(G)}f_{j}(x)$ for any $j\in[7]$, where $g_0=g$ and $f_0=g_5$. Since the charges of vertices of $V_1$ are never altered, it suffices to prove that $\sum_{x\in V(G)\setminus V_1}f_7(x)\ge0$.}
%The basic idea of our proof is similar to the ``discharging method". The function $g$ is the assignment to each vertex $x$ such that an ``initial charge" of $x$ is equal to $g(x)$. We reallocate the ``charge" of each vertex such that the final charge at each vertex except the vertices of $V_1$ is at least zero.  For this purpose, a series of functions (like ``discharging rules") are defined, and in each function, the ``charge" of each vertex would be updated although the ``charges" of most vertices for which we do not involve in each function keep unchanged. Also, the definitions of these functions can be regarded as two stages: at the first stage, the ``initial charge" of $x$ is equal to $g(x)$ and the ``final charge" of $x$ is equal to $g^*(x)$; at the second stage, the ``initial charge" of $x$ is equal to $g^*(x)$ and the ``final charge" of $x$ is equal to $f_7(x)$.

\medskip
In the following, we define the {charge functions} at the first stage. {For any vertex $x$ uninvolved in the definition of $g_i$, $g_i(x)$ is still equal to $g_{i-1}(x)$}. In the following, for $x\in V(G)$, we call $x$ as $j$-vertex if $g(x)=j$; $j^+$-vertex if $g(x)\ge j$.
\begin{defn}\label{law1}

\begin{enumerate}[(1)]

\item Several $\frac{2}{3}^+$-vertices send charges to some $-\frac{1}{3}$-vertices in the previous level and $\frac{1}{6}$-vertices in the same level; several $\frac{1}{6}$-vertices send charges to $-\frac{1}{3}$-vertices in the previous level. We define the concrete laws as follows:

(1.1) For any $u \in V_2^2$, let $g_1(u)=g(u)-\frac{1}{6}n_2^1(u)$, as depicted in Figure~\ref{4.2-1}($a_1$).
For any $u\in V_i^2$ with $i\in\{3,4\}$, let $g_1(u)=g(u)-\frac{1}{6}|A_1(u)|-\frac{1}{3}|B_1(u)|-\frac{1}{6}n_i^1(u)$, as depicted in Figure~\ref{4.2-1}($a_2$). For any $u\in V_i^2$ with $i=5$, let $g_1(u)=g(u)-\frac{1}{6}|A_1(u)|-\frac{1}{3}|B_1(u)|-\frac{1}{6}(n_5^1(u)-|C_1(u)|)$, as depicted in Figure~\ref{4.2-1}($a_3$), where $A_1(u)=\{v\in N_{i-1}^-(u)| n_{i}^2(v)\ge2\}$, $B_1(u)=\{v\in N_{i-1}^-(u)| n_{i}^2(v)=1\}$ and
\vspace{-8pt}

\[C_1(u)=
\begin{cases}
\, \{w\in N_5^1(u)| N_4(w)\subseteq N_4(w') \text{ for some } w'\in N_5(u)\setminus w\} & \text{ if } \,  d(u)=3, n_4(u)=n_4^-(u)=1;\\[2mm]

\, \emptyset & \text{ else}.
\end{cases}
\]
\vspace{-4pt}
% $C_1(u)=\{w\in N_5^1(u)| d(u)=3, n_4(u)=n_4^-(u)=1 \text{ and } N_4(w)\subseteq N_4(w'), \text{ where } w'\in N_5(u)\setminus w\}$.

(1.2) For any $w\in V_2^1$, let $g_1(w)=g(w)+\frac{1}{6}n_{i}^2(w)$, as depicted in Figure~\ref{4.2-1}($b_1$). For any $w\in V_i^1$ with $i\in\{3,4\}$, let $g_1(w)=g(w)+\frac{1}{6}n_{i}^2(w)-\frac{1}{6}|A_2(w)|$, as depicted in Figure~\ref{4.2-1}($b_2$), where $A_2(w)=\{v\in N_{i-1}^-(w)| n_i^2(v)=0\}$. For any $w\in V_5^1$, let $g_1(w)=g(w)+\frac{1}{6}n_{5}^2(w)$ if $w\notin C_1(u)$ for any $u\in V_5^2$, as depicted in Figure~\ref{4.2-1}($b_3$), $g_1(w)=g(w)$ otherwise.

%For any $w\in V_i^1$ with $i\in\{2,3,4\}$, then let $g_1(w)=g(w)+\frac{1}{6}n_{i}^2(w)-\frac{1}{6}|A_2|$, where $A_2=\{v\in N_{i-1}^-(w): n_i^2(v)=0\}$. For any $w\in V_5^1$, then let $g_1(w)=g(w)+\frac{1}{6}n_{5}^2(w)$ if $w\notin C_1$, $g_1(w)=g(w)$ otherwise.

(1.3) For any  $v\in V_i^-$ with $i\in\{2,3,4\}$, let $g_1(v)=g(v)+\frac{1}{6}n_{i+1}^2(v)$ if  $n_{i+1}^2(v)\ge 2$, as depicted in Figure~\ref{4.2-1}($c_1$); $g_1(v)=g(v)+\frac{1}{3}$ if $n_{i+1}^2(v)=1$, as depicted in Figure~\ref{4.2-1}($c_2$); $g_1(v)=g(v)+\frac{1}{6}n_{i+1}^1(v)$ when $i\in\{2,3\}$ (as depicted in Figure~\ref{4.2-1}($c_3$)) and $g_1(v)=g(v)$ when $i=4$ if $n_{i+1}^2(v)= 0$.

%In order to understand the definition, we give the figure as follows.
\begin{figure}[htbp]
	\centering
	\includegraphics[scale=0.55]{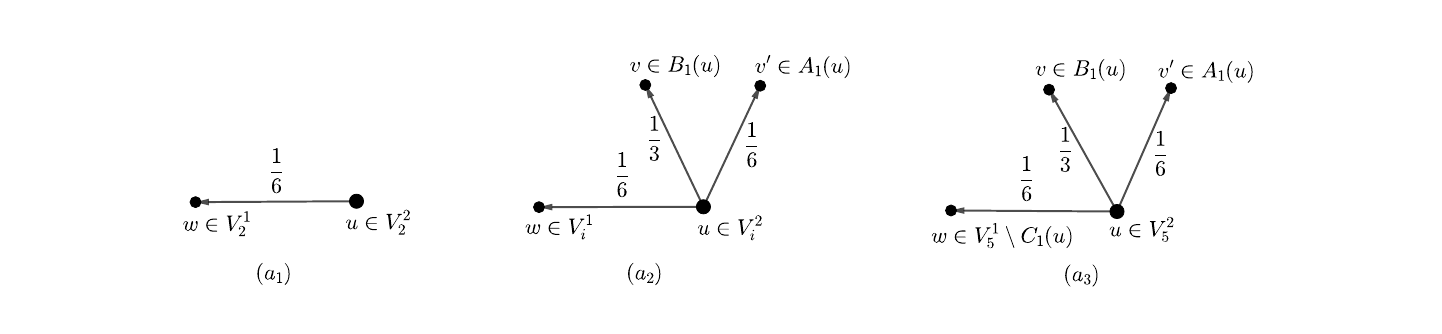}
\includegraphics[scale=0.5]{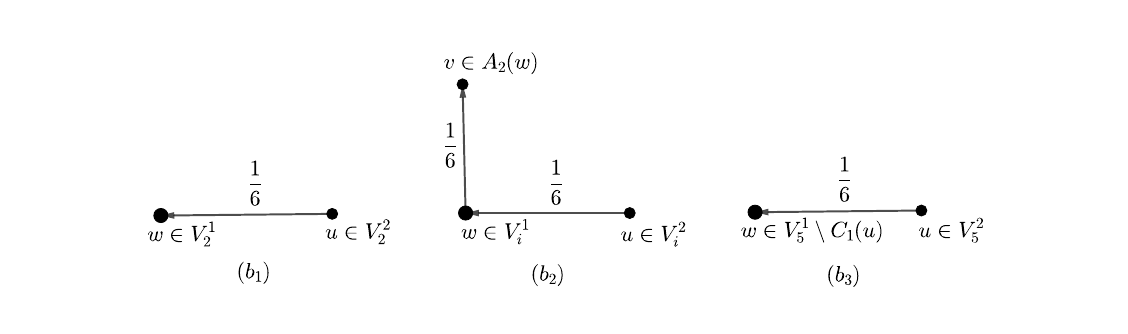}
\includegraphics[scale=0.5]{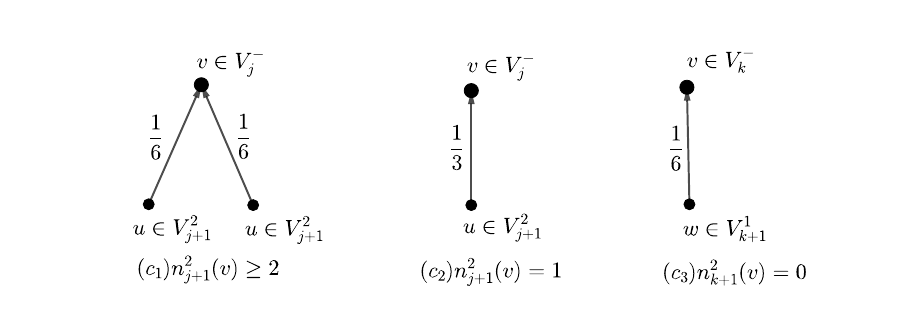}\vspace{-7pt}
	\caption{($a_1,a_2,a_3$), ($b_1,b_2,b_3$), ($c_1,c_2,c_3$) are the explanations of Definition \ref{law1}(1.1), (1.2), (1.3) respectively, where $i\in\{3,4\}$, $j\in\{2,3,4\}$ and $k\in\{2,3\}$.}
	\label{4.2-1}
\end{figure}

\item Several $\frac{1}{6}$-vertices of $V_5$ send charges to some $\frac{1}{6}$-vertices in $V_5$. Concretely, for any two adjacent vertices $w_1,w_2\in V_5^1$ with $N_4(w_i)=\{y_i\}$, let $g_2(w_i)=g_1(w_i)-\frac{1}{6}$ and $g_2(w_{3-i})=g_1(w_{3-i})+\frac{1}{6}$ if $g_1(y_i)\ge0$ and $g_1(y_{3-i})<0$,  as depicted in Figure \ref{4.2-2}($a$); $g_2(w_i)=g_1(w_i)$ otherwise.

\item Several $\frac{1}{6}$-vertices send charges to their neighbours in the previous level. Concretely, for any  $y\in V_i^-$ with $g_2(y)<0$, where $i\in \{2,3,4\}$, let $g_3(y)=g_2(y)+t|A_3(y)|$ and $g_3(z)=g_2(z)-g_2(z)=0$ for $z\in A_3(y)$, where $A_3(y)=\{z\in N_{i+1}^1(y)|g_2(z)=t\}$, and $t=\frac{1}{6}$ when $i\in\{2,3\}$ and $t=\frac{1}{3}$ when $i=4$, as depicted in Figure \ref{4.2-2}($b$).
\vspace{-2mm}

\item Several $\frac{1}{6}$-vertices of $V_5$ send charges to some $\frac{1}{6}$-vertices in $V_5$. Concretely, for any two adjacent vertices $w_1,w_2\in V_5^1$ with $N_4(w_i)=\{y_i\}$ and $g_3(w_1)=g_3(w_2)=\frac{1}{6}$, let $g_4(w_i)=g_3(w_i)-\frac{1}{6}$ and $g_4(w_{3-i})=g_3(w_{3-i})+\frac{1}{6}$ if $g_3(y_i)\ge0$ and $g_3(y_{3-i})<0$, as depicted in Figure \ref{4.2-2}($c$); $g_4(w_i)=g_3(w_i)$ otherwise.

%    For any $z\in V_5^1$ with  $N_4(z)=\{y\}$, $g_3(y)<0$ and $N_5(z)=\{z_1\}$, if $z_1\in N_5^1(z)$ and $g_3(y_1) \ge 0$, where $\{y_1\}=N_4(z_1)$, then $g_4(z)=g_3(z)+\frac{1}{6}$ and $g_4(z_1)=g_3(z_1)-\frac{1}{6}$.

\item Several $\frac{1}{6}$-vertices of $V_5$ send charges to their neighbours in the previous level. Concretely, for any  $y\in V_4^-$ with $g_4(y)<0$, let  $g_5(y)=g_4(y)+\sum_{z\in N_5^1(y)}g_4(z)$ and $g_5(z)=g_4(z)-g_4(z)=0$ for any $z \in N_5^1(y)$, as depicted in Figure \ref{4.2-2}($d$).

\begin{figure}[htbp]
\centering
\includegraphics[scale=0.7]{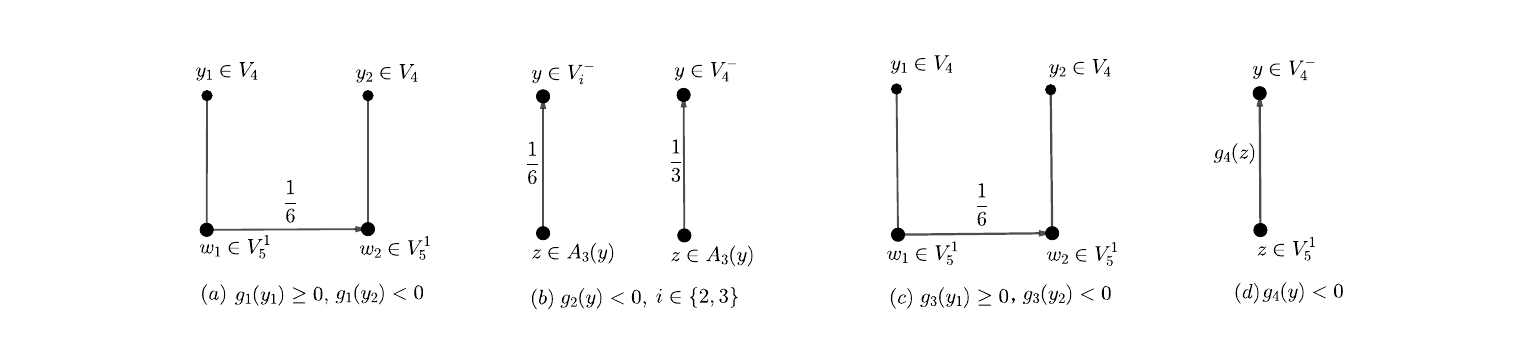}\vspace{-7mm}
\caption{($a$), ($b$), $(c)$ and $(d)$ are the explanations of  Definition \ref{law1}(2), (3), (4) and (5), respectively. }
	\label{4.2-2}
\end{figure}
    \vspace{-2mm}

\end{enumerate}
\end{defn}
	
\noindent For any $x\in V(G)$, let $g^*(x)=g_5(x)$ and $g_0(x)=g(x)$.
Clearly,

\vspace{-12pt}
\begin{align}\label{g}
\sum_{x \in V(G)\setminus V_1} g^*(x)=\sum_{x \in V(G)\setminus V_1}g(x).
\end{align}

\noindent Moreover, for any $0\le i< j \le 5$, if $g_i(x)\ge0$, then $ g_{j}(x)\ge0$; and if $g_i(x)<0$, then $g_0(x)\le \cdots \le g_i(x)<0$.
%($\ast$) if $g_i(x)\ge0$, then $ g_{j}(x)\ge0$; and if $g_i(x)<0$, then $g_0(x)\le \cdots \le g_i(x)<0$.\medskip
For any $i\in [5]$, let $V^{*-}_i:=\{x \in V_i|g^*(x)<0\}$, $N_{i}^{*-}(x):=N_{i}(x)\cap V_i^{*-}$ and $n_i^{*-}(x):=|N_i^{*-}(x)|$.
By Definitions~\ref{gfunc} and \ref{law1}, we see that for any $x\in V_i^2$ and $i\in\{3,4,5\}$,
\begin{align}\label{g*}
g^*(x)\ge&~ g(x)-\frac{1}{3}n_{i-1}^{-2}(x)-\frac{1}{6}n_{i-1}^{-1}(x) -\frac{1}{6}n_{i}^1(x)\notag\\
\ge&~ g(x)-\frac{1}{3}(n_{i-1}(x)-n_{i-1}^+(x)) +\frac{1}{6}n_{i-1}^{-1}(x) -\frac{1}{6}(n_{i}(x)-n_{i}^2(x))\notag\\
=& ~ \frac{2}{3}n_{i-1}(x)+\frac{1}{3}n_{i-1}^+(x)+\frac{1}{6}n_{i-1}^{-1}(x)+ \frac{1}{3}n_{i}(x)+\frac{1}{6}n_{i}^2(x)-\frac{4}{3}.
\end{align}

By Definition~\ref{law1} and Ineq.~(\ref{g*}), we have the following observations.\vspace{-8pt}
%By ($\ast$), we see  $N_5^{*-}(x)=N_5^{-}(x)$.
%From Definition~\ref{law1} and ($\ast$), we have the following Observation.

\begin{ob}\label{ob}
Let $i\in\{2,3,4\}$. For $x\in V_i^{*-}$, we have\\ (1)
 $n_{i+1}^2(x)=0$ and $n_{i+1}^1(x)\le1$; (2)
if $y\in N_{i+1}^1(x)$, then $n_{i+1}^1(y)=1$.\vspace{-8pt}
\end{ob}

\begin{ob}\label{ob2}
Let $i\in\{3,4,5\}$. For any $x\in V_i^2$, $g^*(x)\ge\frac{1}{3}n_{i-1}(x)+\frac{1}{6}n_{i}(x)$ if one of the following holds.\vspace{-8pt}
\begin{enumerate}[(1)]
\item $n_{i-1}^{+}(x)\ge2$ or  $n_{i-1}(x)\ge3$ and $n_{i-1}(x)+n_i(x)\ge5$;\vspace{-8pt}
\item $n_{i-1}^{+}(x)+n_{i-1}(x)\ge4$;\vspace{-8pt}
\item $n_{i-1}^+(x)+n_{i-1}(x)=3$ and $n_i^2(x)\ge1$;\vspace{-8pt}
\item $n_{i-1}^+(x)+n_{i-1}(x)=3$ and $n_{i-1}^{-1}(x)+n_i(x)\ge2$.\vspace{-8pt}
%\item $n_{i-1}^+(x)=1$ and $n_{i-1}(x)=2$ and $n_i(x)\ge2$;\vspace{-8pt}
%\item $n_{i-1}^+(x)=n_{i-1}(x)=1$ and $n_{i}^2(x)\ge2$.\vspace{-8pt}
%\item  $n_{i-1}^{+}(x)=0$ and $n_{i-1}\ge4$;
%\item  $n_{i-1}^{+}(x)=0$ and $n_{i-1}\ge3$ and $n_{i-1}^{-1}(x)\ge2$;
%\item  $n_{i-1}^{+}(x)=0$ and $n_{i-1}\ge3$ and $n_{i-1}^{-1}(x)=1$ and $n_i^2(x)\ge1$;
\end{enumerate}
\end{ob}

\begin{ob}\label{ob3}
Let $i\in\{3,4,5\}$. For any $x\in V_i^2$, $g^*(x)\ge\frac{1}{6}n_{i-1}(x)+\frac{1}{6}n_{i}(x)$ if one of the following holds.\vspace{-8pt}
\begin{enumerate}[(1)]
\item $n_{i-1}^{+}(x)\ge2$ or  $n_{i-1}(x)\ge2$ and $n_{i-1}(x)+n_i(x)\ge4$;\vspace{-8pt}
\item $n_{i-1}^{+}(x)+n_{i-1}(x)\ge3$;\vspace{-8pt}
%\item $n_{i-1}^+(x)=n_{i-1}(x)=1$ and $n_i^2(x)\ge1$;\vspace{-8pt}
%\item  $n_{i-1}^{+}(x)=0$ and $n_{i-1}\ge4$;
\item  $n_{i-1}^{+}(x)+n_{i-1}=2$ and $n_i^2(x)\ge1$;\vspace{-8pt}
\item  $n_{i-1}^{+}(x)+n_{i-1}=2$ and $n_{i-1}^{-1}(x)+n_i(x)\ge3$;\vspace{-8pt}
   % \textcolor{blue}{(4) $ n_{i-1}^+(x)=n_{i-1}(x)=1$ and  $n_i(x)\ge3$;}
\item  $n_{i-1}^{+}(x)+n_{i-1}=1$ and $n_{i-1}^{-1}(x)+n_i^2(x)\ge3$.
   % \textcolor{blue}{(5) $ n_{i-1}^+(x)=0$ and  $n_{i-1}(x)+n_{i-1}^{-1}(x)+n_i^2(x)\ge4$.}
\end{enumerate}
\end{ob}

 %Let $N_j^{*-}(x)=\{v\in N_j^-(x) : g^*(v)<0\}$ for $x\in V_i$ and $i,j\in [5]$.
% In the following, we let each vertex of $G$  change the function accoring to the Definition \ref{law1}. We know that $\sum_{x \in V} g^*(x)=\sum_{x \in V}g(x)$.
%We define $V^-_i=\{x \in V_i:g^*(x)<0\}$ and $P(xy)$ as the path between $x$ and $y$.
%Clearly, if $g(x) \ge 0$, then $g^*(x) \ge 0$.
%Let $V^{*-}_i=\{x \in V_i:g^*(x)<0\}$ for each $i\in [5]$. It is easy to see that for any $x\in V(G)$, if $g(x) >0$, then $g^*(x) \ge 0$.

Now we are ready to define the functions at the second stage. Here, the ``initial charge" of $x$ can be regarded as  $g^*(x)$.  For any vertex uninvolved in the following definition of $f_i$, the value of $f_i$ remains unchanged. For any $x\in V(G)$, let $N_i^{j-}(x):=\{v\in N_i(x)|f_j(v)<0\}$ and $N_i^{j+}(x)=N_i(x)\setminus N_i^{j-}(x)$. We define $n_i^{j-}(x)=|N_i^{j-}(x)|$ and $n_i^{j+}(x)=|N_i^{j+}(x)|$. Let $L(v)\subseteq N(v)$ denote the set of vertices  with degree one.
\begin{defn} \label{law2}
\begin{enumerate}[(1)]	
\item  Every vertex of $V_5$ sends charge equally to its neighbours in $V_4$. %as depicted in Figure \ref{4.6.1}.
 Concretely, for any $z\in V_4$ and $w\in V_5$, $f_1(z)=g^*(z)+\sum_{v\in N_5(z)}{g^*(v)}/{n_4(v)}$ and $f_1(w)=g^*(w)-\sum_{u\in N_4(w)}{g^*(w)}/{n_4(w)}=0$, as depicted in Figure \ref{4.6.1}($a$). With the following two exceptional cases:

\subitem ($1.1$) Suppose $w\in V_5$ with $N(w)=\{z_1,z_2,w_1\}$ such that $n_5^-(z_1)\ne0$, $n_5^-(z_2)=0$ and $N_5^1(w_1)\cap N_5^1(z_2)\ne \emptyset$, where $z_1,z_2\in V_4^-$ and $w_1\in V_5$.
Then $w$ will send  ${2g^*(w)}/{3}$ charge to $z_1$, and  ${g^*(w)}/{3}$ charge to $z_2$, as depicted in Figure \ref{4.6.1}($b$).
%(Hint: We have $n_5^-(z_2 )=0$. Suppose $x\in N_5^-(z_2)$, $y\in N_5^1(w_1)\cap N_5^1(z_2)$. Then $d(x)=1$, $d(y)=2$ and $G+xw_1\in \mathcal{F}_G$. By Proposition \ref{path joint}, there is a (4-p)-path $P(z_2x)$ or $P(w_1x)$ for $p\in [2]$, which contradicts $d(y)=2$. )

% Then let $f_1(z_1)=g^*(z_1)+\sum_{w \in N_5(z_1)\setminus{w_1}}{g^*(w)}/{n_4(w)}+{2g^*(w_1)}/{3}$, $f_1(z_2)=g^*(z_2)+\sum_{w \in N_5(z_2)\setminus w_1} {g^*(w)}/{n_4(w)}+{g^*(w_1)}/{3}\ge0$, and $f_1(w_1)=g^*(w_1)-{2g^*(w_1)}/{3}-{g^*(w_1)}/{3}=0$.
			 %$f_1(w)=0$ for any $w \in N_5(z_1)\cup N_5(z_2)$.}

\subitem ($1.2$) Suppose $w\in V_5$ with $N(w)=N_4^-(w)=\{z_1,z_2,z_3\}$ such that $n_5^-(z_1)\ne0$, $n_5^-(z_2)=n_5^-(z_2)=0$ and $N_3(z_2)=N_3(z_3)$. Then $w$ will send  ${g^*(w)}/{2}$ charge to $z_1$, and  ${g^*(w)}/{4}$ charge to $z_2$, and ${g^*(w)}/{4}$ charge to $z_3$, as depicted in Figure \ref{4.6.1}($c$).\vspace{-20pt}

%Let $P$ be a path with vertices $w_1z_1w_2z_2yz_3$ such that $y\in V_3$, $z_i\in V_4$ and $w_i\in V_5$. Let $H=P+w_2z_3$. If $g(z_i)<0$, $d(w_1)=1$, $d(w_2)=3$ and $n_5^-(z_2)=n_5^-(z_3)=0$, then let $f_1(z_1)=g^*(z_1)+\sum_{w \in N_5(z_1)\setminus{w_2}}\frac{g^*(w)}{n_4(w)}+\frac{g^*(w_2)}{2}$, $f_1(z_j)=g^*(z_j)+\sum_{w \in N_5(z_2)\setminus{w_2}}\frac{g^*(w)}{n_4(w)}+\frac{g^*(w_2)}{4}\ge \frac{1}{6}$ where $j \in \{2,3\}$, and $f_1(w_2)=g^*(w_2)-\frac{g^*(w_2)}{2}-\frac{g^*(w_2)}{4}\times2=0$.
%$f_1(w)=0$ \textcolor{blue}{???} for any $w \in N_5(z_1)\cup N_5(z_2)$.}
			
	\begin{figure}[htbp]
		\centering
		\includegraphics[scale=0.7]{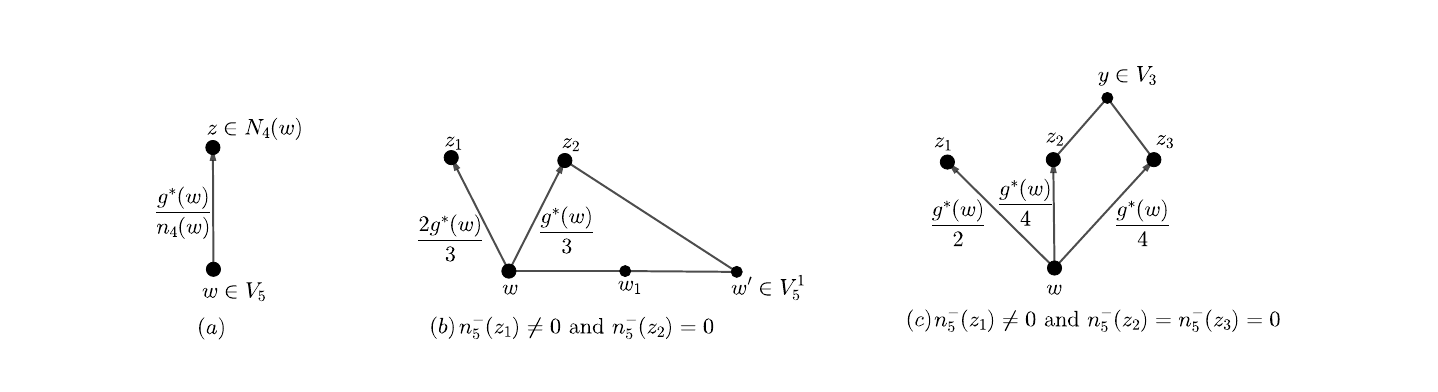}\vspace{-7mm}
		\caption{$(a)$, $(b)$ and $(c)$ are the explanations of  Definition 4.6 (1), (1.1) and (1.2), respectively. }
		\label{4.6.1}
	\end{figure}

	\item Several $\frac{1}{6}^+$-vertices of $V_4$ send charges to their neighbours in $V_4$, as depicted in Figure~\ref{4.6.2}. Concretely, for any $z\in V_4$,
if $f_1(z)\ge \frac{1}{6}n_4^{1-}(z)$, then $f_2(z)=f_1(z)-\frac{1}{6}n_4^{1-}(z)$, and $f_2(z')=f_1(z')+\frac{1}{6}$ for any $z'\in N_4^{1-}(z)$. Moreover,  for any two adjacent vertices $z_1,z_2\in V_4^1$ with $N_3(z_i)=\{y_i\}$ and $f_1(z_j)=0$ and $f_1(y_j)<0$, if $f_1(z_{3-j})\ge\frac{1}{3}$ and $f_1(y_{3-j})<0$, or $f_1(z_{3-j})\ge\frac{1}{6}$ and $f_1(y_{3-j})\ge0$, then let $f_2(z_{3-j})=f_1(z_{3-j})-\frac{1}{6}$ and $f_2(z_j)=f_1(z_j)+\frac{1}{6}$.\vspace{-5mm}

\begin{figure}[htbp]
	\centering
	\includegraphics[scale=0.65]{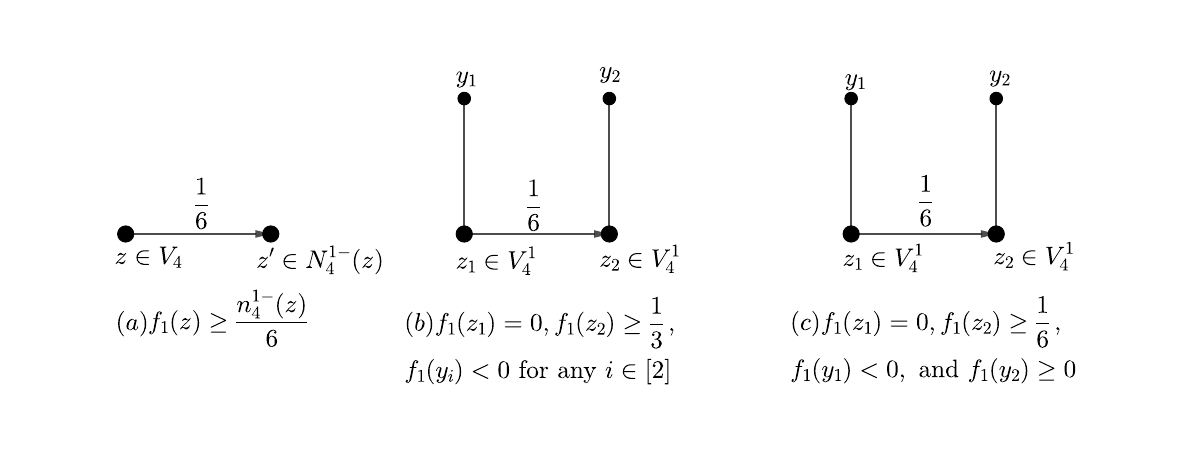}\vspace{-7mm}
	\caption{($a)$, $(b)$ and ($c$) are the explanations of Definition 4.6(2). }
	\label{4.6.2}
\end{figure}

\item Several vertices of $V_4$ send charges equally to their neighbours in $V_3$. Concretely, for any $y\in V_3$, then $f_3(y)=f_2(y)+\sum_{z\in N_4^{2+}(y)}{f_2(z)}/{n_3(z)}$ and $f_3(z)=f_2(z)-f_2(z)=0$ for any $z\in N_4^{2+}(y)$,  as depicted in Figure \ref{4.6.3}($a$). With the exception of  the following case:

	\subitem $(3.1)$ Suppose $z\in V_4$ with $N_3(z)=\{y_1, y_2, y_3\}$, $n_4(z)\le1$ and $n_5^-(z)=0$ such that $L(y_1) \neq \emptyset$ and $L(y_2)=L(y_3)=\emptyset$.
 %where $L(y_i):=\{v\in N_4(y_i):d(v)=1\}$ for any $i\in[3]$.
 Then $z$ will send ${f_2(z)}/{2}$ charge to $y_1$, and ${f_2(z)}/{4}$ charge  to $y_2$,  and ${f_2(z)}/{4}$ charge  to $y_3$, as depicted in Figure \ref{4.6.3}($b$).

\begin{figure}[htbp]
	\centering
	\includegraphics[scale=0.6]{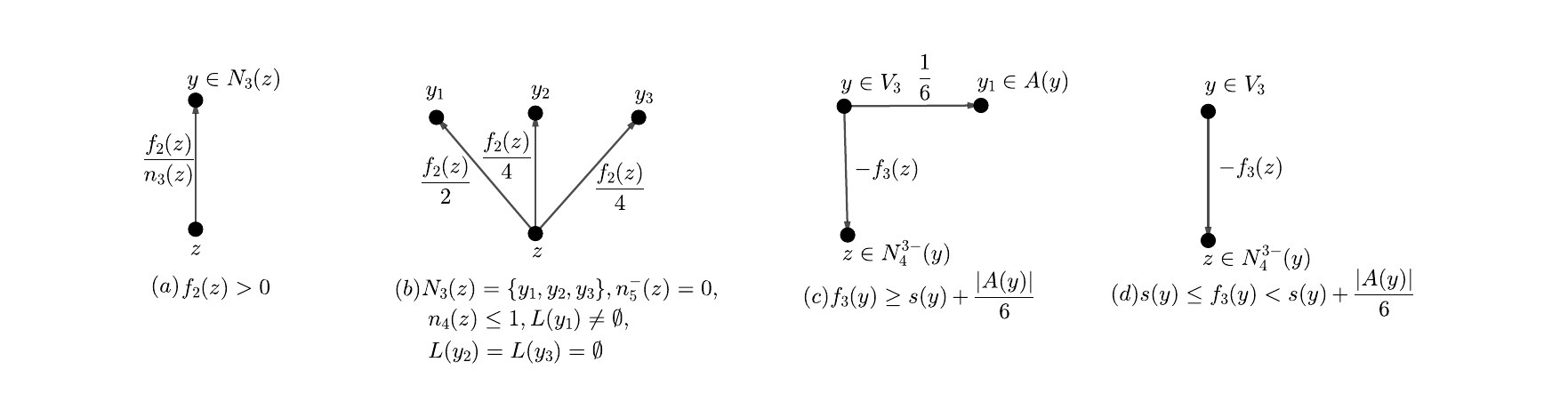}\vspace{-7mm}
	\caption{($a,b$) and ($c,d$) are the explanations of Definition 4.6(3) and (4), respectively. }
	\label{4.6.3}
\end{figure}

%For $z\in V_4$, let $N_3(z)=\{y_1, y_2, y_3\}$, $n_4(z)\le1$ and $n_5^-(z)=0$. Let $L(y_i):=\{v\in N_4(y_i):d(v)=1\}$ for any $i\in[3]$. If $L(y_i) \neq \emptyset$ and $L(y_j)= \emptyset$ for any $j\in[3]\setminus i$, then $f_3(y_i)=f_2(y_i)+\sum_{v\in N_4^{2+}(y_i)\setminus z}\frac{f_2(v)}{n_3(v)}+\frac{f_2(z)}{2}$, $f_3(y_j)=f_2(y_j)+\sum_{v\in N_4^{2+}(y_j)\setminus z}\frac{f_2(v)}{n_3(v)}+\frac{f_2(z)}{4}$ and $f_3(z)=f_2(z)-\frac{f_2(z)}{2}-\frac{f_2(z)}{4}\times 2=0$.

\item Several vertices of $V_3$ send charges to their neighbours in $V_3\cup V_4$. Concretely, for any $y\in V_3$, let $s(y):=\sum_{z \in N_4^{3-}(y)}-f_3(z)$ and $A(y)=\{y_1\in N_3(y)|f_3(y_1)-s(y_1)<0\}$, if $f_3(y)\ge s(y)+\frac{1}{6}|A(y)|$, then $f_4(y)=f_3(y)-s(y)-\frac{1}{6}|A(y)|$, $f_4(z)=f_3(z)+(-f_3(z))$ for any $z\in N_4^{3-}(y)$ and $f_4(y_1)=f_3(y_1)+\frac{1}{6}$ for any $y_1\in A(y)$, as depicted in Figure \ref{4.6.3}($c$); if $s(y)+\frac{1}{6}|A(y)|>f_3(y)\ge s(y)$, then $f_4(y)=f_3(y)-s(y)$ and $f_4(z)=f_3(z)+(-f_3(z))$ for any $z\in N_4^{3-}(y)$, as depicted in Figure \ref{4.6.3}($d$).

%\begin{figure}[htbp]
%	\centering
%	\includegraphics[scale=0.5]{464.pdf}
%	\caption{Definition 4.6(4). }
%	\label{4.6.4}
%\end{figure}

	\item Several vertices of $V_3$ send charges to their neighbours in $V_4$. Concretely, for $y \in V_3$, if $f_4(y)\ge \sum_{z \in N_4^{4-}(y)}-f_4(z)$, then $f_5(y)=f_4(y)-\sum_{z \in N_4^{4-}(y)}-f_4(z)$ and  $f_5(z)=f_4(z)+(-f_4(z))=0$ for any $z \in N_4^{4-}(y)$, as depicted in Figure \ref{4.6.567}($a$).

	%\item For any $y \in V_3$, denote $A'=\{z \in N_4(y): f_4(z)<0\}$,
%	if $f_4(y)\ge -\sum_{z \in A'}f_4(z)$, then $f_5(y)=f_4(y)+\sum_{z \in A'}f_4(z)$ and $f_5(z)=0$ for $z \in A'$.

\item Several $\frac{1}{6}$-vertices of $V_3$ send charges to their neighbours in $V_3$. Concretely, for any two adjacent vertices $y_1,y_2\in V_3^1$ with $N_2(y_i)=\{x_i\}$, if $f_5(x_1)\ge0$, $f_5(x_2)<0$, $f_5(y_1)\ge\frac{1}{6}$ and $f_5(y_2)=0$, then let $f_6(y_1)=f_5(y_1)-\frac{1}{6}$, $f_6(y_2)=f_5(y_2)+\frac{1}{6}$, as depicted in Figure \ref{4.6.567}($b$).

%For $x\in V_2$, $y\in N^1_3(x)$, $y_1\in N^1_3(y)$, $x_1\in N_2(y_1)$ with $f_5(x)<0$, $g(x_1) \ge \frac{1}{6}$, if $n_4^{*-}(y)=0$,  $f_5(y)=0$ and $f_5(y_1) \ge \frac{1}{6}$,  then $f_6(y)=f_5(y)+\frac{1}{6}$, $f_6(y_1)=f_5(y_1)- \frac{1}{6}$.

\item Every vertex of $V_3\cup V_4$ sends charge to some vertices in $V_2$. Concretely, for $x \in V_2$, if  $f_6(x)+ \sum_{y \in N_3^{6+}(x)}{f_6(y)}/{n_2(y)}	\ge \sum_{y \in A_4(x)\cup N_3^{6-}(x)}-f_6(y)$, then $f_7(x)= f_6(x)+\sum_{y \in N_3^{6+}(x)}{f_6(y)}/{n_2(y)}	+\sum_{y \in A_4(x)\cup N_3^{6-}(x)} f_6(y)$, $f_7(y)=f_6(y)-f_6(y)=0$ for any $y \in N_3(x)\cup A_4(x)$, where $A_4(x)=\{y\in N_4(N_3(x))|f_6(y)<0\}$, as depicted in Figure \ref{4.6.567}($c$).\vspace{-5mm}
    \begin{figure}[htbp]
	\centering
	\includegraphics[scale=0.6]{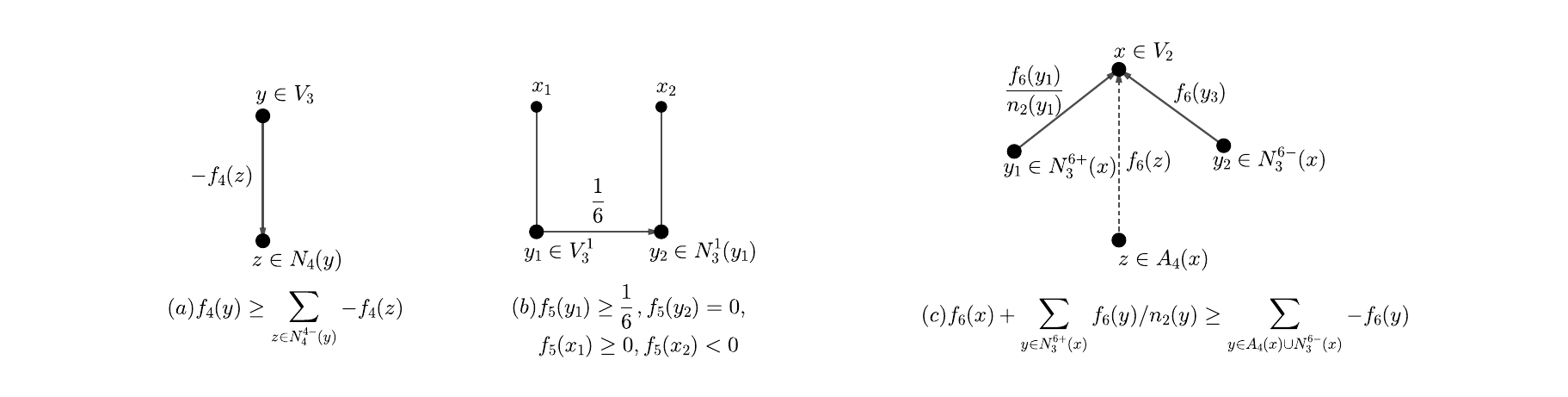}\vspace{-7mm}
	\caption{($a$), ($b$) and $(c)$ are the explanations of Definition 4.6(5), (6) and $(7)$, respectively. }
	\label{4.6.567}
\end{figure}

%\begin{figure}[htbp]
%	\centering
%	\includegraphics[scale=0.5]{467.pdf}
%	\caption{Definition 4.6(7). }
%	\label{4.6.7}
%\end{figure}
		\end{enumerate}
	\end{defn}

 Clearly,
\begin{align}
\sum_{x \in V(G)\setminus V_1} f_7(x)=\sum_{x \in V(G)\setminus V_1}g^*(x).
\end{align}

Observe that for any $x\in V(G)$ and $1\le i<j\le 6$, we have\medskip

($\divideontimes$) if $f_i(x)\ge0$, then $f_j(x)\ge 0$;
 if $f_i(x)<0$, then $f_1(x)\le \cdots \le f_i(x)<0$.\medskip

In order to prove Theorem~\ref{main}, it is enough to prove Theorem~\ref{f7vge0} by inequalities~(\ref{g}) and (4), which will be proved in Section 8. Before we do this, we firstly prove the following two results in Sections 6 and 7, which indicate that the number of vertices in $V_i^{5-}$ for any $i\in\{3,4\}$ are few.
%Their proofs  are stated in Sections 5 and 6.
\vspace{-8pt}

\begin{thm}\label{-320}
Let $x \in V_2$. Then $\sum_{y\in N_3(x)}n_4^{5-}(y)\le1$.\vspace{-4pt}
\end{thm}

\begin{thm}\label{-3little}
Let $x \in V_2$. Then $n_3^{5-}(x)\le1$.  Moreover, if $n_3^{5-}(x)=1$, then $n_4^{5-}(y)=0$ for any $y\in N_3^{5+}(x)$.\vspace{-4pt}
\end{thm}

\begin{thm}\label{f7vge0}
For all $x\notin V_1$, we have $f_7(x)\ge0$.\vspace{-4pt}
\end{thm}

Now we are ready to prove Theorem~\ref{main}.\medskip

\noindent{\bf Proof of Theorem~\ref{main}:} If $\delta(G)\ge3$, then $e(G)\ge3n/2>4n/3-2$. So we may assume that $\delta(G)\le2$. If $\delta(G)=1$, then $\sum_{v \in V_1}g(v)=2\times(\frac{1}{2}-\frac{4}{3})=-\frac{5}{3}$. If $\delta(G)=2$, then $\sum_{v \in V_1}g(v)=2\times(\frac{1}{2}-\frac{4}{3})+(1-\frac{4}{3})=-2$.
By Theorem~\ref{f7vge0} and inequalities~(2) and (4), $\sum_{v\in V(G)\setminus V_1}g(v)\ge0$. Therefore, by Ineq.~(1), $e(G)\ge \sum_{v\in V_1}g(v)+\sum_{v \in V(G)\setminus V_1}g(v)+\frac{4}{3}n\ge \frac{4n}{3}-2$, as desired.\qed

\section{Several Lemmas}
\begin{lem}\label{4-5-0}
Let $y\in V_3$. For $z \in N^{*-}_4(y)$ with $d(z)\ge2$, we have \vspace{-8pt}
		\begin{enumerate}[(a)]
\item $n^{-}_5(z)=0$ and $d(z)=2$,\vspace{-8pt}
	\item  $G$ has a $3$-path $P=zww_1z_1$ such that $d(w)=d(w_1)=d(z_1)=2$ and $yz_1\notin E(G)$.\vspace{-4pt}
	\end{enumerate}			
					
\end{lem}

\begin{lem} \label{34-little}
For any $ y \in V_3$, we have $n^{*-}_4(y) \leq 1$.\vspace{-4pt}
	\end{lem}
%\noindent{\bf Proof of Lemma~\ref{34-little}:} Suppose not. Let $z_1,z_2\in N^{*-}_4(y)$. We claim that $d(z_i)\ge2$ for any $i\in[2]$. W.l.o.g., suppose $d(z_1)=1$. Then $d(z_2)\ge2$, else $G+z_1z_2\notin \mathcal{F}_G$. By Lemma~\ref{4-5-0}($a$), $d(z_2)=2$.
%%{\color{red}By Lemma~\ref{4-5-0}($a$), $d(z_2)=2$, which  contradicts  to the fact $d(z_2)\ge 3$ by Lemma \ref{new}, as claimed.}
%Let $w\in N(z_2)\setminus y$. By Lemma~\ref {path joint}, $G$ has a $4$-path with ends $y$ and $w$ containing $z_2$ because $G+z_1w\in \mathcal{F}_G$. But this is impossible since $d(z_2)=2$,
% By Lemma~\ref{4-5-0}, for any $i\in[2]$, we see $d(z_i)=2$ and $G$ has a $3$-path $P^{i}=z_iw_iw_{i+2}z_{i+2}$ such that $d(w_i)=d(w_{i+2})=d(z_{i+2})=2$ and $yz_{i+2}\notin E(G)$.  This implies $V(P^1)\cap V(P^2)=\emptyset$. But then $G+z_1z_2\notin \mathcal{F}_G$, a contradiction.
%% Hence, $V(P^1)\cap V(P^2)=\emptyset$. Let $C(z_1z_2)$ be a $6$-cycle of $G+z_1z_2$. Then $y\in V(C(z_1z_2))$ because $d_G(z_i)=2$ and $d_G(w_i)=2$ for each $i\in[4]$. Hence, either $yz_4\in E(G)$ or $yz_3\in E(G)$. But then $yz_4\notin E(G)$ and $yz_3\notin E(G)$, by Proposition \ref{4-5-0}(iii), a contradiction.
%\qed\medskip

\noindent{\bf Remark 2.} It is easy to check that $n_5^-(v)\le1$ for any $v\in V_4$. From the above Definition~\ref{law2}, we see that $f_6(v)=\cdots=f_1(v)=g^*(v)$ for any $v\in V_2$, $f_2(v)=f_1(v)=g^*(v)$  for any $v\in V_3$,  $f_1(v)\ge g^*(v)-\frac{1}{3}n_5^-(v)\ge g^*(v)-\frac{1}{3}$ for any $v\in V_4$ and $f_2(v)\ge f_1(v)-\frac{1}{6}n_4^{1-}(v)\ge  f_1(v)-\frac{1}{6}n_4(v)$ for any $v\in V_4^2$.\medskip

%By the Definition \ref{law2}, for any $x \in V(G)$, if $f_i(x)\ge  0$, then $f_{i+1}(x) \ge 0$ and $f(x) \ge 0$ where $i \in \{2,3,4\}$. Hence, if we want to show $f(x)\ge 0$, we just show $f_i(x) \ge 0$ for some $i \in \{2,3,4,5\}$. 	
	
Let $t_v^*:=\frac{g^*(v)}{n_{i-1}(v)}$ for any $v\in V_i$ and $i\in\{3,4,5\}$, $t_v^2:=\frac{f_2(v)}{n_3(v)}$ for any $v\in V_4$ and $t_v^j:=\frac{f_j(v)}{n_2(v)}$ for any $v\in V_3$ and $j\in\{5,6\}$.

\begin{lem}\label{55-0}
For any $z\in V_4$ with $n_5^-(z)\ge 1$, we have $f_2(z)\ge0$; moreover, $f_3(y)\ge0$ for any $y\in N_3(z)$.\vspace{-4pt}
\end{lem}
\medskip

By $(\divideontimes)$, Lemmas~\ref{34-little} and \ref{55-0}, one can easily see the following  corollary.

\begin{cor}\label{f2g}
For any $x\in V_3$, $N_4^{i-}(x)\subseteq N_4^{*-}(x)$ and $n_4^{i-}(x)\le1$ for any $2\le i\le5$.\vspace{-4pt}
\end{cor}

\begin{lem} \label{4-53-0 }
	%Let $y\in V_3$ with $N_5^-(N_4(y))\ne \emptyset$.  If $z\in N_4^{*-}(y)$, then $f_5(z)\ge0$.%\vspace{-12pt}

For any $y\in V_3$ with $N_5^-(N_4(y))\ne \emptyset$, we have $n_4^{5-}(y)=0$.\vspace{-4pt}
%\textcolor{blue}{Let $z_1\in V_4$ with $n_5^-(z_1)\ge1$. If $z\in N_4^{*-}(N_3(z_1))$, then $f_5(z)\ge0$.}
\end{lem}

\begin{lem}\label{law2.3.1}
Let $z\in V_4$ and $N_3(z)=\{y_1, y_2, y_3\}$.
If  $N_2(y_1)=N_2(y_2)=N_2(y_3)$ and $n_2(y_1)=1$, then $n_4^{5-}(y_i)=0$ for any $i\in[3]$.\vspace{-4pt}
\end{lem}

\begin{lem}\label{4-0}
%Let $y\in V_3$. If $z\in L(y)$, then $f_5(y)\ge0$ and $f_5(z)\ge0$.\vspace{-12pt}
For any $y\in V_3$ with $L(y)\ne\emptyset$, we have  $n_4^{5-}(y)=0$, moreover $f_5(y)\ge0$.\vspace{-4pt}
\end{lem}

\begin{lem}\label{yz16-40}
Let $y\in V_3^{5-}$.\vspace{-8pt}
\begin{enumerate}[(a)]
\item If  $z\in N_4^{5-}(y)$, then $n_4^+(y)=0$.\vspace{-8pt}
\item If $z\in N_4^1(y)$ and $N_4(z)=\{z_1\}$, then $n_5^-(z_1)=0$.
\end{enumerate}
\end{lem}

For readability, we have removed all proofs of these lemmas to Appendix.

\section{Proof of Theorem~\ref{-320}}

 Suppose not.  Let $z_i\in N_4^{5-}(y_i)$ for any $i\in[2]$, where $y_i\in N_3(x)$. By ($\divideontimes$) and Corollary~\ref{f2g}, $N_4^{5-}(y_i)=N_4^{3-}(y_i)=N_4^{*-}(y_i)=\{z_i\}$ and $y_1\ne y_2$.
By Lemmas~\ref{4-0} and \ref{4-5-0}($a,b$), for any $i\in[2]$, $G$ has a $3$-path $P^i=z_iw_iw_i'z_i'$ such that $d(v)=2$ for any $v\in V(P^i)$ and $z_i'y_i\notin E(G)$. Clearly, $f_3(z_i)\ge g^*(z_i)\ge-\frac{1}{6}$. Hence, $f_3(y_i)<\frac{1}{6}$ for any $i\in [2]$, else $f_4(z_i)\ge f_3(z_i)+\frac{1}{6}\ge0$ because $N_4^{3-}(y_i)=\{z_i\}$.
Note that either $P_1=P_2$ or $V(P_1)\cap V(P_2)=\emptyset$. Because $d(z_i')=d(w_i')=2$, we see $z_i'\in V_4^{*-}$ and so $z_i'y_{3-i}\notin E(G)$ if $z_i'\neq z_{3-i}$ for any $i \in [2]$.
Since $G+z_1z_2\in \mathcal{F}_G$, we see $G$ has a $3$-path $P(y_1y_2)$. Let $P(y_1y_2)=y_1a_1a_2y_2$. We next prove two claims.
\medskip

\noindent{\bf Claim 1.} $x\notin \{a_1,a_2\}$
\medskip

\noindent{\bf Proof.} W.l.o.g., suppose $x=a_1$. Then
$N_2^+(y_2)\cup (N_3^2(x)\setminus y_2)=\emptyset$, else $f_3(y_2)\ge g^*(y_2)\ge\frac{1}{6}$. Hence, $a_2\in V_3^1$ and $g(x)<0$.
By the choice of $\alpha$, we see $\alpha_1\alpha_2\notin E(G)$ if $\delta(G)=2$. This implies $g(y_1)=\frac{1}{6}$ because $G+z_1\alpha\in\mathcal{F}_G$, $d(z_1,\alpha)=4$ and $g(x)<0$. Clearly, $g(y_2)\le\frac{1}{6}$, else $f_3(y_1)\ge g^*(y_1)\ge\frac{1}{6}$.
Note that $G$ has a $5$-path $P(y_2\alpha)$.  Let $P(y_2\alpha)=y_2b_1b_2b_3\alpha_i\alpha$. By the choice of $\alpha$,  $\alpha$ belongs to no $4$-cycle because $d(z_i)=2$ and $z_i$ belongs to no $4$-cycle. Then $b_1\in V_4$ and $b_2\in V_3$ because $d(y_2,\alpha)=3$, $N_3^2(x)\setminus y_2=\emptyset$ and $g(x)<0$.
By Lemma~\ref{4-53-0 }, $n_5^-(b_1)=0$. By Observation~\ref{ob3}(2), $g^*(b_1)\ge \frac{1}{6}n_3(b_1)+\frac{1}{6}n_4(b_1)$ because $y_2\in N_3^+(b_1)$ and $n_3(b_1)\ge2$, which yields $t_{b_1}^2\ge\frac{1}{6}$. But then $f_3(y_2)\ge f_2(y_2)+t_{b_1}^2\ge \frac{1}{6}$, a contradiction. This proves Claim 1.\qed
\medskip

\noindent{\bf Claim 2.} $|\{a_1,a_2\}\cap V_4|=1$.
\medskip

\noindent{\bf Proof.} Obviously, $|\{a_1,a_2\}\cap V_4|\ge1$, else $f_3(y_i)\ge g^*(y_i)\ge \frac{1}{6}$ because $n_3^2(x)+n_3^2(y_i)\ge1$ for some $i\in[2]$.  Suppose $|\{a_1,a_2\}\cap V_4|=2$. By Lemma~\ref{4-53-0 }, for any $i\in[2]$, $n_5^-(a_i)=0$ and $f_1(a_i)\ge0$, which means $n_4^{1-}(a_i)\le n_4(a_i)-1$.
Then $g(y_i)<0$ for any $i\in[2]$, else let $g(y_1)\ge0$, we see $f_3(y_1)\ge f_2(y_1)+t_{a_1}^2\ge \frac{1}{6}$ because $f_1(a_1)\ge g^*(a_1)\ge\frac{1}{6}n_3(a_1)+\frac{1}{6}(n_4(a_1)-1)$ and $t_{a_1}^2\ge\frac{1}{6}$.
Hence, $G$ has a $3$-path $P(xy_i)$ for some $i\in[2]$ as $G+z_1y_2\in\mathcal{F}_G$. Let $P(xy_i)=xb_1b_2y_i$. Then $b_1\in V_3$ and $b_2\in V_4$. Since $G+z_ib_1\in\mathcal{F}_G$, we have $g(b_1)\ge\frac{1}{6}$ or $G-y_i-b_1$ has a $2$-path with ends $x$ and $b_2$. By Observation~\ref{ob3}(2), $f_1(b_2)\ge g^*(b_2)\ge\frac{1}{6}n_3(b_2)+\frac{1}{6}n_4(b_2)$ because $n_5^-(b_2)=0$ and $n_3^+(b_2)+n_3(b_2)\ge3$, which implies $t_{b_2}^2\ge\frac{1}{6}$. But then $f_3(y_i)\ge f_2(y_i)+t_{b_2}^2\ge\frac{1}{6}$, a contradiction.  \qed
\medskip

By Claims 1 and 2, w.l.o.g., let $a_1\in V_4$. Then $a_2\in V_3$ and so $t_{a_1}^2\ge\frac{1}{6}$ because $n_3^+(a_1)+n_3(a_1)\ge3$ and $n_5^-(a_1)=0$. But then $f_3(y_1)\ge f_2(y_1)+t_{a_1}^2\ge\frac{1}{6}$, a contradiction. This completes the proof of Theorem~\ref{-320}.\qed

\section{Proof of Theorem~\ref{-3little}}
 We first show that $n_3^{5-}(x)\le1$. Suppose not. Let $y_1,y_2\in N_3^{5-}(x)$. By $(\divideontimes)$ and Corollary~\ref{f2g}, $f_j(y_i)<0$ for any $j\in[5]$ and $i\in[2]$, and $y_1,y_2\in N_3^{*-}(x)$.
\medskip

\noindent{\bf Claim 1.} $n_4^2(y_i)=0$, $n_4^1(y_i)\le1$, $N_5^-(N_4(y_i))=\emptyset$
and there is no $s$-path $P(xy_i)$ for any $s \in \{2,3\}$.\medskip

\noindent{\bf Proof.} The desired results hold by Observation~\ref{ob}(1) and Lemma~\ref{55-0}.\qed\medskip

\noindent{\bf Claim 2.} Suppose $N_4^1(y_i)=\{z_i\}$, $N_4(z_i)=\{z_i'\}$, $w_i'\in N_5(z_i')$ and $w_i\in N_5(z_i)$. Then $n_4^2(z_i)=n_5^2(z_i)=n_5^2(w_i)=n_5^-(z_i')=n_5^2(w_i')=0$. Furthermore, $f_1(z_i')<\frac{1}{3}$ and $f_1(v)<\frac{1}{6}$ for any $v\in N_4(y_i)$.
\medskip

\noindent{\bf Proof.} By Observation~\ref{ob}(2), $n_4^2(z_i)=0$.  By Lemma \ref{yz16-40}(2), $n_5^-(z_i')=0$. Clearly, $f_2(y_i)=g^*(y_i)\ge-\frac{1}{6}$. Hence, $f_1(v)<\frac{1}{6}$ for any $v\in N_4(y_i)$, otherwise let $f_1(v_0)\ge\frac{1}{6}$ for some $v_0\in N_4(y_i)$, $t_{v_0}^2\ge\frac{1}{6}$ as $n_4^{1-}(v_0)=0$, which yields $f_3(y_i)\ge f_2(y_i)+t_{v_0}^2\ge0$. Moreover, $f_1(z_i')<\frac{1}{3}$, otherwise by Definition~\ref{law2}(2), $f_2(z_i)\ge\frac{1}{6}$ and so $f_3(y_i)\ge0$. Next we prove $n_5^2(z_i)=n_5^2(w_i)=0$. Suppose $w_i\in N_5^2(z_i)$ or $n_5^2(w_i)\ne0$. By Ineq.~(\ref{g*}), $g^*(w_i)\ge\frac{1}{6}n_4(w_i)$ because $z_i\in N_4^+(w_i)$, which means $t_{w_i}^*\ge\frac{1}{6}$. But then $f_1(z_i)\ge g^*(z_i)+t_{w_i}^*\ge\frac{1}{6}$, a contradiction, as desired.
It remains to show $n_5^2(w_i')=0$. Suppose not. Then $w_i'\in V_5^1$ and $g^*(w_i')<\frac{1}{3}$, else, $t_{w_i'}^*\ge \frac{1}{3}$ and so $f_1(z_i')\ge t_{w_i'}^*\ge\frac{1}{3}$.
%If
%$g^*(w_i')=\frac{1}{3}$, then $f_1(z_i')\ge \frac{1}{3}$, a contradiction.
By Definition~\ref{law1}(1.1), we see there exists a vertex $w_i''\in N_5(z_i')\setminus w_i'$ such that $t^*_{w_i''}\ge \frac{1}{6}$. Note that $g^*(w_i')\ge\frac{1}{6}$ and so $t^*_{w_i'}\ge \frac{1}{6}$, which follows that $f_1(z_i')\ge t^*_{w_i'}+t^*_{w_i''}\ge\frac{1}{3}$, a contradiction.
\qed\medskip

\noindent{\bf Claim 3.}
Suppose $z_i\in N_4^-(y_i)$ and $z_i'\in N_5^2(z_i)$. Then  $d(z_i')\le3$. Moreover, $f_1(z_i)<\frac{1}{3}$.
\medskip

\noindent{\bf Proof.} Clearly, $f_1(z_i)<\frac{1}{3}$, else $t_{z_i}^2\ge\frac{1}{3}$ and so $f_3(y_i)\ge f_2(y_i)+t_{z_i}^2\ge0$. Suppose $d(z_i')\ge 4$.  By Ineq.~(\ref{g*}),  $g^*(z_i')\ge\frac{1}{3}n_4(z_i')$, that is, $t_{z_i'}^*\ge\frac{1}{3}$. By Claim 1, $n_5^-(z_i)=0$.  Hence, $f_1(z_i)\ge g^*(z_i)+t_{z_i'}^*\ge\frac{1}{3}$, a contradiction. \qed\medskip

\noindent{\bf Claim 4.} $G$ has no $4$-cycle containing $y_i$.
\medskip

\noindent{\bf Proof.} W.l.o.g., suppose $G$ has a  $4$-cycle with vertices $y_1,u_1,u_2,u_3$ in order. By Claim 1, we have $u_1,u_3\in V_4$ and $u_2\in V_5$ and $u_1u_3 \notin E(G)$. By Claims 1, 2 and 3, $n_5^-(u_i)=0$, $g(u_i)<0$ for any $i\in\{1,3\}$ and $d(u_2)\le3$. We assert $d(u_2)=2$. Suppose not. By Ineq.~(\ref{g*}), $g^*(u_2)\ge \frac{1}{6}n_4(u_2)$ because $n_4(u_2)+n_5(u_2)=d(u_2)\ge3$ and $n_4(u_2)\ge2$, which means $t_{u_2}^*\ge\frac{1}{6}$. Hence, $f_2(u_i)= f_1(u_i)\ge g^*(u_i)+t_{u_2}^*\ge \frac{1}{6}$ for any $i\in\{1,3\}$, which means $t_{u_i}^2\ge\frac{1}{6}$. But then $f_3(y_1)\ge f_2(y_1)+t_{u_1}^2+t_{u_3}^2\ge 0$, a contradiction, as asserted. By Claim 1, $G-y_1-u_2$ has a $2$-path with vertices $u_1,u_4,u_3$ in order such that $u_4\in V_5$ as $G+u_2x\in\mathcal{F}_G$. By Ineq.~(\ref{g*}),  $g^*(u_2)\ge\frac{1}{6}n_4(u_2)$ since $n_4(u_2)\ge2$ and $u_1,u_3\in N_4^{-1}(u_2)$. Similarly, we have $f_3(y_1)\ge0$, a contradiction.  \qed\medskip

Clearly, $d(y_i)\ge2$ for some $i\in[2]$, else $G+y_1y_2\notin\mathcal{F}_G$. W.l.o.g., assume $d(y_1)\ge2$. Let $z_1\in N_4(y_1)$. Note that $d(y_2)\ge2$, otherwise by Claims 1 and 4, $G$ has no $4$-path with ends $x$ and $z_1$ containing $y_1$, which yields $G+z_1y_2\notin\mathcal{F}_G$. By Claim 1, $z_1y_2\notin E(G)$, which means $G$ has a $5$-path with vertices $z_1,a_1,a_2,a_3,z_2,y_2$ in order. By Claims 1 and 4, $a_1\notin V_3$ and $z_2\ne x$ which implies $z_2\in N_4(y_2)$. Note that $z_1z_2\notin E(G)$, otherwise $a_1,a_3\in V_5$ and $a_2\in V_4\cup V_5$ which violates Claim 2. This yields $G$ has a $5$-path with vertices $z_1,b_1,b_2,b_3,b_4,z_2$ in order. Note that $\{b_1,b_4\}\subseteq V_3\cup V_4\cup V_5$. Clearly, by Lemma~\ref{55-0}, $n_5^-(z_i)=0$ for any $i\in[2]$, which yields $f_1(z_i)\ge g^*(z_i)$.\medskip

\noindent{\bf Claim 5.} $|\{b_1,b_4\}\cap V_3|\le1$.
\medskip

\noindent{\bf Proof. } Suppose not. By Claim 1, $b_1=y_1$, $b_4=y_2$, $b_2,b_3\in V_4^1$ and $b_1b_3\notin E(G)$. This implies $G$ has a $5$-path with vertices $b_1,c_1,c_2,c_3,c_4,b_3$ in order. By Claims 1 and 4, $c_4\ne y_2$ and $c_1\neq x$. Then $c_1\in V_4$ as $g(y_1)<0$.
Moreover, $c_4\notin V_5$, otherwise by Claim 2, $c_3\in V_5$, $c_2\in V_4$, which yields  $g(c_1)\ge\frac{1}{6}$ and so $n_4^2(y_1)\ge1$ or $n_4^1(y_1)\ge2$ violating  Claim 1. Thus, $c_4=b_2$.  By Claims 1 and 2, $c_3,c_2\in V_5$ and $d(c_2)=d(c_3)=2$. But then  $G+c_1c_3\notin\mathcal{F}_G$ by Claims 1 and 4, a contradiction. \qed\medskip

%\textcolor{blue}{In fact, we further prove that $|\{b_1,b_4\}\cap V_3|=0$, which is shown in Appendix A.}\medskip

\noindent{\bf Claim 6.} $|\{b_1,b_4\}\cap V_3|=0$.
\medskip

\noindent{\bf Proof.} Suppose not. By Claim 5, $|\{b_1,b_4\}\cap V_3|=1$. We may assume $b_1\in V_3$. Then $b_1=y_1$ and $b_4\in V_4\cup V_5$.\medskip

{\bf Case A: $b_4\in V_4$.}
Then $a_3=b_4$, otherwise $a_2,a_3\in V_5^1$ and $a_1\in V_4$ by Claim 2, we see $g^*(a_3)\ge\frac{1}{6}$ and so $t_{a_3}^*\ge\frac{1}{6}$, which yields $f_1(z_2)\ge g^*(z_2)+t_{a_3}^*\ge\frac{1}{6}$ violating Claim 2.
Note that $a_2\in V_3\cup V_5$. We next prove \medskip

\noindent{$(a)$ $a_2\in V_5$.}\medskip

To see why $(a)$ is true, suppose $a_2\in V_3$. By Claim 1, $a_1\in V_4$. Observe that $f_1(a_i)<\frac{1}{6}$ for any $i \in \{1,3\}$, otherwise let $f_1(a_1)\ge\frac{1}{6}$, by Definition \ref{law2}(2), $f_2(z_1)\ge \frac{1}{6}$ and so $t_{z_1}^2\ge\frac{1}{6}$ as $g^*(a_2)\ge0$, which implies $f_3(y_1)\ge f_2(y_1)+t_{z_1}^2\ge0$.
It follows that $b_3=a_2$, otherwise $b_3\in V_5$ and $b_2\in V_4$ by Claim 2, we see $g^*(b_3)\ge \frac{1}{6}n_4(b_3)$ and so $t_{b_3}^*\ge\frac{1}{6}$, which yields $f_1(a_3)\ge g^*(a_3)+t_{b_3}^*\ge\frac{1}{6}$. Furthermore, $b_2=x$ by Claim 1. Since $G+z_1a_2\in \mathcal{F}_G$, we set $C_6(z_1a_2)=z_1c_1c_2c_3c_4a_2$. Note that $g(a_2)<0$, else $f_1(a_3)\ge g^*(a_3)\ge\frac{1}{6}$. By Claims 1 and 4, $c_4\in V_4$.
%, otherwise $c_4=x$ and $y_1=c_i$ for some $i\in[3]$, which violates Claims 1 and 4.
We assert $c_1\notin V_3$. Suppose not. Then $c_1=y_1$.  Note  that $c_2=x$, otherwise by Claim 1 and $g(y_1)<0$, $c_2\in V_4$ and $c_3\in V_5$, which yields $G$ has a copy of $C_6$ with vertices $y_1,x,a_2,c_4,c_3,c_2$ in order. So $c_3\in V_3$, which means $g(c_4)\ge\frac{2}{3}$. But then $f_1(a_1)\ge g^*(a_1)\ge \frac{1}{6}$ as $g(c_4)\ge\frac{2}{3}$, a contradiction. Clearly, $c_1\in V_4$, otherwise $c_1, c_2\in V_5^1$ and $c_3\in V_4$ by Claim 2, which implies $g^*(c_1)\ge\frac{1}{6}$ and so $f_1(z_1)\ge \frac{1}{6}$. Then $c_1=a_1$ and $c_2\in V_5$. Moreover, $c_2, c_3\in V_5^1$, otherwise $ g^*(c_2)\ge \frac{1}{6}n_4(c_2)$ and so $t_{c_2}^*\ge\frac{1}{6}$, which implies $f_1(a_1)\ge g^*(a_1)+t_{c_2}^*\ge\frac{1}{6}$. Then $G$ has a 3-path $P(a_2a_1)$ or $P(a_2c_4)$ as $G+c_2c_4\in \mathcal{F}_G$. We have $f_1(a_1)\ge \frac{1}{6}$, a contradiction. This proves $(a)$.\medskip

By $(a)$,  $a_1 \in  V_5$ and $d(a_2)=2$, otherwise $g^*(a_2)\ge \frac{1}{3}n_4(a_2)$ and so $t_{a_2}^*\ge\frac{1}{3}$, which implies $f_1(b_4)\ge g^*(b_4)+t_{a_2}^*\ge\frac{1}{3}$ violating Claim 2. It is easy to check that $d(a_1)=2$ by Claim 2.
Moreover, $g(z_1)<0$, otherwise $g^*(a_1)=\frac{1}{6}$ and so $t_{a_1}^*\ge\frac{1}{6}$, which yields $f_1(z_1)\ge g^*(z_1)+t_{a_1}^*\ge \frac{1}{6}$ violating Claim 2.
Note that $b_2\in V_2\cup V_4$ as $b_1=y_1$. We assert  $b_2\in V_4$. Suppose not.  Then $b_2=x$ and $b_3\in V_3$. Since $G+a_2x\in \mathcal{F}_G$ and $d(a_1)=d(a_2)=2$, $a_2b_4\in E(C_6(a_2x))$, otherwise $a_2a_1\in E(C_6(a_2x))$ and $G$ has a 3-path with ends $x$ and $z_1$, which violates $g(z_1)<0$ and $g(y_1)<0$. Then $G$ has a 4-path with ends $x$ and $b_4$ containing $b_3$, which means $G$ has a 3-path  with ends $v$ and $b_3$, where $v\in\{x, b_4\}$ as $g(z_2)=g(b_4)=\frac{1}{6}$.
%By Claim 2, $f_1(b_4)<\frac{1}{3}$.
By Definition \ref{law2}(2), we have
$n_4^2(b_3)=0$ and $g(b_3)<0$, otherwise $f_2(z_2)\ge \frac{1}{6}$ and so $t_{z_2}^2\ge\frac{1}{6}$, which means $f_3(y_2)\ge f_2(y_2)+t_{z_2}^2\ge0$.
It follows that  $G$ has no 3-path with ends $x$ and $b_3$. Then $G$ must have a 3-path $P(b_3b_4)$. Let $P(b_3b_4)=b_3b_3'b_4'b_4$. Then $b_4'\in V_5$ and $b_3'\in V_4$. We see $d(b_4')=2$, otherwise $g^*(b_4')\ge \frac{1}{3}n_4(b_4')$ and so $t_{b_4'}^*\ge\frac{1}{3}$, which means $f_1(b_4)\ge g^*(b_4)+t_{b_4'}^*\ge\frac{1}{3}$. Since $G+xb_4'\in \mathcal{F}_G$, there exists one vertex of $\{b_4,b_3'\}$, say $b_4$, such that $G-b_4'$ has a $4$-path with ends $x$ and $b_4$ containing $b_3$. It follows $G-b_4'$ has a $3$-path with ends $b_3$ and $b_4$ containing $b_3'$. That is, there is a vertex $b_4''\in (N_5(b_4)\cap N_5(b_3'))\setminus{b_4'}$.  Then $g^*(v)\ge \frac{1}{6}n_4(v)$ for any $v \in \{b_4', b_4''\}$ and so $t_{v}^*\ge\frac{1}{6}$. But then $f_1(b_4)\ge g^*(b_4)+t_{b_4'}^*+t_{b_4''}^*\ge\frac{1}{3}$, a contradiction, as asserted.
Thus, $b_2\in V_4$. Then $b_3\in V_5$. We have $d(b_3)=2$ and $g(b_2)<0$, otherwise $t_{b_3}^*\ge \frac{1}{3}$ and so $f_1(b_4)\ge g^*(b_4)+t_{b_3}^*\ge\frac{1}{3}$.
Since $G+a_1b_3\in \mathcal{F}_G$, we see $b_3b_2\in E(C_6(a_1b_3))$ as $d(a_2)=2$, which implies $G-b_3$ has a 3-path $P(b_2z_1)$ or a 2-path $P(b_2b_4)$. Note that $g^*(b_3)\ge \frac{1}{6}n_4(b_3)$ and so $t_{b_3}^*\ge\frac{1}{6}$. Since $g(b_2)<0$ and $n_5^-(b_2)=0$, we see $t_{b_2}^2\ge\frac{1}{6}$. If the 3-path $P(b_2z_1)$ exists, let $P(b_2z_1)=b_2b_2'z_1'z_1$, then $\{b_2', z_1'\}\subseteq V_5$.
	If  $n_5^2(z_1)\neq 0$, then $g^*(a_1)\ge \frac{1}{6}$ and so $t^*_{a_1}\ge \frac{1}{6}$. If $n_5^2(z_1)=0$,
	then $g_2(a_1)=\frac{1}{3}$, which follows $g^*(z_1')\ge \frac{1}{6}$ and so $t_{z_1'}^*\ge\frac{1}{6}$.
	In both cases,  we see $t^2_{z_1}\ge \frac{1}{6} $ because $n_5^-(z_1)=0$ and $g(z_1)<0$. But then $f_3(y_1)\ge f_2(y_1)+t^2_{z_1}+t^2_{b_2}\ge 0$, a contradiction.
If the 2-path $P(b_2b_4)$ exists, let $P(b_2b_4)=b_2b_2'b_4$. Then $b_2'\in V_5$, otherwise $b_2'=y_1$, we see $f_2(y_1)\ge g^*(y_1)\ge -\frac{1}{6}$, which yields $f_3(y_1) \ge f_2(y_1)+t_{b_2}^2\ge0$. Similarly, $t_{b_2'}^*\ge\frac{1}{6}$ and so $f_1(b_4)\ge g^*(b_4)+t_{b_3}^*+t_{b_2'}^*\ge \frac{1}{3}$ violating Claim 2.\medskip

{\bf Case B:} $b_4\in V_5$. Then $b_2\in V_4$ and $b_3\in V_4\cup V_5$. If $b_3\in V_4$, then $g(z_2)<0$ and $d(b_4)=2$ otherwise $t_{b_4}^*\ge\frac{1}{3}$ and $f_1(b_3)\ge g^*(b_3)+t_{b_4}^*\ge\frac{1}{3}$. Clearly, $a_3\in V_5$.
Note that  $a_3\neq b_4$, otherwise $a_2=b_3$ and $a_1\in V_5$, which means $f_1(b_3)\ge g^*(b_3)+t_{b_4}^*+t_{a_1}^*\ge\frac{1}{3}$ because $t_v^*\ge \frac{1}{6}$ for any $v\in \{a_1, b_4\}$.  It is easy to check
$g^*(v)\ge\frac{1}{6}n_4(v)$ for any $v\in\{b_4,a_3\}$, which means $t_{v}^*\ge\frac{1}{6}$. But then $f_1(z_2) \ge g^*(z_2)+t_{b_4}^*+t_{a_3}^*\ge \frac{1}{3}$, which contradicts to Claim 3. Next we assume that $b_3\in V_5$. By Claim 4, $b_3\neq  a_1$. Moreover, $b_3\neq a_3$ and $b_4\notin \{a_1, a_2\}$, otherwise there is a 6-cycle in $G$, that is, $C_6=z_2y_2xy_1b_2b_3$ when $b_3=a_3$, or $C_6=z_1y_1xy_2z_2b_4$ when $b_4=a_1$, or $C_6=z_1y_1b_2b_3b_4a_1$ when $b_4=a_2$. We next prove that\medskip

\noindent{$(b)$ $|\{b_4, b_3\}\cap \{a_2, a_3\}|\le 1$.}\medskip

To prove $(b)$, suppose $\{b_4, b_3\}=\{a_2, a_3\}$. Then $b_4=a_3$ and $b_3=a_2$. Then $a_1\in V_5$, else $t_{b_3}^*\ge\frac{1}{3}$ and so $f_1(b_2)\ge g^*(b_2)+t_{b_3}^*\ge\frac{1}{3}$, which violates Claim 2. Moreover, $d(a_1)=2$, else $d(a_1)\ge 3$,
$g^*(v)\ge \frac{1}{6}n_4(v)$ as $n_5^2(v)\ge 1$ for any $v \in \{a_1, b_3\}$,  which implies $t_u^2\ge \frac{1}{6}$ for any $u \in \{z_1, b_2\}$ and so $f_3(y_1)\ge f_2(y_1)+t_{z_1}^2+t_{b_2}^2\ge0$. By Claim 3, $d(b_3)=3$.
Since $G+a_1b_2 \in \mathcal{F}_G$ and by Claims 1 and 4, $a_1a_2\in E(C_6(a_1b_2))$ and $G-a_1-a_2$ has a 3-path $P(b_2b_4)$. Let $P(b_2b_4)=b_2b_2'b_4'b_4$. By Claim 2, $g(b_2)<0$ and $g(z_2)<0$.
Then $b_2'\in V_5$, else $b_2'=y_1$, we see $G$ has a copy of $C_6$, that is, $C_6=y_1xy_2z_2b_4b_4'$ when $b_4'=z_1$, or $C_6=y_1z_1a_1b_3b_4b_4'$ when $b_4'\ne z_1$.
We assert $b_4'\in V_5$. Suppose $b_4'\in V_4$.  Then  $b_4'\neq z_2$, else there is a $C_6=b_2'b_2y_1xy_2z_2$ in $G$. By Ineq.~(\ref{g*}), $g^*(v)\ge \frac{1}{6}n_4(v)$ for any $v\in\{b_2',b_3\}$ as $n_4^{-1}(b_2')+n_4^+(b_2')\ge 2$,  $n_4^{-1}(b_3)+n_4^+(b_3)\ge 1$ and $n_4(b_3)+n_5(b_3)\ge3$, that is, $t_{v}^*\ge\frac{1}{6}$. But then $f_1(b_2)\ge g^*(b_2)+t_{b_3}^*+t_{b_2'}^*\ge\frac{1}{3}$, which contradicts to Claim 3, as asserted. So $b_4'\in V_5^2$.  Since $n_5^2(b_4)\ge 2$,  $g^*(b_4)\ge \frac{1}{3}n_4(b_4)$ and so $t_{b_4}^*\ge\frac{1}{3}$. But $f_1(z_2)\ge g^*(z_2)+t_{b_4}^*\ge\frac{1}{3}$, which contradicts to Claim 3. This proves $(b)$.\medskip

\noindent{$(c)$ $|\{b_4, b_3\}\cap \{a_2, a_3\}|=0$.}\medskip

To prove $(c)$, suppose not. By $(b)$, $|\{b_4, b_3\}\cap \{a_2, a_3\}|=1$.
Clearly, $b_4\in \{a_2,a_3\}$, otherwise $b_3=a_2$ and $b_4\neq a_3$, which yields $d(b_3)\ge4$ violating Claims 2 and 3. Then $b_4=a_3$ and $b_3\neq a_2$. By Claim 2, $g(z_2)<0$. By Claim 3,  $d(b_4)=3$.
We assert $d(b_3)=2$. Suppose not. Then $g^*(b_3)\ge \frac{1}{6}n_4(b_3)$ because $b_4\in N_5^2(b_3)$ and $d(b_3)=n_4(b_3)+n_5(b_3)\ge3$ and so $t_{b_3}^*\ge\frac{1}{6}$. By Claim 2, $a_1\in V_5$, otherwise $f_1(b_2)\ge g^*(b_2)+t_{b_3}^*\ge\frac{1}{6}$.  It follows that $g^*(b_4)\ge\frac{1}{3}n_4(b_4)$ because $n_4(b_4)+n_5(b_4)=3$ and $n_5^2(b_4)+n_4^{-1}(b_4)+n_4^+(b_4)\ge2$. %$\{b_3,a_2\}\subseteq N_5^2(b_4)$.
That is, $t_{b_4}^*\ge\frac{1}{3}$. But then $f_1(z_2)\ge g^*(z_2)+t_{b_4}^*\ge\frac{1}{3}$, which contradicts to Claim 3, as asserted.
Since $G+z_1b_3\in \mathcal{F}_G$, let $C_6(z_1b_3)=z_1c_1c_2c_3c_4b_3$. By Claims 1 and 4, $c_4\ne b_2$. Hence, $c_4=b_4$ as $d(b_3)=2$. We assert $c_1\notin V_3$. Suppose not. Then $b_2=c_i$ for some $i\in \{2,3\}$, otherwise $G$ has a copy of $C_6=c_1c_2c_3c_4b_3b_2$. By Claim 1, $b_2=c_2$. Since $d(b_4)=3$, we have $c_3=a_2$.
%Then $g(b_2)\ge0$ or $n_5^2(b_2)\ge1$, which implies
Note that $a_2\in V_4^+\cup V_5^2$. Hence,
$g^*(b_3)\ge\frac{1}{3}$ or $g^*(b_4)\ge \frac{1}{3}n_4(b_4)$. That is, $t_{b_3}^*\ge\frac{1}{3}$ or $t_{b_4}^*\ge\frac{1}{3}$. But then $f_1(b_2)\ge g^*(b_2)+t_{b_3}^*\ge\frac{1}{3}$ or $f_1(z_2)\ge g^*(z_2)+t_{b_4}^*\ge\frac{1}{3}$, a contradiction, as asserted. Since $d(b_4)=3$, we have $c_3\in\{z_2,a_2\}$.

We next show that $c_3\ne a_2$. Suppose not. Then $c_2\ne a_1$, otherwise $c_2,c_1\subseteq V_5$, we see $g^*(v)\ge\frac{1}{6}n_4(v)$  for any $v\in \{c_1,c_2\}$ as $n_5^2(v)\ge1$, which implies $t_{v}^*\ge\frac{1}{6}$ and so $f_1(z_1)\ge g^*(z_1)+t_{c_1}^*+t_{c_2}^*\ge\frac{1}{3}$.
We assert $c_1\ne a_1$. Suppose not. By Claim 2, $c_1\in V_5$.
%Then $c_1\in V_5$, otherwise $c_1\in V_4$ and $\{a_2,c_2\}\subseteq V_5$, we see $g^*(a_2)\ge\frac{1}{3}n_4(a_2)$ because $c_1\in N_4^+(a_2)$ and $d(a_2)\ge3$, which implies $f_1(c_1)\ge g^*(c_1)+t_{a_2}^*\ge\frac{1}{3}$ violating to Claim 2.
Moreover, $a_2\in V_4$, otherwise $a_2\in V_5$, we see $t_{c_1}^*\ge\frac{1}{3}$ because $d(c_1)\ge3$ and $n_4^+(c_1)+n_4^{-1}(c_1)+n_5^2(c_1)\ge2$, which yields $f_1(z_1)\ge g^*(z_1)+t_{c_1}^*\ge\frac{1}{3}$ violating Claims 2 and 3. Clearly, $c_2\in V_5$ and $d(c_2)=2$, otherwise $t_{c_1}^*\ge\frac{1}{3}$  and so $f_1(z_1)\ge g^*(z_1)+t_{c_1}^*\ge\frac{1}{3}$. But then $G+z_1c_2\notin \mathcal{F}_G$ because $d(c_1)=3$ and $d(c_2)=2$, a contradiction, as asserted. By Claim 2, $c_1\in V_5$.
%We say $c_1\in V_5$, otherwise $c_1\in V_4$ and so $a_1\in V_5$, we see $t_{a_1}^*\ge\frac{1}{6}$, which implies $f_1(z_1)\ge g^*(z_1)+t_{a_1}^*\ge\frac{1}{6}$ violating to Claim 2.
Furthermore, $c_2\in V_5$, otherwise $a_1\in V_5$ by Claim 2, we see $t_{v}^*\ge\frac{1}{6}$ for any $v\in\{c_1,a_1\}$  and so $f_1(z_1)\ge g^*(z_1)+t_{a_1}^*+t_{c_1}^*\ge\frac{1}{3}$, which violates Claims 2 and 3.
By Claim 2, $a_1\in V_5$.
Note that $a_2\in V_4\cup V_5$. Then
$d(v)=2$ for any $v\in\{a_1,c_1\}$, otherwise we see $t_{v}^*\ge\frac{1}{6}$ and so $f_1(z_1)\ge g^*(z_1)+t_{a_1}^*+t_{c_1}^*\ge\frac{1}{3}$. By Claim 2, $g(z_1)<0$. If $a_2\in V_5$, then $f_1(z_1)\ge g^*(z_1)\ge\frac{1}{3}$, which violates Claim 3.
If $a_2\in V_4$, then $d(c_2)\ge 3$, otherwise, $G+z_1c_2 \notin \mathcal{F}_G$.   By Claim 3, $g^*(c_1)<\frac{1}{3}$. By Definition~\ref{law1}(1), there exists a vertex $c_1'\in N_5(z_1)\cap N_5(c_2)$. It is easy to check that for any $v\in\{c_1,c_1'\}$, $t_{v}^*\ge\frac{1}{6}$ and so $f_1(z_1)\ge g^*(z_1)+t_{c_1}^*+t_{c_1'}^*\ge\frac{1}{3}$, which violates Claim 3, as desired.

Thus, $c_3\ne a_2$ and so $c_3=z_2$.  By Claim 2, $a_1,c_1\in V_5$. Recall that $g(z_2)<0$. Hence, $c_2\in V_5$.
%Hence, $\{c_1,c_2\}\subseteq V_5$, otherwise $c_1\in V_4$, we see $g^*(c_2)\ge \frac{n_4(c_2)}{6}$ because $n_4^+(c_2)+n_4(c_2)\ge 3$, and $g^*(b_4)\ge \frac{n_4(b_4)}{6}$ because $z_2\in N_4^{-1}(b_4)$ and $n_4(b_4)+n_5(b_4)\ge 3$, which implies that $f_1(z_2)\ge t^*_{c_2}+t^*_{b_4}\ge \frac{1}{3}$ violating to Claim 2.
%We see $a_1\in V_5$, otherwise $a_1\in V_4$ and so $g^*(c_1)\ge\frac{1}{6}n_4(c_1)$ because $z_1\in N_4^+(c_1)$, which implies $f_1(z_1)\ge g^*(z_1)+t_{c_1}^*\ge\frac{1}{6}$ violating to Claim 2.
Note that $a_2\in V_4\cup V_5$. By Ineq.~(\ref{g*}), $t_{b_4}^*\ge\frac{1}{6}$ because $n_4(b_4)+n_5(b_4)\ge3$ and $n_4^{-1}(b_4)+n_5^2(b_4)\ge1$. If $n_5^2(z_1)\neq 0$ or $a_1=c_1$, then $t_{c_2}^*\ge \frac{1}{6}$, which follows that $f_1(z_2)\ge t^*_{c_2}+t^*_{b_4}\ge \frac{1}{3}$ violating Claim 3. If $n_5^2(z_1)=0$ and $a_1\ne c_1$, then $g_1(a_1)=\frac{1}{3}$ as $d(b_4)=3$. That is, $z_1b_4\notin E(G)$.  It follows that $g_2(c_1)=\frac{1}{3}$, or $g_2(c_1)=\frac{1}{6}$ and there exists a vertex $c_1'\in N_5(z_1)\setminus c_1$ such that $g^*(c_1')\ge \frac{1}{6}$. That is, $g^*(z_1)\ge\frac{1}{3}$, or $t_{v}^*\ge\frac{1}{6}$ for any $v\in\{c_1,c_1'\}$. But then $f_1(z_1)\ge g^*(z_1)\ge \frac{1}{3}$ or $f_1(z_1)\ge g^*(z_1)+t^*_{c_1}+t^*_{c_1'}\ge \frac{1}{3}$, which violates Claim 3. This proves $(c)$.\medskip

By $(c)$, $\{b_3,b_4\}\cap\{a_1,a_2,a_3\}=\emptyset$. We first prove  $d(b_4)=2$. Suppose  not.  By Claim 2, $g(z_2)<0$ and so $a_3\in V_5$. By Claim 3, $d(b_4)=3$. Then $a_2\in V_5$, otherwise $a_2\in V_4$ and so $a_1\in V_5$, we see for any $v\in\{a_3,b_4\}$, $t_v^*\ge\frac{1}{6}$ because $d(b_4)=3$, $z_2\in N_4^{-1}(v)$ and $a_2\in N_4^{-1}(a_3)\cup N_4^+(a_3)$, which means $f_1(z_2)\ge g^*(z_2)+t_{a_3}^*+t_{b_4}^*\ge\frac{1}{3}$ violating Claim 2.
Clearly, $a_3\in V_5^1$ and $n_5^2(b_4)=0$, otherwise $t_{v}^*\ge\frac{1}{6}$ for any $v\in\{a_3,b_4\}$, which means $f_1(z_2)\ge\frac{1}{3}$. Moreover, $a_2\in V_5^1$, otherwise  since $b_4\in N_5^2(z_2)$, we see $g^*(a_3)\ge\frac{1}{3}$, or there exists $a_3'\in N_5(z_2)\setminus a_3$ such that $t_{v}^*\ge\frac{1}{6}$ for any $v\in\{a_3,a_3'\}$, which means $f_1(z_2)\ge\frac{1}{3}$.
%which means $t_{a_3}^*\ge\frac{1}{3}$ and so $f_1(z_2)\ge g^*(z_2)+t_{a_3}^*\ge\frac{1}{3}$.
Hence, $a_1\in V_4$. Since $G+a_3b_3\in \mathcal{F}_G$, we see $G$ has a $5$-path $P(a_3b_3)=a_3c_1c_2c_3c_4b_3$. Note that $c_1\in\{a_2,z_2\}$ and $c_4\in\{b_2,b_4\}$ because $d(a_3)=d(b_3)=2$. We assert  $c_4\ne b_4$. Suppose not. Then $z_2=c_i$ for some $i\in[3]$ because $G$ is $C_6$-free. Since $d(a_2)=d(a_3)=2$, we see $z_2=c_1$. By Claim 4, $\{c_2,c_3\}\subseteq V_5$. But then $c_3\in N_5^2(b_4)$, a contradiction, as asserted. Thus, $c_4=b_2$. Then $c_1=z_2$, otherwise $c_1=a_2$, $c_2=a_1$ and $c_3\in V_5$, we see $t_{v}^*\ge\frac{1}{6}$ for any $v\in \{c_3,a_2\}$ because $n_4^+(c_3)+n_4(c_3)\ge3$, which means $f_1(a_1)\ge g^*(a_1)+t_{c_3}^*+t_{a_2}^*\ge\frac{1}{3}$.  Since $G$ is $C_6$-free, $b_4=c_i$ for some $i\in\{2,3,4\}$. Since $d(b_4)=3$ and $b_4\ne c_4$, we see $c_2=b_4$. By Claim 1, $c_3\in V_5^1$ because $n_5^2(b_4)=0$. By Definition~\ref{law1}(1), $g^*(b_4)=\frac{1}{3}$ and so $t_{b_4}^*=\frac{1}{3}$. But then $f_1(z_2)\ge g^*(z_2)+t_{b_4}^*\ge\frac{1}{3}$, a contradiction, as desired. Thus, $d(b_4)=2$.

We next prove $d(b_3)=2$. Suppose not. By Claim 2, $g(v)<0$ for any $v\in\{b_2,z_2\}$, which means $a_3\in V_5$. By Claim 3, $d(b_3)=3$. Then $g_1(b_4)=\frac{1}{3}$, otherwise $g_1(b_4)=\frac{1}{6}$ which implies $g^*(b_3)\ge g_1(b_3)\ge\frac{1}{3}$, that is, $t_{b_3}^*\ge\frac{1}{3}$, and so $f_1(b_2)\ge t_{b_3}^*\ge\frac{1}{3}$. Moreover, $a_3,a_2\in V_5^1$, otherwise $f_1(z_2)\ge\frac{1}{3}$. It follows that $a_1\in V_4$. Then $n_4(b_3)=0$ and $n_5^2(b_3)+n_4^{-1}(b_3)=0$, otherwise by Ineq.~(\ref{g*}),  $t_{b_3}^*\ge\frac{1}{6}$ and so $f_1(b_2)\ge t_{b_3}^*\ge\frac{1}{6}$ violating Claim 2. Since $G+a_2b_4\in\mathcal{F}_G$, we see $G$ has a $5$-path $P(a_2b_4)=a_2c_1c_2c_3c_4b_4$. Then $c_4\in\{z_2,b_3\}$ and $c_1\in\{a_1,a_3\}$ since $a_3,a_2\in V_5^1$. Since $G$ is $C_6$-free, $c_4\ne z_2$ as $a_3\in V_5^1$. Hence, $c_4=b_3$. Since $g_1(b_4)=\frac{1}{3}$ and $g(z_2)<0$, we see $G-b_4$ has no $2$-path, which implies $c_1\ne a_3$. Hence, $c_1=a_1$. By Claims 1 and 2, $c_3=b_2$. By Claim 1 and $g(b_2)<0$, we have $c_2\in V_5$. But then $b_2\in N_4^{-1}(b_3)$, a contradiction, as desired. Thus, $d(b_3)=2$.

%When $b_3\in V_5$, if $g_1(b_4)=\frac{1}{6}$, then $d(b_3)=3$ and $g_1(b_3)=\frac{1}{3}$, which follows that $f_1(b_2)\ge \frac{1}{3}$, a contradiction. If $g_1(b_4)=\frac{1}{3}$, when $a_3\in V_5^2$, then $f_1(z_2)\ge \frac{1}{3}$, a contradiction. When $a_3\in V_5^1$, then $a_2\in V_5^1$, otherwise, $f_1(z_2)\ge \frac{1}{3}$. We see $a_1\in V_4$. Then $d((b_3))=3$, $n_4(b_3)=1$ and $n_5^2(b_3)=0$, else, $f_2(b_2)\ge \frac{1}{6}$ and so $f_3(y_1)\ge 0$, a contradiction. Let $N_5(b_3)=\{b_4, b_3'\}$.  Because $G+a_2b_4 \in \mathcal{F}_G$, then $z_2b_3'\in  E(G)$ or there is a 2-path $P(a_1b_2)$ or there is a 2-path $P(a_1b_3')$. If $z_2b_3'\in E(G)$, then $g_1(b_4)=\frac{1}{6}$ and so $g^*(b_3)\ge \frac{1}{3}$, which follows that $f_1(b_2)\ge \frac{1}{3}$, a contradiction. If the 2-path $P(a_1b_2)$ exists, then $f_3(y_1)\ge 0$, a  contradiction. If the 2-path $P(a_1b_3')$ exists, then $z_1b_3'\in E(G)$ because $d(b_3')=2$ and so $f_1(z_1)\ge \frac{1}{3}$, a contradiction. Thus $d(b_3)=2$.
%W.l.o.g, we assume $g(b_2)<0$.
Since $G+b_4x\in\mathcal{F}_G$, $G$ has a $5$-path $P(xb_4)=xc_1c_2c_3c_4b_4$.  Then $c_4\in\{z_2,b_3\}$ as $d(b_4)=2$. By Claims 1 and 4, $c_4\ne z_2$. Hence, $c_4=b_3$ and so $c_3=b_2$ as $d(b_3)=2$. By Claim 1, $c_2\in V_4$. Then $g(z_2)<0$ and $n_5^2(z_2)=0$, otherwise $t_{b_3}^*\ge\frac{1}{6}$ so $f_1(b_2)\ge t_{b_3}^*\ge\frac{1}{6}$ violating Claim 2. It follows that $a_1,a_2,a_3\in V_5$. Then $g_1(a_3)=\frac{1}{3}$, which implies $t_{b_3}^*\ge\frac{1}{6}$. But then $f_1(b_2)\ge t_{b_3}^*\ge\frac{1}{6}$ violating Claim 2, a contradiction.  \qed\medskip

%$b_4b_3\in E(C_6(b_4x))$, otherwise $b_4z_2\in E(C_6(b_4x))$ and so $G$ has $4$-path with ends $x$ and $z_2$ containing $y_2$, which violates to Claims 1 and 4. Then there is a 4-path $P(xb_3)$. Note that  $d(b_3)=2$, $g(y)<0$ and $d(x, b_3)=3$. Hence $n_3(b_2)\ne 0$, that is $g(b_2)\ge 0$ and  $a_1\in V_5$.  Observe that $n_5^2(z_2)=0$ and $g(z_2)<0$, otherwise, $g^*(b_3)=\frac{1}{6}$ and so $f_1(b_2)\ge \frac{1}{6}$ violating to Claim 2. Thus $a_2\in V_5^2$ and $g_1(a_3)=g_2(b_4)=\frac{1}{3}$, which follows that $f_1(z_2)\ge g^*(z_2)\ge \frac{1}{3}$, a contradiction.

\noindent{\bf Claim 7.} $b_4\ne a_3$ or $b_1\ne a_1$.\medskip

\noindent{\bf Proof.} Suppose not. Then $b_4= a_3$ and $b_1=a_1$. Note that $\{a_1,a_3\}\subseteq V_4\cup V_5$. We first assert $a_2\notin\{b_2,b_3\}$. Suppose not. By symmetry, assume $a_2=b_2$.
Then $a_3\in V_5$, otherwise $\{a_2, b_3\}\subseteq V_5$ and $a_2\in V_5^2$, we see $t_{a_2}^*\ge \frac{1}{3}$ because $n_4(a_2)+n_5(a_2)+n_4^{+}(a_2)\ge4$, which follows $f_1(a_3)\ge g^*(a_3)+t_{a_2}^*\ge\frac{1}{3}$ violating Claim 2. By Claim 2, $g(z_2)<0$. Moreover, $a_2\in V_4$, otherwise $a_2\in V_5$, we see $t_{a_3}^*\ge \frac{1}{3}$ because $n_4(a_3)+n_5(a_3)\ge3$ and $n_4^+(a_3)+n_4^{-1}(a_3)+n_5^2(a_3)\ge 2$, which yields $f_1(z_2)\ge g^*(z_2)+t_{a_3}^*\ge\frac{1}{3}$ violating Claim 3. By Claim 3, $d(a_3)=3$. Clearly, $d(b_3)=2$, otherwise $t_{a_3}^*\ge\frac{1}{3}$ and so $f_1(z_2)\ge t_{a_3}^*\ge\frac{1}{3}$. Since $G+z_2b_3\in \mathcal{F}_G$, $G-b_3-a_3$ has a 3-path with vertices $a_2,a_2',z_2',z_2$ in order. But then $G$ has a copy of $C_6=a_2b_3a_3z_2z_2'a_2'$, a contradiction, as asserted.
%If $a_3, a_2\in V_5$, then $g^*(a_3)\ge \frac{n_4(a_3)}{3}$ becuase $n_5^2(a_3)\ge 2$ or $n_5^2(a_3)\ge 1$ and $n_4(a_3)\ge 2$ with $n_4^+(a_3)+n_4^{-1}(a_3)\ge 1$. If $a_3\in V_4$, then $a_2, b_3\in V_5$ and $a_2\in V_5^2$. We see $g^*(a_2)\ge \frac{n_4(a_2)}{3}$, which follows that $f_1(a_3)\ge \frac{1}{3}$, a contradiction. Then we just consider  the case $a_2\in V_4$ and $a_3\in V_5$.
%We have $b_3\in V_5$ and $d(b_3)=2$ and $d(a_3)=3$, otherwise, $g^*(a_3)\ge \frac{n_4(a_3)}{3}$, which means that $f_1(z_2)\ge \frac{1}{3}$, a contradiction. Becuase $G+z_2b_3\in \mathcal{F}_G$, there is a 3-path $P(a_2z_2)$  in $G-b_3-a_3$. Let $P(a_2z_2)=a_2a_2'z_2'z_2$. But then there is a $C_6=a_2a_2'z_2'z_2a_3b_3a_2$ in $G$, a contradiction. Hence, $a_2\notin\{b_2,b_3\}$.

We next assert $\{a_1,a_3\}\subseteq V_5$.
W.l.o.g., assume $a_3\in V_4$. By Claims 1 and 2, $a_1\in V_4$.
By Claim 1, $\{b_2, b_3\}\subseteq V_5$, otherwise $\{b_2, b_3\} \subseteq V_3$, by Definition \ref{law2}(2), $f_2(z_1)\ge \frac{1}{6}$ and so $t_{z_1}^2\ge\frac{1}{6}$, which yields $f_3(y_1)\ge f_2(y_1)+t_{z_1}^2\ge0$ as $f_2(y_1)\ge g^*(y_1)\ge-\frac{1}{6}$.
We see $a_2\in V_3$ and $d(b_2)=d(b_3)=2$, otherwise $f_1(a_i)\ge \frac{1}{3}$ for some $i \in \{1,3\}$. Then  $G$ has a 3-path $P(a_2a_i)$ for some $i \in \{1,3\}$ since $G+b_2a_3\in \mathcal{F}_G$. W.l.o.g., let $i =1$. Let $P(a_2a_1)=a_2a_2'a_1'a_1$.
Then $a_1'\in V_5$ and so $a_2'\in V_4$. Thus for any $v\in \{b_2, a_1'\}$, $g^*(v)\ge\frac{1}{6}n_4(v)$ as $a_1\in N_4^+(v)$ and so $t_{v}^*\ge \frac{1}{6}$. But then $f_1(a_1)\ge g^*(a_1)+t_{b_2}^*+t_{a_1'}^*\ge\frac{1}{3}$, violating Claim 2, as asserted.
%Furthermore, $a_1\in V_5^2$ because $d(a_1)\ge 3$.  Then $a_2, b_3\in V_5$. Clearly, $g^*(v)\ge \frac{n_4(v)}{6}$ for any $v \in \{a_2, b_3\}$, which follows that $f_1(a_3)\ge \frac{1}{3}$, a contradiction.

By Claim 2, $g(z_i)<0$ since $d(a_1)\ge3$ and $d(a_3)\ge3$. Note that $\{b_2,b_3\}\cap V_4\ne\emptyset$ and $\{b_2,b_3\}\cap V_5\ne\emptyset$, otherwise $t_{a_1}^*\ge \frac{1}{3}$ because $d(a_1)=n_4(a_1)+n_5(a_1)\ge3$ and $n_4^{-1}(a_1)+n_4^+(a_1)+n_5^2(a_1)\ge 2$, which yields $f_1(z_1)\ge t_{a_1}^*\ge\frac{1}{3}$ violating Claim 3. W.l.o.g., assume $b_2\in V_4$ and $b_3\in V_5$.  By Claim 3, $d(a_1)=d(a_3)=3$. Clearly, $d(b_3)=2$, otherwise $t_{a_3}^*\ge \frac{1}{3}$ and so $f_1(z_2)\ge t_{a_3}^*\ge\frac{1}{3}$.
Then $G-a_1$ has a 2-path $P(z_1v)$ for some $v\in\{a_2,b_2\}$, or $G-a_3-b_3$ has a 3-path $P(z_2a_2)$ as $G+b_3a_2\in \mathcal{F}_G$. We assert $G$ has no 3-path $P(z_2a_2)$. Suppose not. Let $P(z_2a_2)=z_2z_2'a_2'a_2$. Then $n_4^1(y_2)=0$, else $g^*(y_2)\ge-\frac{1}{6}$ and so $f_3(y_2)\ge0$ because $t^*_{a_3}\ge \frac{1}{6}$ and so $t^2_{z_2}\ge \frac{1}{6}$. Hence, $z_2'\in V_3\cup V_5$. Moreover, $z_2'\in V_5$, otherwise $z_2'=y_2$,
	$a_2'\in V_4$ and so $a_2\in V_5$, which follows that $f_1(a_2')\ge \frac{1}{3}$ as $t_{a_2}^*\ge \frac{1}{3}$ violating Claim 3. It is easy to verify that $t_{v}^*\ge\frac{1}{6}$ for any $v\in\{a_3,z_2'\}$. But then $f_1(z_2)\ge t_{a_3}^*+t_{z_2'}^*\ge\frac{1}{3}$ violating Claim 3, as asserted. Thus, $G-a_1$ has a 2-path $P(z_1v)$ for some $v\in\{a_2,b_2\}$. Let $P(z_1v)=z_1v'v$. By Claim 4, $v'\notin V_3$. Hence, $v'\in V_5$ as $g(z_1)<0$. Furthermore, one can easily check $t_{u}^*\ge\frac{1}{6}$ for any $u\in\{v',a_1\}$. But then $f_1(z_2)\ge t_{a_1}^*+t_{v'}^*\ge\frac{1}{3}$ violating Claim 3. This completes the proof of Claim 7.\qed\medskip

Note that $\{a_1, a_3\}\subseteq V_4\cup V_5$. By Claim 7 and symmetry, we may assume $b_1\ne a_1$. We first show $g(z_1)<0$. Suppose not.
Note that $a_1\notin V_5$, otherwise $\{a_1,a_2\}\in V_5^1$ and so $a_3\in V_4$ by Claim 2, we see $g^*(a_1)\ge\frac{1}{6}$ and so $t_{a_1}^*\ge\frac{1}{6}$, which yields $f_1(z_1)\ge t_{a_1}^*\ge\frac{1}{6}$ violating Claim 2. By Claim 1, $|\{b_1,a_1\}\cap V_5|\ge 1$. Thus, $b_1\in V_5$. By Claim 2, $b_1\in V_5^1$, $\{b_2,b_4\}\subseteq V_5$ and $b_3\in V_4$. Then $g^*(b_1)\ge\frac{1}{6}$ because $n_5^2(b_3)\ge1$, which means $t_{b_1}^*\ge\frac{1}{6}$. But then $f_1(z_1)\ge t_{b_1}^*\ge\frac{1}{6}$, which contradicts to Claim 2, as desired.
Thus, $\{a_1, b_1\}\subseteq V_5$.
Clearly, $g^*(v)<\frac{n_4(v)}{6}$ for some $v \in \{a_1, b_1\}$, otherwise $f_1(z_1)\ge t_{a_1}^*+t_{b_1}^*\ge\frac{1}{3}$, which contradicts to Claim 2.

Assume first that $g^*(a_1)<\frac{1}{6}n_4(a_1)$. We first assert $d(a_1)=2$.
	Suppose  $d(a_1)\ge 3$. Since $g^*(a_1)<\frac{1}{6}n_4(a_1)$, $n_4(a_1)=1$ and $n_5^2(a_1)+n_4^{-1}(a_1)+n_4^+(a_1)=0$. Hence, $a_2\in V_5^1$, which means $a_3\in V_4$. But then $a_1\in N_5^2(a_2)$ violating Claim 2, as asserted.
We next assert $a_2\in V_5$. Suppose not. Then $a_2\in V_4$.
By Claim 1, $a_3\in V_5$.  Since $g^*(a_1)<\frac{1}{6}n_4(a_1)$, we see $b_1\in V_5^1$, which implies $b_2\in V_5$.  Moreover, $b_2\in V_5^1$, otherwise $g^*(b_1)=\frac{1}{3}$ or $G$ has a vertex $z_1'\in (N_5(z_1)\cap N_5(b_2))\setminus b_1$ such that  $g^*(v)\ge \frac{1}{6}n_4(v)$ for any $v\in \{z_1', b_1\}$, which implies $f_1(z_1)\ge \frac{1}{3}$ violating Claim 3. Hence, $b_3\in V_4$.
By Claim 1, $b_4\in V_5$. Clearly, $b_3\neq a_2$, otherwise $G+b_1a_2\notin \mathcal{F}_G$.
If $b_4\neq a_3$, then $d(b_4)=2$, else $t_u^*\ge\frac{1}{6}$ for any $u \in \{a_3, b_4\}$, which yields $f_1(z_2)\ge t_{a_3}^*+t_{b_4}^*\ge\frac{1}{3}$. Clearly, $d(a_3)=2$, else $t_{a_3}^*\ge\frac{1}{3}$ and so $f_1(z_2)\ge t_{a_3}^*\ge\frac{1}{3}$. Since $G+b_1b_4\in \mathcal{F}_G$, $G-b_4$ has a 2-path $P(b_3z_2)$ or $G$ has a 3-path $P(z_1z_2)$. In both cases, $b_3y_2\in E(G)$ or
there exists $z_2'\in N_5(z_2)\setminus a_3$ such that $t_{v}^*\ge\frac{1}{6}$ for any $v\in\{z_2',a_3\}$. But then $f_3(y_1)+t^2_{b_3}+t^2_{z_2}\ge 0$ because $t^2_{v}\ge\frac{1}{6}$ for any $v \in \{b_3, z_2\}$, or
$f_1(z_2)\ge t_{z_2'}^*+t_{a_3}^*\ge\frac{1}{3}$, a contradiction.
If $b_4=a_3$, then $d(a_3)=3$ by Claim 3. Since $G+a_2b_2\in \mathcal{F}_G$, $G-b_1-b_2$ has a 3-path with vertices $z_1,z_1',a_2',a_2$ in order, or $G-b_4$ has a 2-path
with vertices $z_2u'u$ for some $u\in \{b_3, a_2\}$. For the former, we have $f_1(z_1)\ge t_{z_1'}^*+t_{b_1}^*\ge\frac{1}{3}$ since $t_v^*\ge\frac{1}{6}$ for any $v\in\{z_1',b_1\}$. For the latter, we have $f_1(z_2) \ge t_{b_4}^*\ge\frac{1}{3}$ because $n_4^+(b_4)+n_4^{-1}(b_4)\ge 2$ and $n_4(b_4)=3$ when $u'\in V_5$, or $u'\in N_4^2(z_2)$ when $u'\in V_4$, or $f_3(y_2)\ge f_2(y_2)+t^2_{z_2}+t^2_{u}\ge 0$ since $t^2_{v}\ge \frac{1}{6}$ for any $v \in \{z_2,u\}$ when $u'\in V_3$, a contradiction, as asserted. Thus, $a_2\in V_5$. Then $a_1\in V_5^1$, which yields $n_5^2(z_1)=0$ since $g^*(a_1)<\frac{1}{6}n_4(a_1)$. Hence, $b_1\in V_5^1$ and so $b_2\in V_5$. By Claim 2, $a_3\in V_5$. Then $b_3\in V_4$, otherwise $g_1(v)=\frac{1}{3}$ for any $v\in\{a_1,b_1\}$, which implies $g^*(z_1)\ge \frac{1}{3}$.  Hence, $b_4\in V_5$. But then $f_1(z_2)\ge t_{a_3}^*\ge\frac{1}{3}$, violating Claim 3.

%If $a_3\in V_4$, then $b_4=a_3$, otherwise $b_4\in V_5^1$ and $b_3\in V_5^2$ violating Claim 2.  Thus, $b_3\in V_5$. Clearly, $t_{b_3}^*\ge\frac{1}{3}$ because $b_2\in V_5^2$. But then $f_1(b_4)\ge t_{b_3}^*\ge\frac{1}{3}$, a contradiction.
%If $a_3\in V_5$, then $g_1(a_1)=\frac{1}{3}$, then $g_1(b_1)=g_4(b_1)=\frac{1}{6}$, otherwise, $g^*(z_1)\ge \frac{1}{3}$, a contradiction. Thus $b_2\in V_5^1$, $b_3\in V_4$ and $b_4\in V_5$, which means that $f_1(z_2)\ge \frac{1}{3}$, a contradiction.

Assume next that $g^*(b_1)<\frac{1}{6}n_4(b_1)$.
	We assert  that $d(b_1)=2$. Suppose $d(b_1)\ge3$. Since $g^*(b_1)<\frac{1}{6}n_4(b_1)$, we see $n_4(b_1)=1$ and $n_5^2(b_1)+n_4^{-1}(b_1)+n_4^+(b_1)=0$. Hence,  $a_1,b_2\in V_5^1$, which means $b_3\in V_4$ and $a_2,b_4\in V_5$. By Claim 2, $a_3\in V_5$.
	But then $t_{a_1}^*\ge\frac{1}{3}$ and so $f_1(z_1)\ge t_{a_1}^*\ge\frac{1}{3}$, which  contradicts to Claim  3, as asserted.
Note that $b_2\in V_4$, otherwise $n_5^2(z_1)=0$ and $g(z_1)<0$ as $g^*(b_1)<\frac{1}{6}n_4(b_1)$, which implies $a_2,b_2\in V_5^1$, that is, $a_3,b_3\in V_4$ violating Claims 1 and 2.
By Claim 2, $g(z_1)<0$. Since $g^*(b_1)<\frac{1}{6}n_4(b_1)$, we have $n_4^{+}(b_1)=0$ and $n_4^{-1}(b_1)\le 1$. Hence, $b_3\in V_5$. Note that $b_4\in V_4\cup V_5$. If $b_4\in V_4$, then $a_1\in V_5^1$ as $n_4^{-1}(b_1)\le 1$. Hence, $a_2\in V_5$. By Claims 1 and 2, $a_3=b_4$. Note that $a_2\ne b_3$, otherwise $f_1(b_4)\ge t_{b_3}^*\ge\frac{1}{3}$ violating Claim 2. But then $t_{v}^*\ge\frac{1}{6}$ for any $v\in \{a_2,b_3\}$ and so $f_1(b_4)\ge\frac{1}{3}$, a contradiction. If $b_4\in V_5$, then $g(z_2)<0$, otherwise $f_1(z_2)\ge t_{b_4}^*\ge\frac{1}{6}$ violating Claim 2. Hence, $a_3\in V_5$.  Clearly $a_2\in V_4$, else $t_{a_1}^*\ge\frac{1}{3}$, which implies $f_1(z_1)\ge\frac{1}{3}$ violating Claim 3. Since $n_4^{-1}(b_1)\le 1$, we see $d(b_3)=2$. As $G+z_1b_3\in \mathcal{F}_G$ and $d(b_1)=d(b_3)=2$,  $G$ has a 4-path with vertices $z_1,c_1,c_2,c_3,b_4$ in order. Obviously, $c_3\in V_4$, otherwise $a_3\ne b_4$ and  $t_{v}^*\ge\frac{1}{6}$ for any $v\in\{a_3,b_4\}$, which implies $f_1(z_2)\ge\frac{1}{3}$. Moreover, $g(c_3)<0$, otherwise $c_3\ne z_2$ yields $t_{b_4}^*\ge\frac{1}{3}$ and so $f_1(z_2)\ge\frac{1}{3}$. It follows that $c_1,c_2\in V_5$. Note that $c_1\ne a_1$, otherwise $f_1(a_1)\ge t_{a_1}^*\ge\frac{1}{3}$. It is easy to check that $t_{v}^*\ge\frac{1}{6}$ for any $v\in\{c_1,a_1\}$. But then $f_1(z_1)\ge t_{a_1}^*+t_{c_1}^*\ge\frac{1}{3}$ violating Claim 3.
\medskip

	Finally, let $N_3^{5-}(x)=\{y_1\}$. It remains to prove that $n_4^{5-}(v)=0$ for any $v \in N_3(x)\setminus y_1$. Suppose not. Let $z \in N_4^{5-}(y)$, where $y\in N_3(x)\setminus y_1$.
	By Corollary~\ref{f2g} and Lemmas~\ref{4-0} and \ref{4-5-0}($a,b$), $g^*(z)<0$ and $G$ has a $3$-path $P=zww_2z_2$ such that $d(v)=2$ for any $v\in V(P)$ and $z_2y\notin E(G)$. Let $N_3(z_2)=\{y_2\}$. Note that $g^*(y)\le f_3(y)<\frac{1}{6}$, otherwise $f_4(z)\ge 0$. Then $G$ has a $3$-path $P(xy)$ since $G+y_1z\in \mathcal{F}_G$ and $f_5(y_1)<0$.
	Let $P(xy)=xx_1x_2y$. Then $x_1y_2\in E(G)$ or $G-x_1-y$ has a path of length at most two with ends $x$ and $x_2$ as $G+zx_1\in \mathcal{F}_G$.  Clearly, $x_2\in V_4$, otherwise $g^*(y)\ge\frac{1}{6}$.
Then $x_1\in V_3$.  If $x_1y_2\notin E(G)$, then $G-x_1-y$ has a $2$-path with ends $x$ and $x_2$, that is, $n_3(x_2)\ge3$. If $x_1y_2\in E(G)$, then $n_3^+(x_2)+n_3(x_2)\ge3$. By Lemma~\ref{4-53-0 }, $n_5^-(x_2)=0$. Hence, in both cases,   $f_1(x_2)\ge g^*(x_2)\ge \frac{1}{6}n_3(x_2) +\frac{1}{6}n_4(x_2)$ and so $t_{x_2}^2\ge\frac{1}{6}$. But then $f_3(y)\ge f_2(y)+t_{x_2}^2\ge \frac{1}{6}$, a contradiction.
This completes the proof of Theorem~\ref{-3little}.\qed

%\noindent\textcolor{blue}{{\bf Remark.} By $(\ast\ast)$, Corollary~\ref{f2g} and Lemma~\ref{34-little}, we see  $n_4^{3-}(v)\le n_4^{*-}(v)\le1$. For any $v\in V_3$, $f_5(v)\ge f_3(v)-f_3(v')-\frac{1}{6}|A(v)|$, where $v'\in N_4^{3-}(v)$.}
\section{Proof of Theorem~\ref{f7vge0}}
%\begin{lem}\label{f7vge0}
%For any $2\le i\le5$, $V_i=V_i^{7+}$.\vspace{-12pt}
%\end{lem}
%\pf
By Definition~\ref{law2}(1), for all $x\in V_5$, $f_7(x)\ge0$.  To prove $f_7(x)\ge0$ for all $x\in V_2\cup V_3\cup V_4$, we just need to prove that for any $x\in V_2$, $f_7(x)\ge0$, $n_3^{7+}(x)=n_3(x)$ and $n_4^{7+}(v)=n_4(v)$ for any $v\in N_3(x)$. Let $\alpha_1\in N_2(x)$.
\medskip

\noindent{\bf Case 1.} $n_3^{5-}(x)\ne 0$.
\medskip

Let $y\in N_3^{5-}(x)$. By Theorem~\ref{-3little}, $n_3^{7+}(x)\ge n_3^{5+}(x)\ge n_3(x)-1$ and $n_4^{7+}(v)=n_4^{5+}(v)=n_4(v)$ for any $v\in N_3(x)\setminus y$.
It remains to  prove $f_7(x)\ge0$ and $f_7(y)\ge0$ and $n_4^{7+}(y)=n_4(y)$. If $d(y)=1$, then $V_1=\{\alpha,\alpha_1\}$, which yields $G$ has a $3$-path $P(x\alpha_1)$ because $G+y\alpha\in\mathcal{F}_G$. Let $P(x\alpha_1)=xx_1\alpha_1'\alpha_1$. Then $G-x-\alpha_1'$ has a $t$-path with vertices $\alpha_1,x_t,x_1$ in order for some $t\in[2]$ because $G+y\alpha_1'\in\mathcal{F}_G$. Hence, $G[\{x,x_t,\alpha_1'\}]$ contains at least two edges because $G+\alpha x_1\in\mathcal{F}_G$, which implies $g(x)\ge\frac{2}{3}$ or $n_2^2(x)\ge1$.   Thus, $f_6 (x)\ge g^*(x)\ge\frac{1}{3}$, which means $f_7(x)\ge0$ and $f_7(y)\ge0$. So we next assume $d(y)\ge2$.\medskip

\noindent{\bf Case 1.1 $n_4^{5-}(y)\ne0$. }
\medskip

Let $z\in N_4^{5-}(y)$. By Theorem~\ref{-320}, $n_4^{7+}(y)\ge n_4^{5+}(y)\ge n_4(y)-1$.  So we shall prove $f_7(v)\ge0$ for any $v\in\{x,y,z\}$. Note that  $n_2(x)\ne 0$ as $G+\alpha z\in\mathcal{F}_G$, which yields $f_6(x)\ge g^*(x)\ge\frac{1}{6}$. By  Corollary~\ref{f2g}, Lemma~\ref{4-5-0}($a,b$), $g^*(y)<0$, $g^*(z)<0$ and $G$ has a $3$-path $P=zww_1z_1$ such that $d(v)=2$ for any $v\in V(P)$, and $yz_1\notin E(G)$, where $z_1\in V_4$ and $w,w_1\in V_5$. Let $N_3(z_1)=\{y_1\}$. Clearly, $f_5(y)\ge-\frac{1}{3}$ and $f_5(z)\ge-\frac{1}{6}$. Thus, we need to prove $f_6(x)+\sum_{i\ge1}t_{v_i}^6\ge\frac{1}{2}$ or $\sum_{i\ge1}t_{v_i}^6\ge\frac{1}{3}$, where $v_i\in N_3(x)\setminus \{y\}$.
Since $g^*(y)<0$, we have $y_1x\in E(G)$ and $G$ has a $3$-path $P(yy_1)$ as $\{G+wz_1,G+zz_1\}\subseteq\mathcal{F}_G$. Let $P(y_1y)=y_1a_1a_2y$. Since $g^*(y)<0$, we see $a_2=x$, otherwise $a_2\in V_4$ violating Lemma~\ref{yz16-40}($a$). Thus, $a_1\in V_2\cup V_3$. For any $v\in V_3$, let $A(v)$ be defined as in Definition~\ref{law2}(4), and  $n_4^{3-}(v)\le n_4^{*-}(v)\le1$ by Corollary~\ref{f2g}. Since $d(z_1)=d(w_1)=2$, $f_3(z_1)=-\frac{1}{6}$, which yields $f_5(y_1)\ge f_3(y_1)-(-f_3(z_1))-\frac{1}{6}|A(y_1)|\ge g^*(y_1)-\frac{1}{6}-\frac{1}{6}n_3(y_1)$.

Suppose $a_1\in V_2$. Then $a_1,x\in N_2^+(y_1)$. By Observation~\ref{ob2}(1), $g^*(y_1)\ge\frac{1}{3}n_2(y_1)+\frac{1}{6}n_3(y_1)\ge\frac{1}{6}n_2(y_1)+\frac{1}{6}(n_3(y_1)+1)$ which means $t_{y_1}^6\ge t_{y_1}^5\ge\frac{1}{6}$. If $n_2(y_1)+n_3(y_1)\ge3$, then $t_{y_1}^6 \ge t_{y_1}^5\ge\frac{1}{3}$ as $g^*(y_1)\ge\frac{1}{3}n_2(y_1)+\frac{1}{6}(n_3(y_1)+1)$.
If $g(v) \ge \frac{2}{3}$ for some $v\in\{a_1,x\}$, then $f_6(x)\ge g^*(x)\ge \frac{1}{3}$ and so $f_6(x)+t_{y_1}^6\ge \frac{1}{2}$. So we next assume $n_2(y_1)+n_3(y_1)=2$ and $g(x)=g(a_1)=\frac{1}{6}$. Note that $
N_1(a_1)=\{\alpha_1\}$,  else $\alpha_2a_1\in E(G)$ yields $G$ has a copy of $C_6=\alpha_2\alpha \alpha_1a_2y_1a_1$.
Then $G$ has a $4$-path $P(\alpha_1y_1)$ containing $a_1$ and $x$ or $4$-path $P(\alpha_2y_1)$ because $G+\alpha y_1\in\mathcal{F}_G$. Suppose $P(\alpha_1y_1)$ exists. Then $G-\alpha_1-y_1$ has a $2$-path with vertices $a_1,a_3,x$ in order such that $a_3\in V_3$. Similarly, $t_{a_3}^6\ge\frac{1}{6}$.
Hence, $t^6_{y_1}+t^6_{a_3}\ge \frac{1}{3}$.
Suppose $P(\alpha_2y_1)$ exists. Then $\delta(G)=2$ and let $P(\alpha_2y_1)=\alpha_2b_1b_2b_3y_1$.  Then $b_1\in V_2\setminus\{a_1,x\}$ since $N_1(a_1)=\{\alpha_1\}$ and $g(x)=\frac{1}{6}$. Note that $b_2\notin V_1$, otherwise $b_2=\alpha_1$, which violates the choice of $\alpha$ because $\alpha$ is contained in a $4$-cycle but $z_1$ is a 2-vertex that is not contained in a $4$-cycle. Since $n_2(y_1)+n_3(y_1)=2$, we see $b_3\in V_2\cup V_4$. In fact, $b_3\in V_4$, otherwise $b_3\in V_2$ and so $G$ has a copy of $C_6=b_3b_2b_1\alpha_2\alpha\alpha_1$ because $b_2\notin V_1$. Then $b_2\in V_3$, which implies $n_3(b_3)\ge2$ and $y_1\in N_3^+(b_3)$. Note that $n_5^-(b_3)=0$ as $\delta(G)=2$. By Observation~\ref{ob3}(2), $g^*(b_3)\ge\frac{1}{6}n_3(b_3)+\frac{1}{6}n_4(b_3)$, which implies $t_{b_3}^2\ge\frac{1}{6}$. Hence,  $f_3(y_1)\ge f_2(y_1)+t_{b_3}^2\ge g^*(y_1)+\frac{1}{6}$ and so $t_{y_1}^6\ge\frac{1}{3}$.

Suppose $a_1\in V_3$.
%If $t_{y_1}^5+t_{a_1}^5\ge\frac{1}{3}$, then $f_6(x)+t_{y_1}^6+t_{a_1}^6\ge\frac{1}{2}$ as $t^6_{v}=t^5_{v}$ for any $v \in \{a_1, y_1\}$.
By contradiction, assume $t_{y_1}^5+t_{a_1}^5<\frac{1}{3}$ as $t^6_{v}=t^5_{v}$ for any $v \in \{a_1, y_1\}$.
%, else, we have the desired results.
Then $G$ has a $3$-path $P(y_1x)$ as $G+y_1z\in\mathcal{F}_G$ and $g^*(y)<0$. Let $P(y_1x)=y_1b_1b_2x$. Note that $f_5(a_1)\ge g^*(a_1)-(-f_3(a_1'))-\frac{1}{6}|A(a_1)|\ge g^*(a_1)-\frac{1}{3}-\frac{1}{6}|A(a_1)|$, where $a_1'\in N_4^{3-}(a_1)$. We next prove that\medskip

\noindent{$(a)$ $a_1\notin \{b_1,b_2\}$.}\medskip

To prove $(a)$, suppose first $a_1=b_1$. Then $b_2\in V_2\cup V_3$. In fact $b_2\in V_3$, otherwise $g^*(a_1)\ge\frac{1}{3}n_2(a_1)+\frac{1}{6}(n_3(a_1)+1)$ because $b_2,x\in N_2^+(a_1)$ and $n_2(a_1)+n_3(a_1)\ge3$, which implies $t_{a_1}^5\ge\frac{1}{3}$ as $y_1\notin A(a_1)$.
Moreover, $y_1\in V_3^1$, otherwise for any $v\in\{a_1,y_1\}$,  $g^*(v)\ge\frac{1}{6}n_2(v)+\frac{1}{6}n_3(v)$ since $x\in N_2^+(v)$ and $n_3^2(v)\ge1$
which implies $t_{v}^5\ge\frac{1}{6}$ as $b_2,y_1\notin A(a_1)$ and $a_1\notin A(y_1)$.  Note that $g^*(y_1)\ge\frac{1}{3}$. Clearly, $t_{y_1}^5\ge\frac{1}{6}$ since $a_1\notin A(y_1)$. Then $t_{a_1}^5<\frac{1}{6}$. By Ineq.~(\ref{g*}),
$g^*(a_1)\ge\frac{1}{6}n_2(a_1)+\frac{1}{6}(n_3(a_1)-1)$  since $x\in N_2^+(a_1)$ and $n_3(a_1)\ge 2$. Then $n_4^{3-}(a_1)=1$, otherwise $t_{a_1}^5\ge\frac{1}{6}$ as $b_1,y_1\notin A(a_1)$.
Let $N_4^{3-}(a_1)=\{z_2\}$.
Then $d(z_2)=1$, otherwise $f_3(z_2)\ge-\frac{1}{6}$, which implies $t_{a_1}^5\ge\frac{1}{6}$ as $b_1,y_1\notin A(a_1)$. Then $G-a_1-x$ has a $2$-path with vertices $y_1,b_3,b_2$ in order such that $b_3\in V_4$ as $G+z_2x\in\mathcal{F}_G$ and $y_1\in V_3^1$. Hence, $g^*(b_3)\ge \frac{1}{6}n_3(b_3)+\frac{1}{3}+\frac{1}{6}n_4(b_3)$ since $b_2,y_1\in N_3^+(b_3)$, which yields $t_{b_3}^2\ge\frac{1}{6}$. Hence, $f_3(y_1)\ge f_2(y_1)+t_{b_3}^2\ge\frac{1}{2}$. But then $t_{y_1}^5\ge\frac{1}{3}$ as $a_1\notin A(y_1)$, a contradiction.
%Hence, $t_{v}^6 \ge \frac{1}{6}$ for any $v \in \{y_1, b_2\}$ and so  $f_5(x)+t_{b_2}^6+t_{y_1}^6\ge \frac{1}{2}$, a contradiction.
Suppose now $a_1=b_2$. Then $b_1\in V_2\cup V_3\cup V_4$. In fact $b_1\in V_4$, otherwise  for any $v\in\{a_1,y_1\}$, $g^*(v)\ge\frac{1}{6}n_2(v)+\frac{1}{6}(n_3(v)+1)$ because $x\in N_2^+(v)$ and $n_3^2(v)+n_2(v)\ge3$, which implies $t_{v}^5\ge\frac{1}{6}$ since $y_1\notin A(a_1)$ and $a_1\notin A(y_1)$. For any $v\in\{a_1,y_1\}$, then
$g^*(v)\ge\frac{1}{6}n_2(v)+\frac{1}{6}n_3(v)-\frac{1}{6}$ as $x\in N_2^+(v)$. Then $G-a_1-y_1$ has a $2$-path with ends $x$ and $b_1$ because $G+z_1a_1\in\mathcal{F}_G$. By Ineq.~(\ref{g*}), $g^*(b_1)\ge\frac{1}{3}n_3(b_1)+\frac{1}{6}(n_4(b_1)+2)$ since $a_1,y_1\in N_3^+(b_1)$ and $n_3(b_1)\ge3$, which implies $t_{b_1}^2\ge\frac{1}{3}$ and so $f_3(v)\ge f_2(v)+t_{b_1}^2$.
But then $t_{v}^5\ge\frac{1}{6}$  for any $v \in \{a_1, y_1\}$ as $y_1\notin A(a_1)$ and $a_1\notin A(y_1)$, that is, $t_{y_1}^5+t_{a_1}^5\ge\frac{1}{3}$, a contradiction. This proves $(a)$.\medskip

Then $G-b_2-y_1$ has an $s$-path with ends $x$ and $b_1$ as $G+z_1b_2\in\mathcal{F}_G$ for some $s\in [2]$.  Note that $b_1\notin V_2$, otherwise by Observation~\ref{ob2}(2,4), $g^*(y_1)\ge\frac{1}{3}n_2(y_1)+\frac{1}{6}n_3(y_1)$ because $x\in N_2^+(y_1)$, $n_2(y_1)\ge2$, $b_1\in N_2^+(y_1)\cup N_2^{-1}(y_1)$ and $n_3(y_1)\ge1$, which yields $t_{y_1}^5\ge\frac{1}{3}$ as $a_1\notin A(y_1)$. We assert $b_1\in V_4$. Suppose not. Then $b_1\in V_3$. Hence, $n_2^+(b_1)\ge1$ or $n_2(b_1)+n_3(b_1)\ge4$ which implies $f_3(b_1)\ge g^*(b_1)\ge\frac{1}{3}$ and so $b_1\notin A(y_1)$. Since $x\in N_2^+(a_1)$ and $y_1\in N_3^2(a_1)$, we have $a_1\notin A(y_1)$. Hence, $|A(y_1)|\le n_3(y_1)-2$. By Ineq.~(\ref{g*}),  $g^*(y_1)\ge\frac{1}{3}n_2(y_1)+\frac{1}{6}(n_3(y_1)-1)$ because $x\in N_2^+(y_1)$ and $b_1\in N_3^2(y_1)$.  But then $t_{y_1}^5\ge\frac{1}{3}$, a contradiction, as asserted. Then $y_1\in N_3^+(b_1)$ and $n_3(b_1)\ge3$. Note that when $N_5^-(b_1)=\{w_2\}$, we have  $n_3^+(b_1)\ge2$ as $G+w_2x\in\mathcal{F}_G$. By Ineq.~(\ref{g*}), $g^*(b_1)\ge\frac{1}{3}n_3(b_1)+\frac{1}{6}n_4(b_1)+\frac{1}{3}n_5^-(b_1)$ and so $t_{b_1}^2\ge\frac{1}{3}$. Hence, $f_3(y_1)\ge f_2(y_1)+t_{b_1}^2$. By Ineq.~(\ref{g*}),  $g^*(y_1)\ge\frac{1}{3}n_2(y_1)+\frac{1}{6}(n_3(y_1)-2)$ as $x\in N_2^+(y_1)$. Thus, $a_1\in A(y_1)$, otherwise $t_{y_1}^5\ge\frac{1}{3}$. Hence, $n_4^{3-}(a_1)=1$.
 Let $N_4^{3-}(a_1)=\{z_2\}$.
Clearly, $d(z_2)=1$ and $g(a_1)=\frac{1}{6}$.  Then $G$ has a $3$-path with vertices $a_1,c_1,c_2,x$ in order such that $c_1\in V_4$ as $G+z_2y\in\mathcal{F}_G$. Hence, $a_1\in N_3^+(c_1)$ and $n_3(c_1)\ge3$ since $G+z_2c_2\in\mathcal{F}_G$, which yields $g^*(c_1)\ge\frac{1}{6}n_3(c_1)+\frac{1}{6}n_4(c_1)+\frac{1}{3}$ and so $t_{c_1}^2\ge\frac{1}{6}$. It follows $f_3(a_1)\ge f_2(a_1)+t_{c_1}^2\ge\frac{1}{3}$. But then $a_1\notin A(y_1)$, a contradiction. \medskip

\noindent{\bf Case 1.2} $n_4^{5-}(y)=0$.
\medskip

Note that $n_4^{7+}(y)=n_4^{5+}(y)=n_4(y)$. So we just need to prove $f_7(x)\ge0$ and $f_7(y)\ge0$.
\medskip

\noindent{\bf Case 1.2.1} $n_4^-(y)\ne 0$.
\medskip

Let $z\in N_4^-(y)$. Note that  $n_2(x)\ne 0$ because $G+\alpha z\in\mathcal{F}_G$ and $d(\alpha,z)=4$ and $\alpha_1\alpha_2\notin E(G)$, which yields $f_6(x)\ge\frac{1}{6}$. If $f_5(y)\ge-\frac{1}{6}$, then $f_6(x)+f_6(y)=f_6(x)+f_5(y)\ge0$, which yields $f_7(x)\ge0$ and $f_7(y)\ge0$. So we may assume $f_5(y)<-\frac{1}{6}$. By $(\divideontimes)$, $f_i(y)<-\frac{1}{6}$ for any $i\in[4]$. Then we just need to prove $f_6(x)+\sum_{i\ge1}t_{v_i}^6\ge\frac{1}{3}$, where $v_i\in N_3(x)\setminus y$ and $t_{v_i}^6\ge0$. If $g(x)\ge\frac{2}{3}$ or $n_2^2(x)\ge1$, then $f_6(x)\ge g^*(x)\ge\frac{1}{3}$. So we next assume $g(x)=\frac{1}{6}$ and $n_2^2(x)=0$.  We claim $G-\alpha$ has a $3$-path $P(x\alpha_1)$. Suppose not.
Then $G$ has a $4$-path with ends $\alpha$ and $x$ containing $\alpha_1$ because $G+\alpha y\in\mathcal{F}_G$ and $g^*(y)<0$. It follows $G-x$ has a $3$-path $P(\alpha\alpha_1)$. Let $P(\alpha\alpha_1)=\alpha \alpha_2x_2\alpha_1$,  that is $\alpha \in \mathcal{C}_4$. Then $x_2\in V_2$. Since $G+\alpha_1z\in\mathcal{F}_G$, we see $G$ has a $5$-path with vertices $\alpha_1,x_3,y_3,z_3,w,z$ in order. Then $w\in V_5$, otherwise $w=y$ and $z_3=x$ because $g^*(y)<0$ and $g(z)<0$, which implies $G-\alpha$ has a $3$-path $P(x\alpha_1)$. Hence, $z_3\in V_4$. By Lemma~\ref{55-0}, $n_5^-(z)=0$. Moreover, $d(w)=2$ and $N_5(z)\cap N_5(z_3)=\{w\}$, otherwise
$t_{w}^*\ge\frac{1}{6}$, which yields $f_2(z)\ge g^*(z)+t_{w}^*\ge \frac{1}{6}$ and so $f_3(y)\ge f_2(y)+t_{z}^2\ge-\frac{1}{6}$.
By the choice of $\alpha$, $w \in \mathcal{C}_4$ and $|N(z)\cap N(z_3)|\ge2$. Let $w_1\in N(z)\cap N(z_3)$ and $w_1\ne w$. Then $w_1=y$. Since $G+wx\in\mathcal{F}_G$ and $g^*(y)<0$, $G$ has a $3$-path with ends $y$ and $v\in\{z,z_3\}$ containing $v'\in\{z,z_3\}$ such that $v\ne v'$. But then $|N_5(z)\cap N_5(z_3)|\ge2$ as $g(z)<0$, a contradiction, as claimed.
Let $P(x\alpha_1)=xy_1x_1\alpha_1$. Then $y_1\in V_3$ and $x_1\in V_2$ as $g(x)=\frac{1}{6}$ and $n_2^2(x)=0$. If $n_4^{3-}(y_1)=0$ or $g(x_1)\ge0$,  then $f_3(y_1)\ge g^*(y_1)\ge\frac{1}{6}n_2(y_1)+\frac{1}{6}n_3(y_1)+\frac{1}{3}(n_2^+(y_1)-1)$ since  $x\in N_2^+(y_1)$ and $n_2(y_1)\ge2$, which implies $t_{y_1}^6\ge\frac{1}{6}$ and so $f_6(x)+t_{y_1}^6\ge\frac{1}{3}$. So we next assume $n_4^{3-}(y_1)=1$ and $g(x_1)<0$. Let $N_4^{3-}(y_1)=\{z_1\}$. Then $G$ has a $4$-path with ends $\alpha_1$ and $y_1$ containing $x$ and $x_1$ since $G+z_1\alpha_1\in\mathcal{F}_G$ and $d(z_1,\alpha_1)=3$. Then $G-\alpha_1-y_1$ has a $2$-path  with vertices $x,y_2,x_1$ in order. Then $y_2\in V_3$ as $g(x_1)<0$. Hence, for any $i\in[2]$, $x\in N_2^+(y_i)$ and $x_1\in N_2^{-1}(y_i)$, which yields $f_3(y_i)\ge g^*(y_i)\ge\frac{1}{12}n_2(y_i)+\frac{1}{6}(n_3(y_i)+2)$. Thus, $t_{y_i}^6\ge \frac{1}{12}$ and so $f_6(x)+t_{y_1}^6+t_{y_2}^6\ge \frac{1}{3}$. \medskip

\noindent{\bf Case 1.2.2} $n_4^-(y)=0$.
\medskip

By Observation~\ref{ob}(1),  $N_4(y)=N_4^1(y)$ and $d(y)=2$ which yields no $4$-cycle contains $y$. By the choice of $\alpha$, $\alpha \notin \mathcal{C}_4$.   Clearly, $f_6(y)\ge-\frac{1}{6}$. If $g(x)\ge\frac{1}{6}$, then $f_6(x)+f_6(y)\ge0$, which yields $f_7(x)\ge0$ and $f_7(y)\ge0$. So we next assume $g(x)<0$.
Let $N_4(y)=\{z\}$ and $N_4(z)=\{z_1\}$. By Observation~\ref{ob}(2), $g(z_1)=\frac{1}{6}$. Let $N_3(z_1)=\{y_3\}$. Then $G$ has a $3$-path $P(x\alpha_1)$ because
$G+\alpha y\in\mathcal{F}_G$, $g(x)<0$ and no $4$-cycle contains $\alpha$.
Let $P(x\alpha_1)=xy_1x_1\alpha_1$. Then $y_1\in V_3$ and $x_1\in V_2$, which means $f_5(x)\ge g^*(x)\ge0$. So we just need to prove $\sum_{i\ge1}t_{v_i}^6\ge\frac{1}{6}$, where $v_i\in N_3(x)\setminus y$.
Note that if $n_{4}^{3-}(y_1)=1$, then $G-\alpha_1-y_1$ has a $2$-path with vertices $x,y_2,x_1$ in order such that $y_2\in V_4$ as $G+z_2\alpha_1\in\mathcal{F}_G$, where $N_{4}^{3-}(y_1)=\{z_2\}$.
%In this case, $g(x_1)\ge 0$ or $n_3(y_1)\ge 1$ because $G+z_2\alpha\in\mathcal{F}_G$ and $d(\alpha,z_2)=4$. Then for $v \in \{y_1, y_2\}$, we see $g^*(v)\ge \frac{n_2(v)}{12}+\frac{n_3(v)}{6}+\frac{1}{3}$ because $n_2^{-1}(v)+n_2(v)\ge 3$  and $n_2^+(v)\ge 1$, or $n_2^{-1}(v)+n_2(v)\ge 4$ and $n_3(v)\ge 1$ if $n_4^-(v)\ge 1$;  $g^*(v)\ge \frac{n_2(v)}{12}+\frac{n_3(v)}{6}$ because $n_2^{-1}(v)+n_2(v)+n_2^+(v)\ge 4$ if $n_4^-(v)=0$. Hence, we see $t_{v}^6\ge \frac{1}{12}$ for any $v \in \{y_1,y_2\}$ and so $t_{y_1}^6+t_{y_2}^6\ge \frac{1}{6}$.

Assume $g(x_1)\ge 0$. If $n_4^{3-}(y_1)=0$, then $g^*(y_1)\ge\frac{1}{6}n_2(y_1)+\frac{1}{6}n_3(y_1)$ as $n_2^+(y_1)+n_2(y_1)\ge3$, which yields $t_{y_1}^6\ge\frac{1}{6}$. If $n_4^{3-}(y_1)=1$, then for any $i\in[2]$, $g^*(y_i)\ge\frac{1}{12}n_2(y_i)+\frac{1}{6}(n_3(y_i)+2)$ since $n_2^+(y_i)+n_2(y_i)\ge3$ and $x\in N_2^{-1}(y_i)$, which means $t_{y_i}^6\ge\frac{1}{12}$ and so $t_{y_1}^6+t_{y_2}^6\ge\frac{1}{6}$.

So we further assume $g(x_1)<0$.
Suppose $n_4^{3-}(y_1)=1$. Then $x,x_1\in N_2^{-1}(y_i)$ for any $i\in[2]$.
Note that $n_2(y_1)+n_3(y_1)\ge3$ as $G+z_2\alpha\in\mathcal{F}_G$ and $d(\alpha,z_2)=4$. Similarly, $n_2(y_2)+n_3(y_2)\ge3$ when $n_4^{3-}(y_2)\ne0$. Hence, $g^*(y_i)\ge\frac{1}{12}n_2(y_i)+\frac{1}{6}n_3(y_i)+\frac{1}{3}n_4^{3-}(y_i)$, which implies $t_{y_i}^6\ge\frac{1}{12}$ and so $t_{y_1}^6+t_{y_2}^6\ge\frac{1}{6}$.
 Suppose $n_4^{3-}(y_1)=0$. If $\delta(G)=1$, then $x,x_1\in N_2^{-1}(y_1)$ as $G+\alpha y_1\in\mathcal{F}_G$. Hence, $g^*(y_1)\ge\frac{1}{6}n_2(y_1)+\frac{1}{6}n_3(y_1)$, which yields $t_{y_1}^6\ge\frac{1}{6}$. If $\delta(G)=2$, then $G$ has a $4$-path $P(zx_1)$ or $G-x-x_1$ has a $2$-path with ends $\alpha_1$ and $y_1$ as $G+yx_1\in\mathcal{F}_G$. We claim the later holds. Suppose not. Let $P(zx_1)=za_1a_2a_3x_1$. Then $a_1\in V_4\cup V_5$. If $a_1\in V_4$, then $a_1=z_1$, $a_2=y_3$ and $a_3\in V_3$ as $g(x_1)<0$. Note that $n_5^-(z_1)=0$ as $\delta(G)=2$. Hence, $f_1(z_1)\ge g^*(z_1)\ge\frac{1}{6}$ since $y_3\in N_3^+(z_1)$, which means  $z_1\in N_4^{1+}(z)$. Thus, $t_{z}^2\ge\frac{1}{6}$. If $a_1\in V_5$, then $a_2\in V_4$, which implies $g^*(a_1)\ge\frac{1}{6}n_4(a_1)$ since $z\in N_4^+(a_1)$ and $n_4(a_1)\ge2$. Thus, $t_{z}^2\ge\frac{1}{6}$ as $z_1\in N_4^{1+}(z)$. But then in both cases, $f_3(y)\ge f_2(y)+t_{z}^2\ge\frac{1}{6}\ge0$, a contradiction, as claimed. Thus, $n_2(y_1)\ge3$, which implies $g^*(y_1)\ge\frac{1}{6}n_2(y_1)+\frac{1}{6}n_3(y_1)$. Hence, $t_{y_1}^6\ge\frac{1}{6}$. \medskip

\noindent{\bf Case 2.} $n_3^{5-}(x)=0$.
\medskip

Note that $n_3^{7+}(x)=n_3^{5+}(x)=n_3(x)$. By Theorem~\ref{-320}, $\sum_{v\in N_3(x)}n_4^{5-}(v)\le1$.\medskip

\noindent{\bf Case 2.1} $\sum_{v\in N_3(x)}n_4^{5-}(v)=1$.\medskip

Let $z\in N_4^{5-}(y)$, where $y\in N_3(x)$.   So we just need to prove that $f_7(x)\ge0$ and $f_7(z)\ge0$. By Corollary~\ref{f2g}, Lemmas~\ref{4-0} and \ref{4-5-0}($a,b$), $g^*(z)<0$ and $G$ has a $3$-path $P=zww_1z_1$ such that $d(v)=2$ for any $v\in V(P)$, and $yz_1\notin E(G)$, where $z_1\in V_4$ and $w,w_1\in V_5$. Let $N_3(z_1)=\{y_1\}$. Note that $f_6(z)\ge-\frac{1}{6}$. Then $f_3(y)<\frac{1}{6}$, else $f_5(z)\ge 0$. We claim $g^*(x)\ge\frac{1}{6}$. Suppose not. Then $n_2(x)=0$, otherwise $g^*(x)\ge\frac{1}{6}$. Since $z \notin \mathcal{C}_4$, we see
	 $n_1(x)=1$ by the choice of $\alpha$, $\alpha \notin \mathcal{C}_4$. It follows $g(x)<0$. Clearly, $n_3(y)\ne0$  and $G$ has a $4$-path with ends $y$ and $\alpha_1$ containing $x$ because $\{G+\alpha z,G+\alpha_1z\}\subseteq\mathcal{F}_G$, $g(x)<0$  and $g^*(z)<0$. Note that $G$ has no $3$-path $P(x\alpha_1)$, otherwise let $P(x\alpha_1)=xy_2x_1\alpha_1$, we see $x_1\in V_2$ and $y_2\in V_3$, which yields $f_3(y)\ge g^*(y)\ge\frac{1}{6}$ as $y_2\in N_3^2(x)$. Hence, $G-\alpha_1$ has a $3$-path $P(xy)$. Let $P(xy)=xb_1b_2y$. Then $b_1\in V_3$. Note that $b_2\in V_4$, otherwise $f_3(y)\ge g^*(y)\ge\frac{1}{6}$.  By Lemma~\ref{4-53-0 }, $n_5^-(b_2)=0$. Then $g^*(b_2)\ge \frac{1}{6}n_3(b_2)+\frac{1}{6}n_4(b_2)$ as $y\in  N_3^+(b_2)$ and $n_3(b_2)\ge2$, which yields $t_{b_2}^2\ge\frac{1}{6}$ and so $f_3(y)\ge f_2(y)+t_{b_2}^2\ge\frac{1}{6}$, a contradiction.  Then $f_6(x)+f_6(z)\ge0$, which yields $f_7(x)\ge0$ and $f_7(z)\ge0$. \medskip

%Then $g(b_1)<0$, otherwise $g^*(x)\ge0$ and $g^*(b_2)\ge \frac{1}{3}n_3(b_2)+\frac{1}{6}n_4(b_2)$ because $\{y,b_1\}\subseteq N_3^+(b_2)$, which yields $t_{b_2}^2\ge\frac{1}{3}$ and so $f_3(y)\ge g^*(y)+t_{b_2}^2\ge\frac{1}{3}$. Clearly, $d(v)\ge2$ for any $v\in N_4^{3-}(b_1)$, otherwise let $d(z')=1$ for some $z'\in N_4^{3-}(b_1)$, we see $G+z'\alpha\notin \mathcal{F}_G$ because $g(b_1)<0$, $g(x)<0$ and $d(z',\alpha)=4$. Hence, $f_3(v)\ge-\frac{1}{6}$ for any $v\in N_4^{3-}(b_1)$. By Lemma~\ref{path joint}, $G-y-b_1$ has a $2$-path with ends $x$ and $b_2$. Hence, $g^*(b_2)\ge \frac{1}{3}n_3(b_2)+\frac{1}{6}n_4(b_2)$ because $y\subseteq N_3^+(b_2)$ and $n_3(b_2)\ge3$, which yields $t_{b_2}^2\ge\frac{1}{3}$ and so $f_3(v)\ge g^*(y)+t_{b_2}^2\ge\frac{1}{3}$ for any $v\in\{y,b_1\}$. Hence, $f_5(z)\ge0$ and $f_5(b_1)\ge\frac{1}{6}$. Note that $f_5(x)\ge g^*(x)\ge-\frac{1}{6}$ because $n_3(y)\ne0$. Thus, $f_7(z)\ge0$ and $f_7(x)\ge0$.  \medskip

\noindent{\bf Case 2.2} $\sum_{v\in N_3(x)}n_4^{5-}(v)=0$.\medskip

Note that $n_4^{7+}(v)=n_4^{5+}(v)=n_4(v)$ for any $v\in N_3(x)$. So we shall prove $f_7(x)\ge0$.
We assume $f_7(x)<0$, then
  $g^*(x)<0$.  Since $G+\alpha x\in\mathcal{F}_G$, $G$ has a 3-path $P(\alpha_1\alpha_2)$ where $N(\alpha)=\{\alpha_1, \alpha_2\}$ or a 3-path $P(xx_1)$ for some $x_1\in V_2$. Suppose $P(\alpha_1\alpha_2)$ exists, that is $\alpha \in \mathcal{C}_5$. Let $P(\alpha_1\alpha_2)=\alpha_1b_1b_2\alpha_2$ where $b_1, b_2\in V_2$. In this case,  $\delta(G)=2$ and $v \in \mathcal{C}_4\cup \mathcal{C}_5$
  for any vertex $v \in V(G)$ with $d(v)=2$
by the choice of $\alpha$. By Lemma \ref{4-5-0},  $V_4^{*-}=V_5^{*-}=\emptyset$. Since $G+x \alpha_2\in \mathcal{F}_G$, $G$ has a 5-path $P(x \alpha_2)$. Let $P(x \alpha_2)=x c_1 c_2 c_3 c_4 \alpha_2$. We have $c_1\in V_3$, $c_4\in V_2$, $c_2\in V_3\cup V_4$ and $c_3\in V_3\cup V_2$. When $c_2\in V_3$, then $c_3\in V_2$ as $g^*(x)<0$ and so $g(c_i)= \frac{1}{6}$ for any $i \in [2]$ and $g(c_3)\ge \frac{1}{6}$. It follows that $f_5(c_3)\ge 0$ and $f_5(c_2)\ge \frac{1}{6}$, which yields $f_7(x)\ge 0$ by Definition \ref{law2}(6). When $c_2\in V_4$, then $c_3\in V_3$.
 Observe that $g(c_1)<0$, otherwise $g^*(x)\ge -\frac{1}{6}$ and $t_{c_2}^2\ge t_{c_2}^*\ge \frac{1}{6}$  and so $t_{c_1}^6\ge \frac{1}{6}$
 yielding
 $f_7(x)\ge 0$. It follows that $d(x)\ge 3$ otherwise $x \in \mathcal{C}_4\cup \mathcal{C}_5$. We assert $d(c_1)\ge 3$. Suppose not. Then $c_1\notin \mathcal{C}_4$, otherwise there is another vertex $c_0\in (N(x)\cap N(c_2))\setminus\{c_1, c_3\}$ and so $t_{c_2}^2\ge t_{c_2}^*\ge \frac{1}{6}$, $t_{c_i}^6\ge \frac{1}{6}$ for $i \in \{0,1\}$ yielding $f_7(x)\ge 0$. Then $N_2(\alpha_1)\cap N_2(\alpha_2)= \emptyset$,
  $c_1\in \mathcal{C}_5$ and $G-c_1$ has a 3-path $P(xc_2)$.
Let $P(xc_2)=xc_5c_6c_2$ where  $c_5\in V_3$. Then   $c_6\in V_4$ otherwise $g^*(x)\ge -\frac{1}{6}$ and $t^*_{c_2}\ge \frac{1}{6}$ and so $t_{c_1}^6\ge \frac{1}{6}$ yielding $f_7(x)\ge 0$.
Because $G+c_5\alpha\in \mathcal{F}_G$ and $n_2(x)=0$, there is a 5-path $P(\alpha c_5)$. Let $P(\alpha c_5)=\alpha w_1w_2w_3w_4c_5$ where $w_1\in V_1$ and $w_2\in V_2$.
 When  $w_4\in V_3$, then
 $g(w_4)= \frac{1}{6}$, $w_3\in V_2$ and $g(w_3) \ge 0$, which follows that $f_6(w_4)\ge \frac{1}{6}$ and so $f_7(x)\ge 0$ by Definition \ref{law2}(6).
When $w_4\in V_4$, we see  $w_3\in V_3$ and
$n_3(w_4)\ge 2$, which follows that $t^*_{c_6}\ge \frac{1}{6}$ and $t^*_{c_2}\ge \frac{1}{6}$ yielding $t_{c_i}^6\ge \frac{1}{6}$ for $i \in \{1,5\}$,  which follows that
 $f_7(x)\ge 0$. Thus $d(c_1)\ge 3$. There is another vertex $c_{11}\in N_4(c_1)\setminus\{c_2\}$ and $n_4(c_{11})\ge 1$ because $G+c_{11}\alpha \in \mathcal{F}_G$ and $n_3(c_1)=n_2(x)=0$. We see $t^2_{c_{11}}\ge \frac{1}{6}$.
  When $N_2(\alpha_1)\cap N_2(\alpha_2)\neq  \emptyset$, then  $d(c_2)\ge 3$, otherwise $c_2\in \mathcal{C}_4$ and there is a vertex $v \in (N_4(c_1)\cap N_4(c_3))\setminus c_2$, which follows that $t_{c_1}^6\ge \frac{1}{3}$ yielding $f_7(x)\ge 0$.
    We  claim  $f_2(c_2)\ge \frac{1}{6}n_3(c_2)$.  Suppose not.
 %If $n_3(c_2)+n_4(c_2)\ge 3$, then $f_2(c_2)\ge \frac{1}{6}n_3(c_2)$.
 Then  $n_3(c_2)+n_4(c_2)=2$ and  there is a vertex $c_{21}\in N_5(c_2)$.  When $d(c_{21})\ge 3$ or $g(c_{21})=\frac{1}{6}$, we see $f_2(c_2)\ge \frac{1}{6}n_3(c_2)$. Then we just consider the case $N(c_{21})= N_4(c_{21})=\{v_1, c_2\}$. By the choice of $\alpha$, we see $c_{21}\in \mathcal{C}_4\cap \mathcal{C}_5$.
   If there is a vertex $c_{22}\in (N_5(c_2)\cap N_5(v_1))\setminus c_{21}$ or $g(v_1)\ge 0$, then $f_2(c_2)\ge \frac{1}{6}n_3(c_2)$. Thus $v_1c_3\in E(G)$ otherwise $c_1v_1\in E(G)$ but then  $G+v_1\alpha\notin \mathcal{F}_G$. Note that there is a 3-path $P(c_2v_1)$ with $g(c_2)=\frac{2}{3}$ and $g(v_1)<0$. Let $P(c_2v_1)=c_2u_1u_2v_1$ where $u_i \in V_5$ for $i \in [2]$, which follows $f_2(c_2)\ge \frac{1}{6}n_3(c_2)$, as claimed. Note that  $t^2_{c_{11}}\ge \frac{1}{6}$, we see $t^6_{c_1}\ge \frac{1}{3}$ and so $f_7(x)\ge 0$.
  When $N_2(\alpha_1)\cap N_2(\alpha_2)= \emptyset$,  we see there is a vertex $d_2\in N_4(d_1)$ with $n_3(d_2)\ge 2$ because $G+d_1 \alpha\in \mathcal{F}_G$ for $d_1\in N_3(x)\setminus c_1$.  Similarly we have $g(d_1)<0$ and there is a vertex $d_{11}\in N_4(d_1)\setminus d_2$ with $t^2_{d_{11}}\ge \frac{1}{6}$, together with $t^2_{c_{11}}\ge \frac{1}{6}$ we see $f_7(x)\ge 0$.

Suppose $P(xx_1)$ exists and let $P(xx_1)=xy_1y_2x_1$. By Observation~\ref{ob}(1,2), $y_i\in V_3^1$ for any $i\in[2]$. Hence, $f_6(x)\ge g^*(x)\ge-\frac{1}{6}$. So we just need to prove $\sum_{i\ge1}t_{v_i}^6\ge\frac{1}{6}$, where $v_i\in N_3(x)$.
We shall proceed it by contradiction. Then, $f_6(y_1)<\frac{1}{6}$ and so $f_5(y_1)< \frac{1}{6}$.

Note that $g^*(x)<0$, we see $G$ has no $3$-path with ends $x$ and $\alpha_1$.
We assert that  $G-\alpha_1$ has no $3$-path with ends $x$ and $y_1$. Suppose not.
%Then $x\in V(P(y_1\alpha_1))$ which yields $G-\alpha_1$ has a $3$-path $P(xy_1)$ because $g^*(x)<0$.
Let $P(xy_1)=xy_3z_1y_1$. Then $y_3\in V_3$. By Observation~\ref{ob}(1), $g(y_3)<0$ and so $z_1\in V_4$.
Note that $n_4^{3-}(y_3)=0$, otherwise let $z_3\in N_4^{3-}(y_3)$, $G+z_3\alpha\in\mathcal{F}_G$ yields $n_3(y_3)+n_2(x)\ne0$ and so $g(y_3)\ge0$ or $g(x)\ge0$.
Furthermore, $n_5^-(z_1)\ne0$, otherwise $f_2(z_1)\ge g^*(z_1)-\frac{1}{6}n_4(z_1)\ge\frac{1}{6}n_3(z_1)$ as $y_1\in N_3^+(z_1)$ and $n_3(z_1)\ge2$, which implies $f_6(y_3)\ge f_3(y_3)\ge f_2(y_3)+t_{z_1}^2\ge\frac{1}{6}$ and so $t_{y_3}^6\ge\frac{1}{6}$.
Let $w_1\in N_5^-(z_1)$. Then $G-x-z_1$ has a $2$-path with vertices $y_1,z_2,y_3$ in order such that $z_2\in V_4$ because $G+w_1x\in\mathcal{F}_G$. Hence, for any $i\in[2]$, $f_2(z_i)\ge g^*(z_i)-\frac{1}{3}n_5^-(z_i)-\frac{1}{6}n_4(z_i)\ge\frac{1}{12}n_3(z_i)$ because $y_1\in N_3^+(z_i)$, $y_3\in N_3^{-1}(z_i)$ and $n_5^-(z_i)\le1$, which means $t_{z_1}^2\ge\frac{1}{12}$. But then $f_6(y_3)\ge f_3(y_3)\ge g^*(y_3)+t_{z_1}^2+t_{z_2}^2\ge\frac{1}{6}$ and so $t_{y_3}^6\ge\frac{1}{6}$, a contradiction, as asserted. Thus, $G$ has no $4$-path with ends $y_1$ and $\alpha_1$.
%We first assume that $P(xy_1)$ exists.
% Let $P(xy_1)=xy_3z_1y_1$. Then $z_1\in V_4$ and $y_3\in V_3$. Note that $g(v)<0$ for any $v \in N_3(x)\setminus\{y_1\}$. Note that $N_4^{3-}(y)=\{f_3<0\}$. Hence, $N_5^-(z_1)\cup N_4^{3-}(y_3)\ne \emptyset$, else $f_1(z_1)\ge g^*(z_1)\ge \frac{1}{6}|N_3(z_1)|+\frac{1}{6}|N_4(z_1)|$ which means $f_5(y_3)\ge f_3(y_3)\ge\frac{1}{6}$ and so $f_7(x)\ge0$.
%We claim $N_4^{3-}(y_3)=\emptyset$. Suppose not.
%Let $z_3\in N_4^{3-}(y_3)$. Because $G+\alpha z_3$ has a copy of $C_6$, then $N_3(y_3)\neq \emptyset$ or $N_2(x)\neq \emptyset$, a contradiction.
%Thus, $N_4^{3-}(y_3)=\emptyset$. Then $N_5^-(z_1)\ne \emptyset$.
%Let $w_1\in N_5^-(z_1)$. Then $G-x-z_1$ has a $2$-path $P(y_1y_3)$ because $G+w_1x$ has a copy of $C_6$. Let $P(y_1y_3)=y_1z_2y_3$. Then $z_2\in V_4$. Hence, $g^*(z_i)\ge\frac{1}{12}|N_3(z_i)|+\frac{1}{6}|N_4(z_i)|+\frac{1}{3}$ for any $i\in[2]$.  Then $f_5(y_3)\ge f_3(y_3)\ge\frac{1}{6}$ and so $f_7(x)\ge0$.

Since $G+ y_1\alpha\in\mathcal{F}_G$, $G$ has a $4$-path $P(y_1\alpha_2)$.  Let $P(y_1\alpha_2)=y_1a_1a_2a_3\alpha_2$. Note that $\delta(G)=2$. We claim $n_4^{3-}(v)=0$ for any $v\in N_3(x)$. Suppose not. Let $z_1\in N_4^{3-}(y)$ for some $y\in N_3(x)$. Then $y=y_1$, else $G+\alpha z_1\notin \mathcal{F}_G$ because  $d(z_1,\alpha)=4$, $g(y)<0$ and $g(x)<0$. By  Corollary~\ref{f2g} and Lemma~\ref{4-5-0}($a,b$), $G$ has a $3$-path consisting of vertices of degree $2$ with one end $z_1$. But then $G$ has a $4$-path with ends $\alpha_1$ and $y_1$ containing $x$ because $G+z_1\alpha_1\in\mathcal{F}_G$, a contradiction, as claimed.
%But then $G$ has a $3$-path with ends $x$ and $y_1$ because $g^*(x)<0$, a contradiction.
Note that $a_1\in V_2\cup V_3\cup V_4$.

Assume $a_1\in V_4$. Then $a_2\in V_3$ and $a_3\in V_2$. Then $f_2(a_1)\ge g^*(a_1)-\frac{1}{6}n_3(a_1)\ge\frac{1}{6}n_4(a_1)$ since $y_1\in N_3^+(a_1)$ and $n_3(a_1)\ge2$, which means $t_{a_1}^2\ge\frac{1}{6}$. Hence, $f_3(y_1)\ge f_2(y_1)+t_{a_1}^2\ge\frac{1}{6}$. Then $n_4^{3-}(y_2)\ne0$, otherwise $y_2\notin A(y_1)$ yields $f_6(y_1)\ge f_3(y_1)\ge\frac{1}{6}$ and so $t_{y_1}^6\ge\frac{1}{6}$. Let $z_2\in N_4^{3-}(y_2)$. Similarly, $G$ has a $3$-path $P=z_2w_2w_3z_3$ such that $d(v)=2$ for any $v\in V(P)$ and $z_3y_2\notin E(G)$, where $w_2,w_3\in V_5$ and $z_3\in V_4$.
Since $G+\alpha_iz_2\in \mathcal{F}_G$ for some $i \in [2]$ such that  $\alpha_i\in N_1(x_1)$, $G$ has a $4$-path with ends $y_2$ and $\alpha_i$ containing $x_1$. Hence, $g(x_1)>0$ or $n_3^2(x_1)\ge1$ or $G$ has a 3-path $P(x_1y_2)$.
 If  $g(x_1)>0$ or $n_3^2(x_1)\ge1$, then  $f_3(y_2)\ge g^*(y_2)\ge\frac{1}{6}$. If 3-path $P(x_1y_2)$ exists, then let $P(x_1y_2)=x_1x_1'y_2'y_2$. Since $y_2\in V_3^1$, we see $y_2'\in V_4$ and $x_1'\in V_3$.  Hence, $f_2(y_2')\ge \frac{1}{6}n_3(y_2')$ as $y_2\in N_3^+(y_2')$, $n_3(y_2')\ge 2$ and $n_5^-(y_2')=0$, which yields $t_{y_2'}^2\ge\frac{1}{6}$ and so $f_3(y_2)\ge f_2(y_2)+ t_{y_2'}^2\ge\frac{1}{6}$.
 Thus, in both cases, $y_2\notin A(y_1)$ as $f_3(z_2)\ge-\frac{1}{6}$ and $f_3(y_2)+f_3(z_2)\ge0$. But then $f_6(y_1)\ge f_3(y_1)\ge\frac{1}{6}$ and so $t_{y_1}^6\ge\frac{1}{6}$, a contradiction.
%there is a vertex $y_2'\in N_4(y_2)$ with $n_3(y_2')\ge 2$, then
%But then $f_5(y_1)\ge f_3(y_1)\ge f_2(y_1)+t_{a_1}^2\ge\frac{1}{6}$ and so $t_{y_1}^5\ge\frac{1}{6}$, a contradiction.

Assume $a_1\in V_3$. Then $a_1=y_2$, $a_2=x_1$ and $a_3\in V_2$ as $y_1,y_2\in V_3^1$.  Then $f_5(y_1)\ge g^*(y_1)\ge0$ and $g^*(y_2)\ge\frac{1}{6}$ as $g(x_1)>0$. Note that $n_4^{3-}(y_2)\ne0$, otherwise $f_5(y_2)\ge g^*(y_2)\ge\frac{1}{6}$, which means $f_6(y_1)\ge f_5(y_1)+\frac{1}{6}\ge\frac{1}{6}$ by Definition~\ref{law2}(6).  Let $z_2\in N_4^{3-}(y_2)$. Similarly, $G$ has a $3$-path $P=z_2w_2w_3z_3$ such that $d(v)=2$ for any $v\in V(P)$ and $z_3y_2\notin E(G)$, where $w_2,w_3\in V_5$ and $z_3\in V_4$.
%Similarly, $G$ has a $3$-path consisting of vertices of degree two $z_2,w_2,w_3,z_3$ in order and $z_3y_2\notin E(G)$, where $w_2,w_3\in V_5$ and $z_3\in V_4$.
Let $y_3\in N_3(z_3)$.
We claim $y_3x\notin E(G)$. Suppose not. Then $G$ has a $2$-path with vertices $y_2,z_4,y_3$ in order such that $z_4\in V_4$ because $G+z_2w_3\in\mathcal{F}_G$. Then $f_2(z_4)\ge g^*(z_4)-\frac{1}{6}n_4(z_4)\ge\frac{1}{6}n_3(z_4)$ because $y_2\in N_3^{+}(z_4)$ and $n_3(z_4)\ge2$, which means $t_{z_4}^2\ge\frac{1}{6}$. Hence, $f_3(y_2)\ge f_2(y_2)+t_{z_4}^2\ge\frac{1}{3}$ which means $f_5(y_2)\ge f_3(y_2)-(-f_3(z_2))\ge\frac{1}{6}$ because $f_3(z_2)=-\frac{1}{6}$. But then $f_6(y_1)\ge f_5(y_1)+\frac{1}{6}\ge\frac{1}{6}$, a contradiction, as claimed. Hence, $G-x$ has a $3$-path $P(y_1y_2)$ because $G+xz_2\in\mathcal{F}_G$ and $G-\alpha_1$ has no $3$-path with ends $x$ and $y_1$. Let $P(y_1y_2)=y_1z_5z_6y_2$. Clearly, $z_5,z_6\in V_4$ and $f_2(z_6)\ge g^*(z_6)-\frac{1}{6}(n_4(z_6)-1)\ge\frac{1}{6}n_3(z_6)$ because $y_2\in N_3^+(z_6)$ and $z_5\in N_4^{1+}(z_6)$, which means $t_{z_6}^2\ge\frac{1}{6}$. Similarly, we have $f_3(y_2)\ge f_2(y_2)+t_{z_6}^2\ge\frac{1}{3}$ and $f_5(y_2)\ge\frac{1}{6}$. But then  $f_6(y_1)\ge f_5(y_1)+\frac{1}{6}\ge\frac{1}{6}$, a contradiction.

Assume $a_1\in V_2$. Then $a_1=x$, $a_2=\alpha_1$ and $a_3\in V_2$ as $g^*(x)<0$.  Since $ \alpha \in \mathcal{C}_4$, then $v \in \mathcal{C}_4$ for any $v \in V(G)$ with degree two. Then $G$ has a $5$-path with vertices $x,b_1,b_2,b_3,x_2,\alpha_2$ in order as $G+x\alpha_2\in\mathcal{F}_G$. Since $G$ has no $3$-path with ends $x$ and $\alpha_1$, we see $x_2\in V_2$. Note that $b_1\in V_3$ because $g^*(x)<0$ and $d(\alpha)=2$. Then $b_1\in V_3^-$, otherwise $b_1=y_1$ and we get a contradiction by similar analysis with the case $a_1\in V_3\cup V_4$.  Hence, $b_2\in V_4$ and $b_3\in V_3$.
We claim $t_{b_2}^2\ge\frac{1}{6}$. Suppose not. Then $n_3^{+}(b_2)=0$ and $n_3(b_2)+n_4(b_2)=2$ and $n_3^{-1}(b_2)\le1$, otherwise $f_2(b_2)\ge g^*(b_2)\ge\frac{1}{6}n_3(b_2)$ since $n_5^-(b_2)+n_4^{1-}(b_2)=0$. By the choice of $\alpha$, $n_5(b_2)\ne0$. Let $w\in N_5(b_2)$. Note that $g^*(b_2)\ge\frac{1}{6}n_3(b_2)-\frac{1}{3}$. Then $n_4^+(w)\le1$ and $d(w)=2$, otherwise $t_w^*\ge\frac{1}{3}$ as $b_2\in N_4^+(w)$ which yields $f_2(b_2)\ge g^*(b_2)+t_w^*\ge\frac{1}{6}n_3(b_2)$. Let $b_4\in N(w)\setminus b_2$. Since $d(w)=2$, we have $w$ belongs a $4$-cycle, which means $G-w$ has a $2$-path $P(b_2b_4)$. Let $P(b_2b_4)=b_2b_5b_4$. Note that $b_5\in V_3\cup V_5$. In fact, $b_5\in V_3$, that is, $b_5\in\{b_1,b_3\}$, otherwise $t_v^*\ge\frac{1}{6}$ for any $v\in\{w,b_5\}$ which implies $f_2(b_2)\ge g^*(b_2)+t_w^*+t_{b_5}^*\ge\frac{1}{6}n_3(b_2)$. W.l.o.g., let $b_5=b_3$. Then $b_4\in V_4^-$ as $n_4^+(w)\le1$. We further may assume $N(b_2)\cap N(b_4)=\{w,b_5\}$.  Then for any $i\in[2]$, $G$ has a $2$-path with vertices $\alpha_i,\alpha_i',x_2$ in order because $G+\alpha_iw\in\mathcal{F}_G$. Note that $\alpha_1'=\alpha_2'$, else $G$ has a copy of $C_6$ with vertices $x_2,\alpha_1',\alpha_1,\alpha,\alpha_2,\alpha_2'$ in order. This means $\alpha$ belongs to a graph $\Theta_5$. By the choice of $\alpha$, each vertex of degree two belongs to a  graph $\Theta_5$ which yields $G$ has a $2$-path with ends $b_3$ and $b_i$ for some $i\in\{2,4\}$. But then $b_4\in V_4^+$ or $n_3(b_2)+n_4(b_2)\ge3$, a contradiction, as claimed. Thus, $t_{b_2}^2\ge\frac{1}{6}$.
We see $n_4^-(b_1)=0$, otherwise $G+b_1'\alpha\notin \mathcal{F}_G$ because $g(x)<0$ and $g(b_1)<0$ for any $b_1'\in N_4^-(b_1)$.
 Hence,  $f_6(b_1)\ge f_3(b_1)\ge f_2(b_1)+t_{b_2}^2\ge\frac{1}{6}$ and so $t_{b_1}^6\ge\frac{1}{6}$, a contradiction.\qed\medskip

\noindent{\bf Acknowledgments.}
Lan was partially supported by National Natural Science Foundation of China (No.12001154), Natural Science Foundation of Hebei Province\,(No.A2021202025). Shi and Zhang were partially supported by the National
Natural Science Foundation of China (No.11922112), the Natural Science Foundation of Tianjin (Nos.20JCJQJC00090 and 20JCZDJC00840), the Fundamental Research Funds for the Central Universities, Nankai University (No.63213037) and the Funds for International Cooperation and Exchange of the National Natural Science Foundation of China\,(12161141006). Wang was partially supported by the National Natural Science Foundation of China (No.12071048) and the Science and Technology Commission of Shanghai Municipality (No.18dz2271000).

\vspace{-2mm}

\section*{Appendix}

\noindent{\bf Proof of Lemma~5.1}  Suppose $n^{-}_5(z)\ne0$. Let $w\in N^{-}_5(z)$. Then $d(w)=1$. Hence, $G$ has a $4$-path with vertices $y,z_1,w_2,w_1,z$ in order since $G+yw\in\mathcal{F}_G$. By Observation~\ref{ob}(1), we see that
$w_1\in V_5^1$ as $z \in V^{*-}_4$. That is, $d(w_1)=2$.
%$w_1,w_2\in V_5$, $z_1\in V_4$ and $d(w_1)=2$ since $z \in V^{*-}_4$.
%If $g(w_1) \ge \frac{2}{3}$, then $g^*(z) \ge 0$. Thus $g(w_1) \le \frac{1}{6}$, $d(w_1)=2$.
	Then $G$ has a $4$-path with ends $z$ and $w_2$ containing $w_1$ because $G+ww_2\in\mathcal{F}_G$. But this is impossible since $d(w_1)=2$. So we derive that $n^{-}_5(z)=0$. Note that $n_3(z)+n_4(z)=1$. By Observation~\ref{ob}(1), $n_5^1(z)=n_5(z)\le1$. Thus, $d(z)=2$, since otherwise $n_5^1(z)=n_5(z)=d(z)-n_3(z)-n_4(z)\ge2$. This proves $(a)$.\medskip
					
%Since $z \in V^{*-}_4$, we see that $u\in N_5(z)$ for every $u\in N(z)\setminus\{y\}$. By Proposition \ref{4-5-0}(i), we have $g(u)>0$ i.e., $g(u)\ge \frac{1}{6}$ for each $u\in N_5(z)$. But then $g^*(z)\ge0$ by Definition \ref{law1}(4), a contradiction. Thus, $d(z)=2$, as desired.
	
To prove ($b$), let $P=zww_1z_1$ be a $3$-path of $G$ such that $w\in V_5$.
By Observation~\ref{ob}, $w,w_1\in V_5^1$. This implies $z_1\in V_4$ and $d(w)=d(w_1)=2$. Note that $g(z_1)<0$ and $n_5^2(z_1)=0$, else $g_2(w)=\frac{1}{3}$ which yields $g^*(z)\ge g_3(z)\ge0$.

Now we shall show that $yz_1\notin E(G)$. Suppose not. Then $G$ has a 3-path with vertices $y,z_2,w_2,z_1$ in order because $G+wz_1\in \mathcal{F}_G$ and $d(z)=2$. Then $w_2\in V_5$ and $z_2\in V_4$ since $g(z_1)<0$.
But then  $w_2\in N_5^2(z_1)$, a contradiction. Hence $yz_1\notin E(G)$.
					
It remains to show that $d(z_1)=2$. Let $y_1\in N_3(z_1)$. Then $y_1\ne y$. Suppose $d(z_1)\ge3$.
	Let $w_2\in N(z_1)\setminus\{y_1,w_1\}$. Then, $w_2\in V_5$  and $d(w_2)\le2$ because $g(z_1)<0$ and $n_5^2(z_1)=0$. Moreover, $d(w_2)=2$, otherwise $G$ has a 4-path with ends $w$ and $z_1$ containing $w_1$ because $G+ww_2\in \mathcal{F}_G$ which means $d(w_1)\ge3 $.
	Let $w_3\in N(w_2)\setminus z_1$.
%\textcolor{red}{Then $w_3\in V_5$ and $d(w_3)=2$, else $g^*(z)\ge0$ by Definition \ref{law1}.}

Since $n_5^2(z_1)=0$, we see $w_3\in V_5$. Note that $n_4(w_3)=1$, else $n_4(w_3)\ge2$, which means $g_1(w_2)\ge\frac{1}{3}$, $g_3(z_1)\ge0$, $g_4(w)\ge\frac{1}{3}$ and $g^*(z)=g_5(z)\ge0$. Let $N_4(w_3)=\{z_2\}$. We claim that $yz_2\in E(G)$. Suppose not.  Note that $G$ has a $5$-path $P(w_1w_3)$ since $G+w_1w_3\in \mathcal{F}_G$. Let $P(w_1w_3)=w_1a_1a_2a_3a_4w_3$. Then $a_1\in\{w,z_1\}$. If $a_1=w$, then $a_2=z$, $a_3=y$. And so $a_4=z_2$ since $n_4(w_3)=1$, which contradicts to $yz_2\notin E(G)$. Hence, $a_1=z_1$. Then $G$ has a $4$-path with ends $z_1$ and $w_3$ containing $w_2$. But this is impossible because $d(w_2)=2$, as claimed.  Since $yz_1\notin E(G)$, we see $z_1\ne z_2$. Clearly, $g(z_2)<0$, else $g_2(w_2)\ge\frac{1}{3}$, $g_3(z_1)\ge0$, $g_4(w)\ge\frac{1}{3}$ and so $g^*(z)=g_5(z)\ge0$. Note that $G$ has a $5$-path $P(w_1w_2)$ because $G+w_1w_2\in \mathcal{F}_G$. Let $P(w_1w_2)=w_1b_1b_2b_3b_4w_2$. Then $b_1\in \{w,z_1\}$. If $b_1=w$, then $b_2=z$, $b_3=y$. And so $b_4=z_1$ since $d(w_2)=2$, which contradicts to $yz_1\notin E(G)$. Thus $b_1=z_1$ and so $b_4=w_3$. Since $g(z_1)<0$, we see that $b_2\in V_3\cup V_5$. In fact, $b_2\in V_5$, else $b_2=y_1$ and $b_3=z_2$ as $n_4(w_3)=1$, which contradicts to $g(z_2)<0$. Then $b_3\in V_5$ and $b_3z_1\notin E(G)$ because $n_5^2(z_1)=0$. But then $g_1(w_2)\ge\frac{1}{3}$, $g_3(z_1)\ge0$, $g_4(w)\ge\frac{1}{3}$ and so $g^*(z)=g_5(z)\ge0$, a contradiction, as desired.\qed

%We next claim $y_1z_2\in E(G)$. Suppose not. Then $G$ has a $5$-path $P(w_1w_2)$ because $G+w_1w_2\in \mathcal{F}_G$. Let $P(w_1w_2)=w_1b_1b_2b_3b_4w_2$. Then $b_1\in \{w,z_1\}$. If $b_1=w$, then $b_2=z$, $b_3=y$ and so $b_4=z_1$ because $d(w_2)=2$, which contradicts to $yz_1\notin E(G)$.
%Thus $b_1=z_1$ and so $b_4=w_3$. Since $g(z_1)<0$, we see that $b_2\in V_3\cup V_5$. Note that $b_2\in V_3$, else $b_3\in V_5$ because $n_5^2(z_1)=0$, which yields $g_1(w_2)\ge\frac{1}{3}$, $g_3(z_1)\ge0$, $g_4(w)\ge\frac{1}{3}$ and so $g^*(z)=g_5(z)\ge0$. Then $b_2=y_1$ and $b_3=z_2$ because $n_4(w_3)=1$, which contracts $y_1z_2\notin E(G)$. Thus, $y_1z_2\in E(G)$. But then $g_1(w_2)\ge\frac{1}{3}$, $g_3(z_1)\ge0$, $g_4(w)\ge\frac{1}{3}$ and so $g^*(z)=g_5(z)\ge0$, a contradiction. \qed

\noindent{\bf Proof of Lemma~\ref{34-little}:} Suppose not. Let $z_1,z_2\in N^{*-}_4(y)$. We claim that $d(z_i)\ge2$ for any $i\in[2]$. W.l.o.g., suppose $d(z_1)=1$. Then $d(z_2)\ge2$, else $G+z_1z_2\notin \mathcal{F}_G$. By Lemma~\ref{4-5-0}($a$), $d(z_2)=2$.
%{\color{red}By Lemma~\ref{4-5-0}($a$), $d(z_2)=2$, which  contradicts  to the fact $d(z_2)\ge 3$ by Lemma \ref{new}, as claimed.}
Let $w\in N(z_2)\setminus y$. Then $G$ has a $4$-path with ends $y$ and $w$ containing $z_2$ because $G+z_1w\in \mathcal{F}_G$. But this is impossible since $d(z_2)=2$,
 By Lemma~\ref{4-5-0}, for any $i\in[2]$, we see $d(z_i)=2$ and $G$ has a $3$-path $P^{i}=z_iw_iw_{i+2}z_{i+2}$ such that $d(w_i)=d(w_{i+2})=d(z_{i+2})=2$ and $yz_{i+2}\notin E(G)$.  This implies $V(P^1)\cap V(P^2)=\emptyset$. But then $G+z_1z_2\notin \mathcal{F}_G$, a contradiction.
% Hence, $V(P^1)\cap V(P^2)=\emptyset$. Let $C(z_1z_2)$ be a $6$-cycle of $G+z_1z_2$. Then $y\in V(C(z_1z_2))$ because $d_G(z_i)=2$ and $d_G(w_i)=2$ for each $i\in[4]$. Hence, either $yz_4\in E(G)$ or $yz_3\in E(G)$. But then $yz_4\notin E(G)$ and $yz_3\notin E(G)$, by Proposition \ref{4-5-0}(iii), a contradiction.
\qed\medskip

\noindent{\bf Proof of Lemma~\ref{55-0}:}
Note that $n_5^{-}(z)=1$.
%else let $v_1,v_2\in N_5^{-}(z)$, $G+v_1v_2\notin \mathcal{F}_G$ because $d(v_1)=d(v_2)=1$.
Let $N_5^{-}(z)=\{w\}$ and $y_1\in N_3(z)$. Clearly, $d(z)\ge3$, else $G+y_1w\notin \mathcal{F}_G$.

Firstly, we shall show that $f_2(z)\ge0$. Suppose $f_2(z)<0$. By ($\divideontimes$), we have $f_1(z)<0$. Clearly, $g^*(z)<\frac{1}{3}$, else $f_1(z)\ge g^*(z)-\frac{1}{3}n_5^-(z)\ge0$.  We next prove several claims.\medskip

\noindent{\bf Claim 1.} $n_3(z)=1$. \medskip

\noindent{\bf Proof.} Suppose not. Let $y_2\in N_3(z)\setminus y_1$.  Hence, for any $i \in [2]$, $n_2(y_i)=1$ and $N_2(y_1)=N_2(y_2)$, else let $x_i\in N_2(y_i)$, $G$ has a copy of $C_6$ with vertices $z,y_1,x_1,\alpha_1,x_2,y_2$  since $\delta(G)=1$, where $\alpha_1\in V_1$. Let $N_2(y_1)=N_2(y_2)=\{x\}$.
%By the partition of $V(G)$, we see $|V_1|=2$ because $\delta(G)=1$. Let $b\in V_1$ with $d(b)\ge2$. Then $N_1(v)=\{b\}$ for any $v\in V_2$. Hence, $N_2(y)=N_2(y_1)=\{x\}$, else $G$ has a copy of $C_6$.
Then $G$ has a $4$-path with ends $x$ and $z$ containing $y_1$ and $y_2$ since $G+wx\in\mathcal{F}_G$. Then $G-x-z$ has an $s$-path with ends $y_1$ and $y_2$ for some $s\in[2]$. Hence, $n_3^+(z)\ge1$ or $n_3^{-1}(z)\ge2$. But then by Ineq.~(\ref{g*}), $g^*(z)\ge\frac{1}{3}$, a contradiction. \qed
\medskip

\noindent{\bf Claim 2. $n_5(z)=1$. }\medskip

\noindent{\bf Proof.} Suppose not.  We first assert $d(v)=3$ for any $v\in N_5(z)\setminus w$. If $d(u)=2$ for some $u\in N_5(z)$, then let $N(u)=\{z,u'\}$, $G+wu'\notin\mathcal{F}_G$ because $G$ has no $4$-path with ends $z$ and $u'$ containing $u$, a contradiction.
%{\color{red} Since $d(w)=1$, $d(u)\ge 3$ for any $u\in N_5(z)\setminus\{w\}$ by Lemma \ref{new}. We first assert $d(v)=3$ for any $v\in N_5(z)\setminus w$. }
If $d(u)\ge4$ for some $u\in N_5(z)$, then $g^*(u)\ge\frac{1}{3}n_4(u)$ as $d(u)=n_4(u)+n_5(u)$, which implies $t_{u}^*\ge\frac{1}{3}$ and $f_1(z)\ge g^*(z)-\frac{1}{3}n_5^-(z)+t_{u}^*\ge0$, a contradiction, as asserted.
%Note that $d(v)\ge3$ for any $v\in N_5(z)\setminus w$, otherwise let $d(u)=2$ for some $u\in N_5(z)$ and $N(u)=\{z,u'\}$, $G+wu'\notin\mathcal{F}_G$ because $G$ has no $4$-path with ends $z$ and $u'$ containing $u$.
Hence, $N_5(z)\setminus w=N_5^2(z)$.
%Let $t_1=|N_4(w_1)|$ and $t_2=|N_5(w_1)|$. Then $g(w_1)=t_1+\frac{1}{2}t_2-\frac{4}{3}$.  If $d(w_1)\ge4$, then $t_1+t_2\ge4$ and so $g^*(w_1)\ge\frac{2}{3}t_1+\frac{1}{3}t_2-\frac{4}{3}\ge\frac{1}{3}t_1$, which implies that $f_1(z)\ge0$.
Note that $n_5(z)=2$, otherwise let $w_1,w_2\in N_5(z)\setminus w$, $g^*(w_i)\ge\frac{1}{6}n_4(w_i)$ because $n_4^{-1}(w_i)+n_4^+(w_i)\ge1$ and $n_4(w_i)+n_5(w_i)=3$ which yields $t_{w_i}^*\ge\frac{1}{6}$ and $f_1(z)\ge g^*(z)-\frac{1}{3}n_5^-(z)+t_{w_1}^*+t_{w_2}^*\ge0$. Let $N_5(z)=\{w,w_1\}$.
%Note that $d(w_1)\ge3$. Actually, we see $d(w_1)=3$ and $n_4^+(w_1)=0$, else $g^*(w_1)\ge\frac{1}{3}n_4(w_1)$ because $d(w_1)=n_4(w_1)+n_5(w_1)$ which implies $t_{w_1}^*\ge\frac{1}{3}$ and $f_1(z)\ge g^*(z)-\frac{1}{3}n_5^-(z)+t_{w_1}^*\ge0$.
%So we next assume that $n_5(z)=2$, $d(w_1)=3$ and $n_4^+(w_1)=0$.
Note that $n_4^+(w_1)=0$, otherwise $t_{w_1}^*\ge\frac{1}{3}$ and so $f_1(z)\ge0$. We next prove that\medskip
%Suppose $n_5(z)\ge3$. Let $w_2\in N_5(z)\setminus\{w,w_1\}$. Then $n_4^{-1}(w_i)\ge1$ and so $g^*(w_i)\ge\frac{1}{6}|N_4(w_i)|$.
%Hence, $f_1(z)\ge0$. So we next assume that $|N_5(z)|=2$.
%Suppose $G$ has a $2$-path $P(zw_1)=zz'w_1$. Then $z'\in V_4$. Hence, $g(z)\ge0$ and so $g^*(w_1)\ge\frac{2}{3}t_1+\frac{1}{3}+\frac{1}{3}t_2-\frac{4}{3}\ge\frac{1}{3}t_1$, which implies that $f_1(z)\ge0$.
%So we next assume $G$ has no $2$-path between $z$ and $w_1$.
%Suppose $g(u)\ge0$ for some $u\in N_4(w_1)$. Then $g^*(w_1)\ge\frac{2}{3}t_1+\frac{1}{3}+\frac{1}{3}t_2-\frac{4}{3}\ge\frac{1}{3}t_1$, which implies $f_1(z)\ge0$. So we next assume $g(u)<0$ for any $u\in N_4(w_1)$.

\noindent($a$) $n_4(w_1)=1$.\medskip

To see why ($a$) is true, suppose $n_4(w_1)\ge2$. Let $z_1\in N_4(w_1)\setminus z$. Then $g(z_1)<0$ as $n_4^+(w_1)=0$. Then $G$ has a $4$-path with ends $z$ and $z_1$ containing $w_1$ because $G+wz_1\in\mathcal{F}_G$. Since $d(w_1)=3$, $G-z_1$ has a $3$-path  $P(zw_1)$ or $G-z$ has a $3$-path $P(z_1w_1)$. We assert that $G-z_1$ has no $3$-path  $P(zw_1)$. Suppose not. Let $P(zw_1)=za_1b_1w_1$. Since $g(z)<0$ and $n_5(z)=2$, we see that $a_1=y_1$ and $b_1\in V_4\setminus z_1$. Then $G$ has a $4$-path with ends $z$ and $b_1$ containing $y_1$ and $w_1$ because $G+wb_1\in\mathcal{F}_G$, which implies $G-z-b_1$ has a $2$-path between $y_1$ and $w_1$. Hence, $y_1z_1\in E(G)$ as $d(w_1)=3$. Since $n_5(z)=2$, $n_3(z)=1$ and $g(z)<0$, $G$ has a $2$-path with ends $y_1$ and $v$ because $G+ww_1\in\mathcal{F}_G$, where $v\in\{b_1,z_1\}$. But then $g(v)>0$ and so $n_4^+(w_1)\ne0$, a contradiction, as asserted.
%Since $G+ww_1\in\mathcal{F}_G$, we see $G$ has a $3$-path $P(zv)=zz''z'''v$, where $v\in\{w_1',z_1\}$. Since $|N_5(z)|=2$ and $g(z)<0$, $z''=y$. This implies that $z'''\in V_3\cup V_4$, that is, $g(v)>0$, a contradiction.
Thus, $G-z$ has a $3$-path $P(z_1w_1)=z_1a_1b_1w_1$.  We see that $a_1\in V_3\cup V_5$ because $g(z_1)<0$. Assume first $a_1\in V_3$. Then $b_1\in V_4$.
Moreover, $n_5^-(z_1)+n_5^-(b_1)=0$, otherwise let $w_2\in N_5^-(z_1)$, we see $y_1=a_1$ because $G+w_2b_1\in\mathcal{F}_G$ which yields $G-z_1$ has a 3-path $P(zw_1)$ violating the above assert.
But then  $f_1(z)\ge g^*(z)-\frac{1}{3}n_5^-(z)+\frac{1}{3}\ge0$ by Definition~\ref{law2}(1.2), a contradiction.
Assume now $a_1\in V_5$. Then $b_1\in V_4\cup V_5$. Note that $b_1\in V_5$, otherwise by Ineq.~(\ref{g*}), we see $g^*(w_1)\ge\frac{1}{3}n_4(w_1)$ because $z_1,b_1\in N_4^{-1}(w_1)$, which implies $t_{w_1}^*\ge\frac{1}{3}$ and $f_1(z)\ge g^*(z)-\frac{1}{3}n_5^-(z)+t_{w_1}^*\ge0$. We see $d(a_1)\ge3$, else $f_1(z)\ge g^*(z)-\frac{1}{3}n_5^-(z)+\frac{1}{3}\ge0$ by Definition~\ref{law2}(1.1). By Ineq.~(\ref{g*}), $g^*(w_1)\ge\frac{1}{3}n_4(w_1)$ because $z_1\in N_4^{-1}(w_1)$ and $b_1\in N_5^2(w_1)$. That is, $t_{w_1}^*\ge\frac{1}{3}$. But then $f_1(z)\ge g^*(z)-\frac{1}{3}n_5^-(z)+t_{w_1}^*\ge0$, a contradiction. This proves ($a$).\medskip
%If $b_1\in V_5$ and $d(a_1)=2$, then $f_1(z)\ge g^*(z)-\frac{1}{3}n_5^-(z)+\frac{1}{3}\ge0$ .
% Note that $d(w_1)=3$ and $d(z_1')\ge2$.
% If $d(z_1')\ge3$, then $g^*(w_1)\ge\frac{2}{3}$ because $g(w_1')\ge\frac{2}{3}$ and $g(z_1')\ge\frac{2}{3}$, which implies that $f_1(z)\ge0$. If $d(z_1')=2$, then $g^*(w_1)\ge\frac{2}{3}$ by Definition~\ref{law2}(1.1), which means $f_1(z)\ge0$.

By ($a$), $N(w_1)\setminus z=N_5(w_1)$. Let $w_2\in N(w_1)\setminus z$.  Then $G$ has a $4$-path with ends $z$ and $w_2$ containing $w_1$ because $G+ww_2\in\mathcal{F}_G$. Since $g(z)<0$,  $n_4(w_1)=1$ and $n_5(z)=2$, we see $G-w_2$ has no $s$-path with ends $z$ and $w_1$ for any $s\in\{2,3\}$. Hence, $G-z$ has a $3$-path $P(w_1w_2)$.  % However, $P(zw_1)$ does not exist because $g(z)<0$ and $n_4(w_1)=1$.
% Suppose not. Let $P(zw_1)=zz'w_1'w_1$. Since $g(z)<0$ and $|N_5(z)|=2$, we see that $z'=y$. Hence, $w_1'\in N_4(w _1)$ which contracts to $N_4(w_1)=\{z\}$.
Let $P(w_1w_2)=w_1a_2b_2w_2$. Then $a_2\in V_5$. Note that $\{w_2,a_2\}\cap N_5^1(w_1)\ne\emptyset$, otherwise $g^*(w_1)\ge\frac{1}{3}n_4(w_1)$ because $w_2,a_2\in N_5^2(w_1)$ which yields $f_1(z)\ge0$. Hence, $b_2\in V_4$. By Definition~\ref{law1}(1.1), we see $g^*(w_1)=\frac{1}{3}$. But then $f_1(z)\ge 0$, a contradiction. \medskip\qed

%\textcolor{red}{If $b_2\in V_4$, ?? If $b_2\in V_4$, when $d(w_1)\ge 4$ or $\{w_2, a_2\} \subseteq N_5^2(b_2)$, then $g^*(w_1)\ge \frac{|N_4(w_1)|}{3}$ and so $f_2(z)\ge 0$; when $d(w_1)=3$ and $|\{w_2, a_2\} \cap N_5^1(b_2)|\ge 1$, by Definition \ref{law1} (1.1), then  $g^*(w_1)=\frac{1}{3}$ and so $f_2(z)\ge 0$. } If $b_2\in V_5$, then $g(w_2)\ge\frac{2}{3}$, $g(a_2)\ge\frac{2}{3}$ and so $g^*(w_1)\ge \frac{1}{3}$, which implies that $f_1(z)\ge0$.\medskip

By Claims 1 and 2,  $N(z)\setminus\{y_1,w\}=N_4(z)\ne\emptyset$ since $d(z)\ge3$. Let $z_1\in N_4(z)$ and $y_2\in N_3(z_1)$. We assert $G-y_2$ has no $2$-path with ends $z$ and $z_1$. Suppose not.  Let $P(zz_1)=za_3z_1$. Then $a_3\in V_3\cup V_4$. Clearly, $a_3\in V_3$, otherwise $a_3\in V_4$, $\{z_1,a_3\}\subseteq N_4^2(z)$ yields $g^*(z)\ge\frac{1}{3}$. Then $a_3=y_1$.
Hence, $z_1\in N_4^2(z)\cap N_4^2(y_1)$. But then $g^*(z)\ge\frac{1}{3}$, a contradiction, as asserted.
%Hence, $y_1\ne y_2$.
%If $a_3\in V_4$, then  $\{z_1,a_3\}\subseteq N_4^2(z)$, which means $g^*(z)\ge\frac{1}{3}$. But then in both cases, $f_1(z)\ge0$, a contradiction, as asserted.
Then $G$ has a $4$-path with ends $z$ and $y_2$ containing $z_1$ because $G+wy_2\in\mathcal{F}_G$. Hence, $G-y_2$ has a $3$-path $P(zz_1)$ or $G-z$ has a $3$-path $P(y_2z_1)$. %It is possible that $y_1=y_2$.% If $y_1=y_2$, we say $y_1y_2\notin E(G)$.

\medskip

%Hence, $G$ has a $t$-path $P(zz_1)$ and $(4-t)$-path $P(y_2z_1)$ for some $t\in[3]$.

%{\bf Case 2.1 $t=1$.}

\noindent{\bf Claim 3.} $G-z$ has no $3$-path $P(y_2z_1)$.\medskip

\noindent{\bf Proof.} Suppose not. Let $P(y_2z_1)=y_2a_4b_4z_1$. We assert $g(y_1)<0$. Suppose $g(y_1)\ge0$. Then $g(z)=g(z_1)=\frac{1}{6}$, otherwise $g^*(z)\ge \frac{1}{3}$.
% which implies $f_1(z)\ge g^*(z)-\frac{1}{3}n_5^-(z)\ge0$.
Hence, $d(z)=3$, $b_4\in V_5$ and $a_4\in V_4$. Then $y_1y_2\notin E(G)$(if  $y_1\ne y_2$), else $G$ has a copy of $C_6$ with vertices $y_1,z,z_1,b_4,a_4,y_2$ in order. Note that $n_5^-(z_1)=0$, otherwise let $w_1\in N_5^-(z_1)$, we see $G-z$ has a $2$-path with ends $y_1$ and $z_1$ because $G+ww_1\in\mathcal{F}_G$ and $d(z)=3$ which yields $g(z_1)\ge\frac{2}{3}$. By Observation~\ref{ob3}(2), $g^*(b_4)\ge\frac{1}{6}n_4(b_4)$ because $z_1\in N_4^{+}(b_4)$ and $n_4(b_4)\ge2$. Hence, $t_{b_4}^*\ge\frac{1}{6}$ and $f_1(z_1)\ge g^*(z_1)+t_{b_4}^*\ge\frac{1}{6}$. Note that $g^*(z)=\frac{1}{6}$ %第四层的1/6点不横向给
and $f_1(z)=g^*(z)-\frac{1}{3}=-\frac{1}{6}$. But then $f_2(z)\ge f_1(z)+\frac{1}{6}\ge0$, a contradiction, as asserted. Let $x_1\in N_2(y_1)$. Then $G$ has a $4$-path with ends $z$ and $x_1$ containing $y_1$ because $G+wx_1\in\mathcal{F}_G$. Since $g(y_1)<0$, we see $G-z$ has a $3$-path $P(x_1y_1)$ or  $G-x_1$ has a $3$-path $P(y_1z)$. We next prove that\medskip

\noindent($b$) $G-z$ has no $3$-path $P(x_1y_1)$. \medskip

To see why ($b$) is true, suppose $P(x_1y_1)$ exists and let $P(x_1y_1)=x_1c_4d_4y_1$. Then $d_4\in V_4$ and $c_4\in V_3$ as $g(y_1)<0$.
Hence, $d_4\in N_4^2(y_1)$, which means $g^*(z)\ge\frac{1}{6}$. Hence, $f_1(z)\ge g^*(z)-\frac{1}{3}n_5^-(z)\ge-\frac{1}{6}$. Note that $g(z_1)=\frac{1}{6}$, otherwise $g(z_1)\ge\frac{2}{3}$ yields $g^*(z)\ge\frac{1}{3}$.  Then $b_4\in V_5$ and $a_4\in V_4$. Similarly, we have $t_{b_4}^*\ge\frac{1}{6}$ and $y_1y_2\notin E(G)$(if $y_1\ne y_2$). Then $n_5^-(z_1)\ne0$, otherwise $f_1(z_1)\ge g^*(z_1)+t_{b_4}^*\ge\frac{1}{6}$ yields $f_2(z)\ge f_1(z)+\frac{1}{6}\ge0$.
Let $w_1\in N_5^-(z_1)$. Then $G-z_1-a_4$ has a $2$-path $P(y_2b_4)=y_2z_2b_4$ and $G$ has a $3$-path $P(zz_1)=zz'z_1'z_1$ because $\{G+w_1a_4,G+ww_1\}\in\mathcal{F}_G$. Then $z'\in V_4$ because $n_3(z)=n_5(z)=1$, $g(z_1)=\frac{1}{6}$ and $y_1y_2\notin E(G)$(if $y_1\ne y_2$). Moreover, $z'\notin\{z_2,a_4\}$, else let $z'=z_2$, we see $G$ has a copy of $C_6=b_4z_1zz'y_2a_4$.
Clearly, $z_1'\in V_5$ because $g(z_1)=\frac{1}{6}$. Then $z_1'\ne b_4$, else $G$ has a copy of $C_6=zz'z_1'z_2y_2z_1$. By Observation~\ref{ob2}(1), $g^*(z_1')\ge\frac{1}{3}n_4(z_1')$ because $z',z_1\in N_4^+(z_1')$. Then $t_{z_1'}^*\ge\frac{1}{3}$.  Hence, $f_1(z_1)\ge g^*(z_1)-\frac{1}{3}n_5^-(z_1)+t_{z_1'}^*+t_{b_4}^*\ge\frac{1}{6}$. But then $f_2(z)\ge f_1(z)+\frac{1}{6}\ge0$, a contradiction. This proves ($b$).\medskip

By ($b$), $G-x_1$ has a $3$-path $P(y_1z)$. Let $P(y_1z)=y_1c_{4}d_{4}z$. Then $c_4,d_4\in V_4$ because $g(y_1)<0$. We assert $d_4=z_1$. Suppose not. Note that $d_4\in N_4^2(z)$. Then $g(z_1)=\frac{1}{6}$,  $N_4^2(y_1)=\{z\}$ and $d(z)=4$, otherwise $n_3^{-1}(z)+n_4^2(z)\ge2$ or $n_3(z)+n_4(z)=d(z)-1\ge4$ which yields $g^*(z)\ge\frac{1}{3}$.
%
%If $g(z_1)\ge\frac{2}{3}$ or  $n_4^2(y_1)\ge2$ or $d(z)\ge5$, then $n_4^2(z)\ge2$ or $n_3^{-1}(z)=1$ or $n_3(z)+n_4(z)\ge4$, which means $g^*(z)\ge\frac{1}{3}$ and so $f_1(z)\ge0$. So we next assume $g(z_1)=\frac{1}{6}$,  $N_4^2(y_1)=\{z\}$ and $d(z)=4$.
This means $b_4\in V_5$ and $a_4\in V_4$. Similarly, $t_{b_4}^*\ge\frac{1}{6}$. Note that $n_5^-(z_1)=0$, otherwise let $w_1\in N_5^-(z_1)$, we see $G$ has a $3$-path $P(zz_1)=zz'z_1'z_1$ such that $z'=d_4$ and $z_1'\notin \{y_1,c_4,z\}$ because $G+ww_1\in\mathcal{F}_G$, which yields  $G$ has a copy of $C_6=y_1c_4d_4z_1'z_1z$. Hence, $f_1(z_1)\ge g^*(z_1)+t_{b_4}^*\ge\frac{1}{6}$. Note that $g^*(z)\ge\frac{1}{6}$ and so $f_1(z)\ge g^*(z)-\frac{1}{3}\ge-\frac{1}{6}$. But then $f_2(z)\ge f_1(z)+\frac{1}{6}\ge0$, a contradiction, as asserted.  Note that $z_1\in N_4^2(z)$.
%Then $z_1\in N_4^2(z)$, else $c_4=y_2$ yields $G$ has a copy of $C_6=y_1zz_1b_4a_4y_2$ .
Moreover, $N_4^2(y_1)\setminus z=\emptyset$, else  $g^*(z)\ge\frac{1}{3}$ because $z_1\in N_4^2(z)$.
Hence, $y_1z_1\notin E(G)$. Then $G-z-c_4$ has a $2$-path $P(y_1z_1)=y_1y_1'z_1$ because $G+wc_4\in\mathcal{F}_G$.
We see $c_4,y_1'\in N_4^1(y_1)$ because $g(y_1)<0$.
Note that $n_5^-(z_1)=0$, otherwise let $w_1\in N_5^-(z_1)$, we see $G$ has a $4$-path with ends $y_1$ and $z_1$ containing each of $\{z,c_4,y_1'\}$ because $G+y_1w_1\in\mathcal{F}_G$, which yields $zv\in E(G)$ for some $v\in\{c_4,y_1'\}$ and so $v\in N_4^2(y_1)$. If $b_4\in V_5$, then $t_{b_4}^*\ge \frac{1}{6}$.
By Ineq.~(\ref{g*}), $f_1(z_1)\ge g^*(z_1)+t_{b_4}^*\ge\frac{1}{6}n_4(z_1)$ because $n_4(z_1)\ge3$.  If $b_4\in V_3\cup V_4$, then $b_4\notin \{ y_1', c_4\}$ because $y_1', c_4 \in V_4^1$ and $y_1y_2\notin E(G)$. Note that $b_4\neq z$. Thus $n_4(z_1)\ge 4$ and by Ineq.~(\ref{g*}), $f_1(z_1)\ge g^*(z_1)\ge\frac{1}{6}n_4(z_1)$.
 Since  $n_4^2(z)\ge1$ and $n_4^1(y_1)\ge 2$, we have $g^*(z)\ge\frac{1}{6}$ and so $f_1(z)\ge g^*(z)-\frac{1}{3}\ge-\frac{1}{6}$.
But then $f_2(z)\ge f_1(z)+\frac{1}{6}\ge0$, a contradiction.  \medskip \qed

By Claim 3, $G-y_2$ has a $3$-path $P(zz_1)$. Let $P(zz_1)=za_5b_5z_1$. Then $a_5\in V_3\cup V_4$. We next prove that\medskip

\noindent($c$) $a_5\in V_4$.\medskip

To prove ($c$), suppose $a_5\in V_3$. Then $a_5=y_1$ because $n_3(z)=1$. Moreover, $z_1\in N_4^2(z)$ as $b_5\ne y_2$.
%Hence, $g^*(z)\ge\frac{1}{6}$ and $f_1(z)\ge g^*(z)-\frac{1}{3}n_5^-(z)\ge-\frac{1}{6}$.
%This implies $z_1\in N_4^2(z)$. Hence, $f_1(z)\ge g^*(z)-\frac{1}{3}n_5^-(z)\ge-\frac{1}{6}$.
Hence, $y_1z_1\notin E(G)$, $g(y_1)<0$ and $N_4^2(y_1)\setminus z=\emptyset$, else $g^*(z)\ge\frac{1}{3}$.  Then $G-z-b_5$ has a $2$-path $P(y_1z_1)$ because $G+wb_5\in \mathcal{F}_G$.
Let $P(y_1z_1)=y_1c_5z_1$. We see
$b_5, c_5\in V_4^1$. It is easy to see $g^*(z)\ge\frac{1}{6}$ and $f_1(z)\ge g^*(z)-\frac{1}{3}n_5^-(z)\ge-\frac{1}{6}$. Then $n_5^-(z_1)=0$, otherwise let $w_1\in N_5^-(z_1)$, we see $G$ has a $4$-path with ends $y_1$ and $z_1$ containing each of $\{z,b_5,c_5\}$ because $G+y_1w_1\in\mathcal{F}_G$, which yields $zv\in E(G)$ for some $v\in\{b_5,c_5\}$ and so $v\in N_4^2(y_1)$.
Note that $n_4^2(z_1)=0$, else $f_1(z_1)\ge g^*(z_1)\ge \frac{1}{6}n_4(z_1)$ because $n_3(z_1)+n_4(z_1)\ge4$ and $n_4^2(z_1)\ge1$, which implies $f_2(z)\ge f_1(z)+\frac{1}{6}\ge0$. This follows $d(z)=3$. Moreover, $n_5^-(b_5)\ne0$ or $n_5^-(c_5)\ne0$, else $f_1(z_1)\ge g^*(z_1)\ge \frac{1}{6}(n_4(z_1)-1)\ge \frac{1}{6}n_4^{1-}(z_1)$ because $n_3(z_1)+n_4(z_1)\ge4$ and $b_5\in N_4^{1+}(z_1)$, which means $f_2(z)\ge f_1(z)+\frac{1}{6}\ge0$. W.l.og, let $w_2\in N_5^-(b_5)$. Then $G$ has a 3-path $P(zb_5)$ because $G+ww_2\in \mathcal{F}_G$.
Let $P(zb_5)=zz'b_5'b_5$. We see $z'=z_1$, else $z'=y_1$ because $d(z)=3$, which means $b_5\in N_4^2(z_1)$ or $y_1z_1\in E(G)$. Moreover, $b_5'\in V_5$. By Observation~\ref{ob2}(1), we have $t_{b_5'}^*\ge\frac{1}{3}$ because $z_1,b_5\in N_4^+(b_5')$. This implies $b_5\in N_4^{1+}(z_1)$. Similarly, we have $f_1(z_1)\ge \frac{1}{6}n_4^{1-}(z_1)$ and so $f_2(z)\ge0$, a contradiction. This proves $(c)$.\medskip

%But then $b_5\in N_4^2(z_1)$ because $d(z)=3$, a contradiction. This proves $(c)$.\medskip

By $(c)$, $a_5\in V_4$. Then $\{z_1,a_5\}\cap N_4^1(z)\ne\emptyset$, else $g^*(z)\ge\frac{1}{3}$. W.l.o.g., let $z_1\in N_4^1(z)$. Then $b_5\in V_5$ because $b_5\ne y_2$.
Note that $z_1a_5\notin E(G)$. Hence, $G-z-b_5$ has a $2$-path $P(z_1a_5)=z_1d_5a_5$ because $G+wb_5\in \mathcal{F}_G$.
%Let $P(z_1a_5)=z_1d_5a_5$.
%Hence, $\{b_5,d_5\}\cap V_5\ne\emptyset$. W.l.o.g, let $b_5\in V_5$.
Clearly, $n_5^-(z_1)=0$, otherwise let $w_1\in N_5^-(z_1)$ we see $G$ has a $4$-path with ends $z_1$ and $a_5$ containing each of $\{z,b_5,d_5\}$ because $G+w_1a_5\in \mathcal{F}_G$, which implies $d_5\in V_4$ and so $z_1\notin N_4^1(z)$. By Observation~\ref{ob2}(1), $g^*(b_5)\ge\frac{1}{3}n_4(b_5)$ and so $t^*_{b_5}\ge \frac{1}{3}$ because $z_1,a_5\in N_4^+(b_5)$, which yields $f_1(z_1)\ge g^*(z_1)+t^*_{b_5}\ge\frac{1}{3}$.
If $a_5\notin N_4^1(z)$, then $g^*(z)\ge\frac{1}{6}$, which yields $f_2(z)\ge f_1(z)+\frac{1}{6}\ge0$, a contradiction. If $a_5\in N_4^1(z)$, then we have $f_1(a_5)\ge\frac{1}{3}$ by similar analysis with $z_1$, which implies $f_2(z)\ge f_1(z)+\frac{1}{6}\times 2\ge0$, a contradiction.
%
%
%Since $\{z_1, a_5\}\cap N_4^1(z)\neq \emptyset$, we see $f_1(v)\ge g^*(v)+t^*_{b_5}\ge \frac{1}{3}$ for any $v \in \{z_1, a_5\}\cap N_4^1(z)$ and so
% $f_2(z)\ge f_1(z)+\frac{1}{6}\ge0$, a contradiction.
 This completes the proof of $f_2(z)\ge0$.\medskip

Finally, we shall show that $f_3(y) \ge 0$ for any $y \in N_3(z)$. Suppose not. Let $f_3(y_1)<0$ for some $y_1\in N_3(z)$. By ($\divideontimes$), $g^*(y_1)<0$. By Observation~\ref{ob}(1), $n_4^2(y_1)=0$ and $n_4^1(y_1)\le1$. Let $x_1\in N_2(y_1)$. Then $G$ has a $4$-path with ends $x_1$ and $z$ containing $y_1$ because $G+wx_1\in\mathcal{F}_G$. This implies that $G$ has a $3$-path $P(y_1z)=y_1z_1w_1z$ such that $w_1\in V_5$ and $z_1\in V_4$ because $g^*(y_1)<0$. Since $G+wz_1\in\mathcal{F}_G$, we have $G-z-z_1$ has a $2$-path $P(y_1w_1)=y_1z_2w_1$. Clearly, $g(z_i)<0$ for some $i\in[2]$. W.l.o.g., let $g(z_1)<0$. We claim that $n_5^-(z_1)+ n_5^-(z_2)\ne0$. Suppose not. We say $g(v)<0$ for any $v\in\{z,z_2\}$, otherwise by Observation~\ref{ob2}(2), $t_{w_1}^*\ge\frac{1}{3}$ because $n_4^+(w_1)\ge1$ and $n_4(w_1)\ge3$, which implies $f_2(z_1)\ge f_1(z_1)\ge t_{w_1}^*$, that is, $t_{z_1}^2\ge\frac{1}{3}$, and so $f_3(y_1)\ge f_2(y_1)+t_{z_1}^2\ge0$. By Observation~\ref{ob3}(2), $t_{w_1}^*\ge\frac{1}{6}$ because $n_4(w_1)\ge3$.
% If $g(v)>0$ for some $v\in\{z,z_2\}$, then $f_1(y')\ge\frac{1}{3}$ and so $f_3(y)\ge0$ because $g(y')<0$, a contradiction.
Hence, $f_2(z_i)\ge f_1(z_i)\ge t_{w_1}^*\ge\frac{1}{6}$ for any $i\in[2]$.  That is, $t_{z_i}^2\ge\frac{1}{6}$. But then $f_3(y_1)\ge f_2(y_1)+t_{z_1}^2+t_{z_2}^2\ge0$, a contradiction, as claimed. Thus, $n_5^-(z_1)+n_5^-(z_2)\ne0$. W.l.o.g., let $w_2\in N_5^-(z_2)$. Then $G$ has a $3$-path $P(zz_2)=za_6b_6z_2$ because $G+ww_2\in\mathcal{F}_G$.
We claim that $a_6,b_6\in V_5$. Suppose not. W.l.o.g., let $a_6\notin V_5$.   Then $a_6\in V_4$ and $b_6\in V_5$ because $n_4^2(y_1)=0$ and $n_4^1(y_1)\le1$. Clearly, $f_2(y_1)\ge g^*(y_1)\ge-\frac{1}{6}$. Moreover, $g(z_2)<0$. We see $b_6\ne w_1$, otherwise $t_{w_1}^*\ge\frac{1}{2}$ because $n_4^+(w_1)\ge2$ and $n_4(w_1)\ge4$, which implies $t_{z_2}^2\ge\frac{1}{6}$ and $f_3(y_1)\ge f_2(y_1)+t_{z_2}^2\ge0$. By Observations~\ref{ob2}(2) and \ref{ob3}(2),  $t_{w_1}^*\ge\frac{1}{3}$ because $n_4(w_1)\ge3$ and $z\in N_4^+(w_1)$, and $t_{b_6}^*\ge\frac{1}{6}$ because $n_4(b_6)\ge2$ and $a_6\in N_5^+(b_6)$.  Hence, $f_2(z_2)\ge f_1(z_2)\ge g^*(z_2)-\frac{1}{3}n_5^-(z_2)+t_{w_1}^*+t_{b_6}^*\ge\frac{1}{6}$, that is, $t_{z_2}^2\ge\frac{1}{6}$.  But then $f_3(y_1)\ge f_2(y_1)+t_{z_2}^2\ge0$, a contradiction, as claimed. Thus, $a_6,b_6\in V_5$. Obviously, $d(a_6)\ge3$ and $d(b_6)\ge3$, else $G+wb_6\notin\mathcal{F}_G$ or $G+w_2a_6\notin\mathcal{F}_G$. Hence, $n_4(v)+n_5(v)=d(v)\ge3$ for any $v\in\{a_6,b_6\}$. By Ineq.~(\ref{g*}), $t_{v}^*\ge\frac{1}{3}$ because $n_5^2(v)\ge1$ and $n_4^+(v)+n_4^{-1}(v)+d(v)\ge4$.
We assert $g(u)<0$ for any $u\in\{z,z_2\}$. Suppose not. Then $f_2(y_1)\ge g^*(y_1)\ge-\frac{1}{6}$. Similarly, $t_{w_1}^*\ge\frac{1}{3}$. W.l.o.g., we assume $w_1\ne a_6$.
%Moreover, $n_4^+(w_1)\ge1$ and $n_4(w_1)\ge3$, which implies $g^*(w_1)\ge\frac{1}{3}n_4(w_1)$ by Observation~\ref{ob2}(4).
Then $f_1(z)\ge g^*(z)-\frac{1}{3}n_5^-(z)+t_{a_6}^*+t_{w_1}^*\ge\frac{1}{3}$, which means $f_2(z)\ge\frac{1}{6}$ because $g(z)\le\frac{1}{6}$. That is, $t_{z}^2\ge\frac{1}{6}$. But then $f_3(y_1)\ge f_2(y_1)+t_z^2\ge0$, a contradiction, as asserted.
If $w_1\in\{a_6,b_6\}$, say $w_1=b_6$, then $t_{w_1}^*\ge\frac{1}{3}$ because $n_4(w_1)+n_5(w_1)\ge4$. Hence, $f_2(z)\ge f_1(z)\ge g^*(z)-\frac{1}{3}n_5^-(z)+t_{w_1}^*+t_{a_6}^*\ge\frac{1}{3}$. That is, $t_{z}^2\ge\frac{1}{3}$. But then $f_3(y_1)\ge f_2(y_1)+t_z^2\ge0$, a contradiction. If $w_1\notin\{a_6,b_6\}$, then $t_{w_1}^*\ge\frac{1}{6}$ because $n_4(w_1)\ge3$. Hence, $f_2(u)\ge f_1(u)\ge g^*(u)-\frac{1}{3}n_5^-(u)+t_{u'}^*+t_{w_1}^*\ge\frac{1}{6}$ for any $u\in\{z,z_2\}$, where $u'\in \{a_6,b_6\}\cap N(u)$. That is, $t_{u}^2\ge\frac{1}{6}$. But then $f_3(y_1)\ge f_2(y_1)+t_z^2+t_{z_2}^2\ge0$, a contradiction. Thus, $f_3(y) \ge 0$ for any $y \in N_3(z)$. \qed\medskip

\noindent{\bf Proof of Lemma~\ref{4-53-0 }:} Let $w\in N^{-}_5(z_1)$ and $z_1\in N_4(y)$. Note that $d(w)=1$.  By Lemma~\ref{55-0} and $(\divideontimes)$, $f_5(y)\ge0$ and $f_2(z_1)\ge0$.
To establish the desired result, suppose $f_5(z)<0$ for some $z\in N_4(y)\less z_1$. Then $g^*(y)\le f_3(y)<\frac{1}{3}$. By ($\divideontimes$), $f_i(z)<0$ for any $i\in[5]$. By Corollary~\ref{f2g}, $g^*(z)<0$. By Lemma~\ref{4-5-0}($a,b$), $d(z)\le2$ and $G$ has a $3$-path $P=zw_1w_2z_2$ such that $d(u)=2$ for each $u\in V(P)$ and $yz_2\notin E(G)$ when $d(z)=2$, where $w_1,w_2\in V_5$, $z_2\in V_4$.
When $d(z)=2$, let $N_3(z_2)=\{y_2\}$.
We first prove that \medskip

\noindent{$(a)$ for any $v\in N_4^2(y)$, $n_3^+(v)\le1$, or $n_3^+(v)=2$ and $n_4^2(v)=0$, or $n_3^+(v)=2$ and $N_4^1(v) \cap N_4(N_3(v))=\emptyset$.}\medskip

To see why $(a)$ is true, suppose first $n_3^+(v)\ge3$, or $n_3^+(v)=2$ and $n_4^2(v)\ge1$. By Ineq.~(\ref{g*}), $g^*(v)\ge\frac{1}{3}n_3(v)+\frac{1}{6}n_4(v)+\frac{1}{3}$.  Thus, $f_2(v)\ge \frac{1}{3}n_3(v)$.   Suppose $n_3^+(v)=2$  and $N_4^1(v) \cap N_4(N_3(v))\ne \emptyset$. Let $y_1\in N_3(v)$ and $v_1\in N_4^1(v)\cap N_4(y_1)$. By Ineq.~(\ref{g*}), $g^*(v)\ge\frac{1}{3}n_3(v)+\frac{1}{6}n_4(v)+\frac{1}{6}$ and $g^*(v_1)=\frac{1}{3}$, which means $f_1(v_1)\ge 0$ and so $n_4^{1-}(v)\le n_4(v)-1$. Hence, $f_2(v)\ge f_1(v)-\frac{1}{6}n_4^{1-}(v)\ge \frac{1}{3}n_3(v)$. In both cases, we see $t_{v}^2\ge\frac{1}{3}$. But then $f_3(y)\ge f_2(y)+t_{v}^2\ge\frac{1}{3}$, a contradiction.
This proves $(a)$.\medskip

Then $G$ has a $3$-path $P(yz_1)$ because $G+zw\in\mathcal{F}_G$ and $g^*(z)<0$. Let $P(yz_1)=ya_1b_1z_1$.
We then prove that
\medskip

\noindent{($b$) $a_1z_1\notin E(G)$.}\medskip

To prove $(b)$, suppose $a_1z_1\in E(G)$. We  claim that $b_1\in V_5$.  Suppose not. Then $b_1\in V_3\cup V_4$.
Note that $b_1\in V_4$, otherwise
by ($a$) and Ineq.~(\ref{g*}),  we see $a_1\in V_4$ and so $g^*(v)\ge \frac{1}{6}n_3(v)+\frac{1}{6}n_4(v)+\frac{1}{3}$ for any $v\in \{z_1, a_1\}$ because $n_3^{-1}(v)+n_3^+(v)\ge 2$ and $n_4^2(v)\ge 1$, which means $t_{v}^2\ge \frac{1}{6}$ and so $f_3(y)\ge f_2(y)+t_{z_1}^2+t_{a_1}^2\ge\frac{1}{3}$.
By $(a)$, $a_1\in V_4$.
Hence, $\{a_1, z_1,  b_1\}\subseteq V_4^2$, which means that $g^*(v')\ge \frac{1}{3}$ and $f_1(v')\ge 0$ for any $v' \in \{a_1, z_1, b_1\}$. By Observation~\ref{ob3}(2,3,5), $g^*(v)\ge \frac{1}{6}n_3(v)+\frac{1}{6}n_4(v)$ for any $v\in \{z_1, a_1\}$ because $n_3^{-1}(v)+n_3^+(v)\ge 1$ and $n_4^2(v)\ge 2$. This implies $t_{v}^2\ge\frac{1}{6}$ as $n_4^{1+}(v)\ge2$. But then $f_3(y)\ge f_2(y)+t_{z_1}^2+t_{a_1}^2\ge\frac{1}{3}$, a contradiction, as claimed.
Then $G-z_1-a_1$ has a $2$-path $P(yb_1)$ because $G+a_1w\in\mathcal{F}_G$. Let $P(yb_1)=yc_1b_1$. Then $G$ has a $4$-path with ends $y$ and $b_1$ containing each of $\{z_1,c_1,a_1\}$ since $G+zb_1\in\mathcal{F}_G$, which means $c_1a_1\in E(G)$ or $c_1z_1\in E(G)$.  W.l.o.g., let $c_1a_1\in E(G)$. Then $a_1\in V_4^2$. Note that $g^*(b_1)\ge \frac{1}{2}n_4(b_1)$ as $n_4^+(b_1)\ge 3$, which means $t_{b_1}^*\ge\frac{1}{2}$.
Then $v\in V_4^2$ for any $v\in \{z_1, c_1\}$, otherwise let $z_1\in V_4^1$, we see $g^*(z_1)= \frac{1}{3}$ and so $f_2(z_1)\ge g^*(z_1)-\frac{1}{3}-\frac{1}{6}+t_{b_1}^*\ge \frac{1}{3}$, which yields $t_{z_1}^2\ge\frac{1}{3}$ and $f_3(y)\ge f_2(y)+t_{z_1}^2\ge\frac{1}{3}$.
Hence, for any $v\in\{a_1,z_1,c_1\}$, $g^*(v)\ge \frac{1}{6}n_3(v)+\frac{1}{6}(n_4(v)-1)$ because $n_3^+(v)+n_3^{-1}(v)\ge 1$ and $n_4^2(v)\ge 1$, which implies $t_{v}^2\ge\frac{1}{6}$ because $f_1(v)\ge g^*(v)-\frac{1}{3}+t_{b_1}^*$. But then $f_3(y)\ge f_2(y)+t_{z_1}^2+t_{a_1}^2+t_{c_1}^2\ge\frac{1}{2}$, a contradiction. This proves $(b)$.
\medskip

\noindent{$(c)$ $yb_1\notin E(G)$.}\medskip

%Furthermore, there are no such  two vertice $\{a_1, b_1\}\subseteq N(y)$ and $\{z_1, a_1\} \subseteq N(b_1)$ in $G$.
To see why $(c)$ is true, suppose $yb_1\in E(G)$.  Then $b_1\in V_3\cup V_4$. We assert $b_1\in V_4$. Suppose not. Then $b_1\in V_3$. Then $a_1\notin V_4$, otherwise $g^*(v)\ge \frac{1}{6}n_3(v)+\frac{1}{6}n_4(v)+\frac{1}{3}$ for any $v \in \{z_1, a_1\}$ because $n_3^+(v)\ge 2$, which implies $t^2_{v}\ge \frac{1}{6}$ and so $f_3(y)\ge f_2(y)+t^2_{z_1}+t^2_{a_1}\ge \frac{1}{3}$. Since $g^*(y)<\frac{1}{3}$, we see $a_1\notin V_3$. Hence, $a_1\in V_2$. Note that $y_2b_1\notin E(G)$ when $d(z)=2$, otherwise $g^*(y)\ge\frac{1}{3}$. Then $G-y-b_1$ has a $2$-path with ends $a_1$ and $z_1$ because $G+zb_1\in\mathcal{F}_G$. This implies $n_3(z_1)\ge3$. By Ineq.~(\ref{g*}), $g^*(z_1)\ge \frac{1}{3}n_3(z_1)+\frac{1}{6}n_4(v)+\frac{1}{3}$, which yields $t_{z_1}^2\ge\frac{1}{3}$. But then $f_3(y)\ge f_2(y)+t_{z_1}^2\ge0$, a contradiction, as asserted. By $(a)$, $a_1\in V_4$.
%Because $g^*(y)<\frac{1}{3}$, then $\{a_1,b_1\}\cap V_4\neq \emptyset$. By (a), $a_1\in V_4$. If $b_1\in V_3$, then $g^*(v)\ge \frac{n_3(v)}{6}+\frac{n_4(v)}{6}+\frac{1}{3}$ for any $v \in \{z_1, a_1\}$ because $n_3^+(v)\ge 2$. Thus $f_2(v)\ge \frac{n_3(v)}{6}$, $t^2_{v}\ge \frac{1}{6}$ for any $v \in \{z_1, a_1\}$ and so $f_3(y)\ge g^*(y)+t^2_{z_1}+t^2_{a_1}\ge \frac{1}{3}$, a contradiction, as claimed.
 We assert  $y_2b_1\notin E(G)$. Suppose not. Then $d(z)=2$ and $g^*(z)\ge -\frac{1}{6}$. Hence, $f_3(y)<\frac{1}{6}$, otherwise $f_4(y)\ge0$ and $f_4(z)\ge0$. By Observation~\ref{ob3}(1), $g^*(b_1)\ge \frac{1}{6}n_3(b_1)+\frac{1}{6}n_4(b_1)$ because $n_3(b_1)\ge2$ and $n_3(b_1)+n_4(b_1)\ge4$.  Clearly, $g^*(v)\ge \frac{1}{3}$ for any $v\in \{z_1, a_1\}$ because $b_1\in V_4^2$, which implies $\{z_1,a_1\}\subseteq N_4^{1+}(b_1)$.
Thus $t_{b_1}^2\ge\frac{1}{6}$ as $n_4^{1-}(b_1)\le n_4(b_1)-2$, which yields that $f_3(y)\ge f_2(y)+t_{b_1}^2\ge\frac{1}{6}$, a contradiction, as asserted.
By $(b)$, $G-y-b_1$ has a 2-path $P(z_1a_1)$ because $G+zb_1\in \mathcal{F}_G$.  Let $P(z_1a_1)=z_1c_1a_1$. We assert $c_1\in V_5$. Suppose not.  Then $c_1\in V_4$, otherwise $g^*(v)\ge \frac{1}{6}n_3(v)+\frac{1}{6}(n_4(v)+2)$ because $n_3^{-1}(v)+n_3^+(v)\ge 2$ and $n_4^2(v)\ge 1$ for any $v\in \{z_1, a_1\}$, which means that $t_{v}^2\ge\frac{1}{6}$ and so $f_3(y)\ge f_2(y)+t_{z_1}^2+t_{a_1}^2 \ge\frac{1}{3}$. By Observation~\ref{ob3}(2,3,5),  $g^*(v)\ge \frac{1}{6}n_3(v)+\frac{1}{6}n_4(v)$ for any $v \in \{z_1, a_1, b_1\}$ and $g^*(c_1)\ge \frac{1}{3}$ because $n_4^2(c_1)\ge 2$, which implies $n_4^{1+}(v)\ge 2$. Hence, $t_v^2\ge\frac{1}{6}$ because $n_4^{1-}(v)\le n_4(v)-2$ for any $v \in \{z_1, a_1, b_1\}$. But then $f_3(y)\ge f_2(y)+t_{z_1}^2+t_{a_1}^2+t_{b_1}^2\ge\frac{1}{2}$, a contradiction, as asserted.  By Observation~\ref{ob2}(1),  $g^*(c_1)\ge \frac{1}{3}n_4(c_1)+\frac{1}{6}n_5(c_1)$ because $n_4^+(c_1)\ge2$, which means $t_{c_1}^*\ge\frac{1}{3}$.
We assert $t_{v}^2\ge \frac{1}{6}$ for any  $v \in \{z_1, a_1\}$. W.l.o.g., suppose $t_{z_1}^2<\frac{1}{6}$. Then $z_1\in V_4^2$, otherwise $g^*(z_1)\ge \frac{1}{3}$ and so $f_1(z_1)\ge g^*(z_1)-\frac{1}{3}+t_{c_1}^*\ge \frac{1}{3}$, which means $t_{z_1}^2\ge \frac{1}{6}$. By Ineq.(3), $g^*(z_1)\ge \frac{1}{6}n_3(z_1)+\frac{1}{6}(n_4(z_1)-1)\ge\frac{1}{3}$ because $n_3^+(z_1)+n_3^{-1}(z_1)\ge1$ and $b_1\in N_4^2(z_1)$. Similarly, $g^*(b_1)\ge\frac{1}{3}$, which means $b_1\in N_4^{1+}(z_1)$ and
$n_4^{1-}(z_1)\le n_4(z_1)-1$. But then
$f_1(z_1)\ge g^*(z_1)-\frac{1}{3}n_5^-(z_1)+t_{c_1}^*\ge g^*(z_1)$, which yields that $t_{z_1}^2\ge \frac{1}{6}$, a contradiction, as asserted. Hence,
 $f_3(y)\ge f_2 (y)+t_{z_1}^2+t_{a_1}^2\ge\frac{1}{3}$, a contradiction.
This proves $(c)$. \medskip
% For any  $v \in \{z_1, a_1\}$, we claim that  $t_{v}^2\ge \frac{1}{6}$.
%If $v\in V_4^1$, then $g^*(v)\ge \frac{1}{3}$ and so $f_1(v)\ge g^*(v)-\frac{1}{3}+t_{c_1}^*\ge \frac{1}{3}$, which means $t_{v}^2\ge \frac{1}{6}$.
%If $v \in V_4^2$,
% let $v\in z_1$,
%then  $g^*(u)\ge \frac{1}{6}n_3(u)+\frac{1}{6}(n_4(u)-1)$ for any $u \in \{z_1, b_1\}$ and so $g^*(u)\ge \frac{1}{3}$, which means $f_1(u)\ge 0$ and
%$n_4^{1-}(u)\le n_4(u)-1$. Then
%$f_1(z_1)\ge g^*(z_1)-\frac{1}{3}+t_{c_1}^*\ge g^*(z_1)$, which yields that $t_{z_1}^2\ge \frac{1}{6}$. Thus $t_{v}^2\ge \frac{1}{6}$.
%Hence,
% $f_3(y)\ge f_2 (y)+t_{z_1}^2+t_{a_1}^2\ge\frac{1}{3}$, a contradiction.
%This proves $(c)$. \medskip

By $(b)$ and $(c)$, we may assume any $4$-cycle containing the edge $yz_1$ is an induced cycle in $G-z$. By  $(c)$, $G-z_1-a_1$ has a $2$-path $P(yb_1)$ since $G+wa_1\in\mathcal{F}_G$. Let $P(yb_1)=yc_1b_1$.  We assert $y_2b_1\notin E(G)$. Suppose not. Then $d(z)=2$ and $b_1\in V_3\cup V_4$. Hence, $f_3(y)<\frac{1}{6}$ as $f_3(z)\ge g^*(z)\ge-\frac{1}{6}$. If $b_1\in V_3$, then $g(y)<0$, otherwise by Observation~\ref{ob2}(1),  $g^*(z_1)\ge \frac{1}{3}n_3(z_1)+\frac{1}{6}n_4(z_1) \ge\frac{1}{6}n_3(z_1)+\frac{1}{6}n_4(z_1)+\frac{1}{3}$ because $\{y,b_1\}\subseteq N_3^+(z_1)$, which means that $t_{z_1}^2\ge\frac{1}{6}$ and so $f_3(y)\ge f_2(y)+t_{z_1}^2\ge\frac{1}{6}$. Then $\{a_1,c_1\}\cap V_4\ne\emptyset$. W.l.o.g., let $c_1\in V_4$. By Ineq.~(\ref{g*}), for $v\in\{z_1,c_1\}$, $g^*(v)\ge\frac{1}{12}n_3(v)+\frac{1}{6}(n_4(v)+2)$ because $b_1\in N_3^+(v)$ and $y\in N_4^{-1}(v)$, which yields $t_{v}^2\ge\frac{1}{12}$. But then $f_3(y)\ge f_2(y)+t_{z_1}^2+t_{c_1}^2\ge\frac{1}{6}$, a contradiction.
If $b_1\in V_4$, then $y_2z_1\notin E(G)$, otherwise by Ineq. (\ref{g*}), $b_1\in N_4^{1+}(z_1)$ and $g^*(z_1)\ge \frac{1}{6}n_3(z_1)+\frac{1}{6}(n_4(z_1)-1)+\frac{1}{3}$ because $n_3(z_1)\ge2$, $n_3^{-1}(z_1)+n_3^+(z_1)\ge1$ and $b_1\in N_4^2(z_1)$, which implies $t_{z_1}^2\ge\frac{1}{6}$ and so $f_3(y)\ge\frac{1}{6}$.

Hence,
$G$ has 2-path $P(yz_1)$ and $P(yy_2)$ because $\{G+ww_1,G+w_2z\}\subseteq \mathcal{F}_G$ and $y_2z_1\notin E(G)$. Let $P(yz_1)=yz_{11}z_1$ and $P(yy_2)=yy_1y_2$. By $(c)$, $z_{11}\neq b_1$. Then $y_1\in \{z_1, z_{11}\}$, else there is a $C_6=yy_1y_2b_1z_1z_{11}$ in $G$. But then $G$ has a copy of  $C_6=yz_{11}z_1y_2b_1a_1$ when $y_1=z_1$, or $C_6=yz_1z_{11}y_2b_1a_1$ when $y_1=z_{11}$, a contradiction, as asserted.
Because $G+zb_1\in \mathcal{F}_G$ and $y_2b_1\notin E(G)$, we see there are at least two edges in $G[\{z_1, a_1, c_1\}]$.
By $(b)$, $z_1c_1\in E(G)$. But this is impossible because $4$-cycle $C_4=yc_1b_1z_1$ is not an induced cycle in $G$, a contradiction. This completes the proof of Lemma~\ref{4-53-0 }.  \qed\medskip

\noindent{\bf Proof of Lemma~\ref{law2.3.1}:} By  Corollary~\ref{f2g}, $N_4^{5-}(y_i)\subseteq N_4^{*-}(y_i)$ and $n_4^{5-}(y_i)\le1$ for  any $i\in[3]$. We proceed it by contradiction. W.l.o.g., suppose $n_4^{5-}(y_1)\ne0$. Let $z_1\in N_4^{5-}(y_1)$. By $(\divideontimes)$, $f_j(z_1)<0$ for any $j\in[5]$. By Lemma~\ref{4-53-0 }, $N_5^-(N(y_1))=\emptyset$. By Observation~\ref{ob3}(2), $f_1(z)\ge g^*(z)\ge\frac{1}{6}n_3(z)+\frac{1}{6}n_4(z)$ because $n_3(z)\ge3$, which implies  $t_{z}^2\ge\frac{1}{6}$.  We claim $d(z_1)=1$. Suppose not.
By Lemma~\ref{4-5-0}($a$), $d(z_1)=2$.
One can easily check $f_3(z_1)\ge g^*(z_1)\ge-\frac{1}{6}$ as $n_5^-(z_1)=0$. Note that $f_3(y_1)\ge f_2(y_1)+t_{z}^2\ge\frac{1}{6}$.
But then $f_4(z_1)\ge f_3(z_1)+\frac{1}{6}\ge0$, a contradiction, as claimed.  Hence,  $f_3(z_1)=g(z_1)=-\frac{1}{3}$ which means $f_3(y_1)<\frac{1}{3}$.
Clearly, $g(y_i)<0$ for any $i\in[3]$ and $n_3^{-1}(z)+n_4(z)\le1$, otherwise by Observation~\ref{ob2}(2,4), we see $g^*(z)\ge\frac{1}{3}n_3(z)+\frac{1}{6}n_4(z)$ which implies $t_z^2\ge\frac{1}{3}$ and so $f_3(y_1)\ge f_2(y_1)+t_z^2\ge\frac{1}{3}$.
Note that $L(y_i)\ne\emptyset$ for some $i\in\{2,3\}$, otherwise by Definition~\ref{law2}(3.1), $f_3(y_1)\ge f_2(y_1)+\frac{1}{2}f_2(z)\ge\frac{1}{3}$ because $f_2(z)\ge\frac{2}{3}n_3(z)-\frac{4}{3}\ge\frac{2}{3}$.  W.l.o.g., let $L(y_2)=\{z_2\}$.
Then $G$ has a $3$-path $P(y_1y_2)$ because $G+z_1z_2\in\mathcal{F}_G$. Let $P(y_1y_2)=y_1a_1a_2y_2$. We see $a_1,a_2\in V_4$ because $g(y_i)<0$. Moreover, $\{a_1,a_2\}\cap V_4^1\ne\emptyset$ as $n_3^{-1}(z)+n_4^2(z)\le1$. We next prove several claims.\medskip

\noindent{\bf Claim 1.} $a_1\ne z$.
 \medskip

\noindent{\bf Proof.} Suppose not. Then $a_2\in V_4^1$. Then $g^*(a_2)=\frac{1}{3}$ as $z\in N_4^2(y_2)\cap N_4^2(a_2)$, which means $f_1(a_2)\ge0$ and so $a_2\in N_4^{1+}(z)$. Hence, $n_4^{1-}(z)\le n_4(z)-1$. By Ineq. (\ref{g*}), $g^*(z)\ge\frac{1}{3}n_3(z)+\frac{1}{6}(n_4(z)-1)$ since $n_3(z)\ge3$ and $n_4(z)\ge1$. Hence, $t_z^2\ge\frac{1}{3}$. Then $f_3(y_1)\ge f_2(y_1)+t_z^2\ge\frac{1}{3}$, a contradiction.\qed\medskip
%\noindent{\bf Proof.}  Suppose $z=a_i$ for some $i \in [2]$. Then $g^*(a_{3-i})\ge\frac{1}{3}$ because $z\in N_4^2(y_{3-i})\cap N_4^2(a_{3-i})$, which means $f_1(a_{3-i})\ge0$ and so $a_{3-i}\in N_4^{1+}(a_i)$. Hence, $n_4^{1-}(a_i)\le n_4(a_i)-1$. By inequality~(\ref{g*}), $g^*(z)\ge\frac{1}{3}n_3(z)+\frac{1}{6}(n_4(z)-1)$ because $n_3(z)\ge3$ and $n_4(z)\ge1$. Hence, $f_2(z)\ge\frac{1}{3}n_3(z)$ and $t_z^2\ge\frac{1}{3}$. But then $f_3(y_1)\ge g^*(y_1)+t_z^2\ge\frac{1}{3}$, a contradiction.\qed

By Claim 1, we see $t_{a_1}^2<\frac{1}{6}$, otherwise $f_3(y_1)\ge f_2(y_1)+t_z^2+t_{a_1}^2\ge\frac{1}{3}$.\medskip

\noindent{\bf Claim 2.} $a_2\in V_4^1$ and $a_2\in N_4^{1+}(a_1)$.\medskip

\noindent{\bf Proof.} Suppose $a_2\notin V_4^1$. Then $a_1\in V_4^1$. Hence, $g^*(a_1)\ge\frac{1}{3}$ because $a_2\in N_4^2(a_1)$ and $z\in N_4^2(y_1)$. But then $f_2(a_1)\ge\frac{1}{6}$ and so $t_{a_1}^2\ge\frac{1}{6}$, a contradiction.   Thus, $a_2\in V_4^1$.  Clearly, $g^*(a_2)\ge\frac{1}{6}$.
 Now we prove $a_2\in N_4^{1+}(a_1)$. Suppose not. Clearly, $g(a_1)=\frac{1}{6}$ and $n_5^-(a_2)\ne0$. Let $w_2\in N_5^-(a_2)$. Then $G$ has a $3$-path $P(y_2a_2)$ because $G+z_2w_2\in\mathcal{F}_G$. Let $P(y_2a_2)=y_2b_1b_2a_2$. Then $b_2\in V_5$ and $b_1\in V_4$ as $g(a_2)=g(a_1)=\frac{1}{6}$ and $g(y_i)<0$ for any $i\in[2]$. By Observation~\ref{ob3}(2), $g^*(b_2)\ge\frac{1}{6}n_4(b_2)+\frac{1}{6}n_5(b_2)$ as $n_4^+(b_2)+n_4(b_2)\ge3$, which means $t_{b_2}^*\ge\frac{1}{6}$. But then $f_1(a_2)\ge g^*(a_2)-\frac{1}{3}n_5^-(a_2)+t_{b_2}^*\ge0$, a contradiction. \qed\medskip

By Claims 1 and 2, we see $n_3(a_1)=1$, otherwise $g^*(a_1)\ge\frac{1}{6}n_3(a_1)+\frac{1}{6}(n_4(a_1)-1)$ because $n_3(a_1)\ge2$ and $y_1\in N_3^{-1}(a_1)$, which implies $t_{a_1}^2\ge\frac{1}{6}$. Then $G$ has a $4$-path  with ends $y_1$ and $a_2$ containing $a_1$ because $G+z_1a_2\in\mathcal{F}_G$. Note that $N(y_1)\cap N(a_1)=\emptyset$, else let $b_1\in N(y_1)\cap N(a_1)$, we see $G$ has a copy of $C_6=y_1b_1a_1a_2y_2x$ when $b_1\ne a_2$ or $C_6=y_1b_1y_2xy_3z$ when $b_1=a_2$. Thus, $G-a_2$ has a $3$-path $P(y_1a_1)$ or $G-y_1$ has a $3$-path $P(a_1a_2)$. We next prove that \medskip

\noindent{$(a)$ $G-a_2$ has no $3$-path $P(y_1a_1)$.}\medskip

To prove $(a)$, suppose $3$-path $P(y_1a_1)$ exists and let $P(y_1a_1)=y_1c_1c_2a_1$. Then $c_2\in V_4\cup V_5$ as $n_3(a_1)=1$. We assert $c_2\in V_5$. Suppose not. Then $c_1\in V_4$ because $g(y_1)<0$. By Ineq.~(\ref{g*}), $g^*(a_1)\ge \frac{1}{6}n_3(a_1)+\frac{1}{6}(n_4(a_1)-1)$ because $c_2\in N_4^2(a_1)$ and $y_1\in N_3^{-1}(a_1)$. By Claim 2, $n_4^{1-}(a_1)\le n_4(a_1)-1$. But then  $t_{a_1}^2\ge\frac{1}{6}$, a contradiction, as asserted. Then $c_1\in V_4$. By Observation~\ref{ob3}(2), $g^*(c_2)\ge\frac{1}{6}n_4(c_2)+\frac{1}{6}n_5(c_2)$ because $n_4^+(c_2)+n_4(c_2)\ge3$, which means $t_{c_2}^*\ge\frac{1}{6}$. Thus, $f_1(a_1)\ge g^*(a_1)+t_{c_2}^*$. One can easily check $g^*(a_1)\ge \frac{1}{6}n_3(a_1)+\frac{1}{6}(n_4(a_1)-2)$. But then $t_{a_1}^2\ge\frac{1}{6}$ because $n_4^{1-}(a_1)\le n_4(a_1)-1$, a contradiction. This proves $(a)$. \medskip

By ($a$), $G-y_1$ has a $3$-path $P(a_1a_2)$ and let $P(a_1a_2)=a_1c_1c_2a_2$. Then $c_1\in V_4\cup V_5$ as $n_3(a_1)=1$. We assert $c_1\in V_5$. Suppose not. Then $c_2\in V_3\cup V_5$ since $a_2\in V_4^1$. If $c_2\in V_3$, then $c_2=y_2$, which implies $g^*(c_1)\ge\frac{1}{3}$ as $a_1\in N_4^2(c_1)$ and $z\in N_4^2(y_2)$. Hence, $f_1(c_1)\ge0$. If $c_2\in V_5$, then by Observation~\ref{ob2}(2), $g^*(c_2)\ge \frac{1}{3}n_4(c_2)+\frac{1}{6}n_5(c_2)$ as $\{c_1,a_2\}\subseteq N_4^+(c_2)$, which means $t_{c_2}^*\ge\frac{1}{3}$. Hence, $f_1(c_1)\ge0$. By Claim 2, $\{c_1,a_2\}\subseteq N_4^{1+}(a_1)$, which yields $n_4^{1-}(a_1)\le n_4(a_1)-2$. But then $t_{a_1}^2\ge\frac{1}{6}$ because $g^*(a_1)\ge \frac{1}{6}n_3(a_1)+\frac{1}{6}(n_4(a_1)-2)$, a contradiction, as asserted. Then $c_2\in V_5$ as $a_2\in V_4^1$. By Ineq.~(\ref{g*}), $g^*(c_1)\ge\frac{1}{6}n_4(c_1)$ since $a_1\in N_4^+(c_1)$, which means $t_{c_1}^*\ge\frac{1}{6}$. Hence, $f_1(a_1)\ge g^*(a_1)+t_{c_1}^*$. But then $t_{a_1}^2\ge\frac{1}{6}$ because $g^*(a_1)\ge \frac{1}{6}n_3(a_1)+\frac{1}{6}(n_4(a_1)-2)$, a contradiction. This completes the proof of Proposition~\ref{law2.3.1}.\qed

\noindent{\bf Proof of Lemma~\ref{4-0}:}
Clearly, we only need to prove $n_4^{5-}(y)=0$. Suppose not. Then $f_3(y)<\frac{1}{3}$ and $f_4(y)<\frac{1}{3}$. By Lemma~\ref{4-53-0 }, $N_5^-(N_4(y))=\emptyset$.
Let $z\in L(y)$, $x \in N_2(y)$ and $N_1(x)=\{\alpha_1\}$. For any $v\in V_3$, $N_4^{3-}(v)\subseteq N_4^{*-}(v)$ and $n_4^{3-}(v)\le1$ by Corollary~\ref{f2g}. We first prove several Claims. \medskip

\noindent{\bf Claim 1.} If $g(y)=\frac{1}{6}$, then $n_3^2(y)=0$,  or  $g(x)<0$ and $n_3^{2}(x)=0$. If $g(y)\ge\frac{2}{3}$, then $n_2^+(y)=0$ and $n_3^2(y)+n_2^{-1}(y)\le1$.
\medskip

\noindent{\bf Proof.} Obviously, the results hold since otherwise  $f_3(y)\ge g^*(y)\ge\frac{1}{3}$.\qed\medskip

\noindent{\bf Claim 2.} $G$ contains no  $2$-path $P(yv)$ for some $v\in\{x,\alpha_1\}$.
\medskip

\noindent{\bf Proof.} Suppose $G$ contains $2$-path $P(yv)$ for any $v\in\{x,\alpha_1\}$. Then $g(x)\ge0$ and $g(y)\ge0$. Let $P(xy)=xa_1y$.  By Claim 1, $y,a_1\in V_3^1$. Hence, $f_3(y)\ge f_2(y)=g^*(y)=\frac{1}{6}$.
Then $G-a_1$ has a $3$-path $P(xy)$ or $G-y$ has a 3-path  $P(xa_1)$ because $G+za_1\in\mathcal{F}_G$.
We assert $G-a_1$ has no $3$-path $P(xy)$. Suppose not. Let $P(xy)=xb_1b_2y$.  Then $b_1 \in V_3$ and $b_2\in V_4$ because $y\in V_3^1$. By Observation~\ref{ob3}(2), $g^*(b_2) \ge \frac{1}{6}n_3(b_2)+\frac{1}{6}n_4(b_2)$ because $n_3(b_2)\ge2$ and $y\in N_3^+(b_2)$. Since $n_5^-(b_2)=0$, we have  $t_{b_2}^2\ge\frac{1}{6}$. But then  $f_3(y) \ge f_2(y)+t_{b_2}^2\ge\frac{1}{3}$, a contradiction, as asserted. Hence, $G-y$ has a $3$-path $P(xa_1)$. Let $P(xa_1)=xb_1b_2a_1$. Similarly, $b_1\in V_3$, $b_2\in V_4$. Then $n_4^{3-}(a_1)\ne 0$, else $f_4(y)\ge f_3(y)+\frac{1}{6}\ge\frac{1}{3}$ since $f_3(a_1)\ge g^*(a_1)=\frac{1}{6}$.  Let $z_1\in N_4^{3-}(a_1)$.  Then $G$ has a $4$-path  with ends $y$ and $z_1$ containing $a_1$
as $G+zz_1\in\mathcal{F}_G$. Then $G$ has a 3-path $P(ya_1)$ because $g^*(z_1)<0$.
 Let $P(ya_1)=yc_1c_2a_1$. Then $c_1,c_2\in V_4$. Clearly, $G-y-c_2$ has an $s$-path with ends $c_1$ and $a_1$ for some $s\in[2]$ since $G+zc_2\in\mathcal{F}_G$. Then $n_3(c_1)+n_4(c_1)\ge3$ and $y\in N_3^+(c_1)$. By Ineq.~(\ref{g*}), $f_1(c_1)\ge g^*(c_1)\ge\frac{1}{6}n_3(c_1)+\frac{1}{6}(n_4(c_1)-1)$ as $n_5^-(c_1)=0$. Note that $g^*(c_2)\ge\frac{1}{3}$ because $a_1\in N_3^+(c_2)$ and $c_1\in N_4^2(c_2)$, which means $f_1(c_2)\ge0$ and so $c_2\in N_4^{1+}(c_1)$. Hence, $n_4^{1-}(c_1)\le n_4(c_1)-1$. Then, $t_{c_1}^2\ge\frac{1}{6}$. But then $f_3(y)\ge f_2(y)+t_{c_1}^2\ge\frac{1}{3}$, a contradiction. This proves Claim 2.\qed\medskip

Note that $G$ has a $4$-path with ends $\alpha_1$ and $y$ containing $x$ since $G+z\alpha_1\in\mathcal{F}_G$. By Claim 2, $G-y$ has a $3$-path with ends $x$ and $\alpha_1$ or $G-\alpha_1$ has a $3$-path with ends $x$ and $y$. We next prove the following claim.\medskip

\noindent{\bf Claim 3.} $G-\alpha_1$ has no $3$-path with ends $x$ and $y$.\medskip

\noindent{\bf Proof.} Suppose not. Let $P(xy)=xx_1x_2y$. By Claim 1, $xx_2\notin E(G)$. Then $G$ has a $4$-path with ends $y$ and $x_1$ containing $x$ and $x_2$ as $G+zx_1\in\mathcal{F}_G$. We thus see that $G-y-x_1$ has a $2$-path $P(xx_2)$.
Let $P(xx_2)=xx_3x_2$. By Claim 1, $x_1,x_3\in V_3$ and  $x_2\in V_3\cup V_4$.   We next prove that \medskip

\noindent{$(a)$ $x_2\in V_3$ and $L(x_2)=\emptyset$.}\medskip

To see why $(a)$ is true, suppose first $x_2\in V_4$.
Clearly, $N_3(x_2)=\{y,x_1,x_3\}$, else by Observation~\ref{ob2}(2), $ g^*(x_2)\ge\frac{1}{3}n_3(x_2)+\frac{1}{6}n_4(x_2)$ which implies $t_{x_2}^2\ge\frac{1}{3}$ and so $f_3(y)\ge f_2(y)+t_{x_2}^2\ge\frac{1}{3}$.  Moreover, we see $N_2(y)=N_2(x_1)=N_2(x_3)=\{x\}$, else we assume $x'\in N_2(x_1)\setminus x$, we see $G$ has a copy of $C_6=\alpha_1xx_3x_2x_1x'$. But then $n_4^{5-}(y)=0$ by Lemma~\ref{law2.3.1}, a contradiction. Thus, $x_2\in V_3$.
  We now show that
$L(x_2)=\emptyset$.
Suppose not. Let $z_1\in L(x_2)$. Then $z_1\in V_4$. Then $G$ has a $4$-path with ends $x$ and $x_2$ containing each of $\{y,x_1,x_3\}$ because $G+xz_1\in\mathcal{F}_G$. Hence, $G[\{y, x_1, x_3\}]$ has at least two edges. W.l.o.g., let $yx_3\in E(G)$. Then $y\in V_3^2$ and $\{x_2, x_3\}\subseteq N_3^2(y)$, which violates Claim 1.
%
%By Claim 1, $x_2\in V_4$, else $g(y)\ge\frac{2}{3}$ and  $n_3^2(y)+n_2^+(y)+n_2^{-1}(y)\ge3$. Then $x_1\in V_3$ and $n_3^+(x_2)\ge2$. By Observation~\ref{ob2}(1), $f_1(x_2)\ge g^*(x_2)\ge\frac{1}{3}n_3(x_2)+\frac{1}{6}n_4(x_2)$ because $n_5^-(x_2)=0$. Hence, $f_2(x_2)\ge\frac{1}{3}n_3(x_2)$ and so $t_{x_2}^2\ge\frac{1}{3}$. But then $f_3(y)\ge g^*(y)+t_{x_2}^2\ge\frac{1}{3}$, a contradiction.
 This proves $(a)$.\medskip

 By $(a)$ and Claim 1, $g(x)<0$ and
$x_1, x_3\in V_3^1$. For any $v\in V_3$, let $A(v)$ be defined as in Definition~\ref{law2}(4).
By $(a)$ and Lemma~\ref{4-5-0}($a,b$),  $f_3(v)\ge-\frac{1}{6}$ for any $v\in N_4^{3-}(x_2)$.  Hence, $f_3(x_2)<\frac{1}{6}|A(x_2)|+\frac{1}{6}n_4^{3-}(x_2)$, else $f_4(y)\ge f_3(y)+\frac{1}{6}\ge\frac{1}{3}$.
%Thus $f_3(x_2)<\frac{1}{6}|A(x_2)|+\frac{1}{6}$ because $n_4^{3-}(x_2)\le 1$.
Since $g(x_2)\ge\frac{2}{3}$, we see $f_3(v)\ge g^*(v)\ge\frac{1}{6}$ for any $v\in\{y,x_1,x_3\}$. By Ineq.~(\ref{g*}), $g^*(x_2)\ge\frac{1}{6}(n_3(x_2)-1)$ because  $n_2(x_2)+n_3(x_2)\ge4$.
We assert $g(y)=\frac{1}{6}$. Suppose not.
For any $i\in\{1,3\}$, $f_3(x_i)\ge g^*(x_i)=\frac{1}{3}$ because $y\in N_3^2(x)$ and $x_2\in N_3^2(x_i)$, which yields $f_3(x_i)+f_3(x_i')\ge0$ for any $x_i'\in N_4(x_i)$. Hence, $x_1,x_3\notin A(x_2)$, which means $|A(x_2)|\le n_3(x_2)-2$. But then $f_3(x_2)\ge g^*(x_2)\ge\frac{1}{6}|A(x_2)|+\frac{1}{6}\ge\frac{1}{6}|A(x_2)|+\frac{1}{6}n_4^{3-}(x_2)$, a contradiction, as asserted.
We further assert $L(x_1)\cup L(x_3)\ne\emptyset$.
Suppose not. By Lemma~\ref{4-5-0}($a,b$), $f_3(v)\ge-\frac{1}{6}$ for any $v\in N_4^{3-}(x_i)$, which yields $f_3(x_i)+f_3(x_i')\ge0$ for any $x_i'\in N_4(x_i)$ and $i\in\{1,3\}$. Hence, $x_1,x_3\notin A(x_2)$. Similarly, we also get a contradiction, as asserted.
W.l.o.g., let $z_1\in L(x_1)$.
Because $G+zz_1\in\mathcal{F}_G$, we see $G$ has a 3-path $P(yx_1)$. Let $P(yx_1)=ya_1a_2x_1$.  We next show that\medskip

\noindent{$(b)$ $a_1,a_2\in V_4^1$.}\medskip

To prove $(b)$, suppose first $a_1\notin V_4^1$. We assert $a_1\in V_4^2$. Suppose not. Then $a_1=x_2$ and $a_2\in V_4$ because $g(y)=g(x_1)=\frac{1}{6}$ and $xx_2\notin E(G)$. By Ineq.~(\ref{g*}), $g^*(a_2)\ge\frac{1}{6}n_3(a_2)+\frac{1}{6}(n_4(a_2)+2)$ because $x_1,x_2\in N_3^+(a_2)$, which means $t_{a_2}^2\ge\frac{1}{6}$. Hence,  $f_3(x_1)\ge f_2(x_1)+t_{a_2}^2\ge \frac{1}{3}$. Then $x_1\notin A(x_2)$, which means $|A(x_2)|\le n_3(x_2)-1$. But then $f_3(x_2)\ge f_2(x_2)+t_{a_2}^2\ge g^*(x_2)+\frac{1}{6}\ge \frac{1}{6}|A(x_2)|+\frac{1}{6}n_4^{3-}(x_2)$, a contradiction, as asserted. Thus, $a_1\in V_4^2$.
Note that $a_2\in V_4$, else $a_2=x_2$ yields $t_{a_1}^2\ge\frac{1}{6}$ and so $f_3(y)\ge f_2(y)+t_{a_1}^2\ge\frac{1}{3}$.  Then $g^*(a_2)\ge\frac{1}{3}$ because $x_1\in N_3^+(a_2)$ and $a_1\in N_4^2(a_2)$, which means $a_2\in N_4^{1+}(a_1)$ and so $n_4^{1-}(a_1)\le n_4(a_1)-1$. By Ineq.~(\ref{g*}),
$g^*(a_1)\ge\frac{1}{6}n_3(a_1)+\frac{1}{6}(n_4(a_1)-1)$ because $y\in N_3^+(a_1)$.
Then, $t_{a_1}^2\ge\frac{1}{6}$ because $n_5^-(a_1)=0$. But then $f_3(y)\ge f_2(y)+t_{a_1}^2\ge\frac{1}{3}$, a contradiction. Thus $a_1\in V_4^1$. Then $a_2\in V_4$. Moreover, $a_2\in V_4^1$, otherwise  $f_1(a_1)\ge g^*(a_1)\ge\frac{1}{3}$ because $y\in N_3^+(a_1)$ and $n_5^-(a_1)=0$, which means $t_{a_1}^2\ge \frac{1}{6}$ and so $f_3(y)\ge f_2(y)+t_{a_1}^2\ge\frac{1}{3}$. This proves $(b)$.\medskip

By $(b)$, $g^*(a_1)=\frac{1}{6}$ as $y\in N_3^+(a_1)$. Then $G$ has a $4$-path with ends $y$ and $a_2$ containing $a_1$ because $G+za_2\in\mathcal{F}_G$. Since $g(a_1)=\frac{1}{6}$, we see $G-a_2$ has a $3$-path $P(a_1y)$ or $G-y$ has a $3$-path $P(a_1a_2)$. Let $P(a_1v)=a_1b_1b_2v$, where $v\in \{y,a_2\}$. Then $b_1\in V_5$ because $g(a_1)=\frac{1}{6}$. Clearly, $g^*(b_1)\ge\frac{1}{6}n_4(b_1)$ since $a_1\in N_4^+(b_1)$, which means $t_{b_1}^*\ge\frac{1}{6}$. Thus, $f_1(a_1)\ge g^*(a_1)+t_{b_1}^*\ge\frac{1}{3}$ as $n_5^-(a_1)=0$, which means $t_{a_1}^2\ge\frac{1}{6}$. But then $f_3(y)\ge f_2(y)+t_{a_1}^2\ge\frac{1}{3}$, a contradiction.\qed

\vspace{0.35cm}

By Claim 3, $G-y$ has a $3$-path $P(x\alpha_1)$ and let $P(x\alpha_1)=xa_2a_1\alpha_1$. Then $a_1\in V_2$ since $\delta(G)=1$. We assert $n_2(y)=1$. Suppose not. Let $x'\in N_2(y)\setminus x$. By Claim 1, $x'\neq a_2$.  Hence, $x'=a_1$, else $G$ has a copy of $C_6=\alpha_1a_1a_2xyx'$. But then $n_2^+(y)+n_2^{-1}(y)\ge2$, which violates  Claim 1, as asserted. We next prove that \medskip

\noindent{$(c)$ $n_3(y)=0$.}\medskip

To see why $(c)$ is true, suppose $n_3(y)\ne0$. By Claim 1, $n_3^2(y)=0$. So $N_3(y)=N_3^1(y)$. Let $y_1\in N_3(y)$ and $x_1\in N_2(y_1)$.
Clearly, $f_3(y)\ge g^*(y)\ge\frac{1}{6}$ because $x\in N_2^+(y)$ or $N_3^2(x)\setminus y\ne\emptyset$.

We first show that $G-x_1$ has no $3$-path with ends $y$ and $y_1$. Suppose not. Let $P(yy_1)=yc_1c_2y_1$. Then $c_2\in V_4$ because $g(y_1)=\frac{1}{6}$.
We claim $c_1\in V_4$. Suppose not. Then $c_1\in V_3$. Clearly, $N_2(c_1)=\{x_1\}$. Hence, $f_2(y_1)\ge g^*(y_1)\ge\frac{1}{6}$ because $y\in N_3^2(y_1)$.  By Ineq.~(\ref{g*}), $g^*(c_2)\ge\frac{1}{6}n_3(c_2)+\frac{1}{6}(n_4(c_2)+2)$ because $c_1,y_1\subseteq N_3^+(c_2)$, which implies $t_{c_2}^2\ge\frac{1}{6}$. Hence, $f_3(y_1)\ge f_2(y_1)+t_{c_2}^2\ge\frac{1}{3}$. Moreover, $L(y_1)\ne\emptyset$, otherwise by
Lemma~\ref{4-5-0}($a,b$),  $f_3(v)\ge-\frac{1}{6}$ for any $v\in N_4^{3-}(y_1)$, which yields $f_4(y)\ge f_3(y)+\frac{1}{6}\ge\frac{1}{3}$. Let $z'\in L(y_1)$. We see $G-c_1-y_1$ has a path of length at most two with ends $y$ and $c_2$ because $G+z'c_1\in \mathcal{F}_G$, which implies $t_{c_2}^2\ge\frac{1}{3}$ and so $f_3(y_1)\ge f_2(y_1)+t_{c_2}^2\ge\frac{1}{2}$. But then $f_4(y)\ge f_3(y)+\frac{1}{6}\ge\frac{1}{3}$, a contradiction, as claimed.
%\textcolor{blue}{We assert $c_1\in V_4$. Suppose not. Then $c_1\in V_3$. Then there is another vertex $c_3\in N_4(c_1)\cap N_4(y_1)\setminus c_2$ because
% $G+zc_2\in \mathcal{F}_G$ and $y_1\in V_3^1$.
%  By inequality~(\ref{g*}), $g^*(c_i)\ge\frac{1}{6}n_3(c_i)+\frac{1}{6}(n_3(c_i)+2)$ because $\{c_1,y_1\}\subseteq N_3^+(c_i)$ for any $i \in \{2,3\}$, which implies $t_{c_i}^2\ge\frac{1}{6}$. Hence, $f_3(y_1)\ge f_2(y_1)+t_{c_2}^2+t_{c_3}^2\ge\frac{1}{3}$. Moreover, $L_4(y_1)\ne\emptyset$, otherwise by
%Lemma~\ref{4-5-0}($a,b$),  $f_3(z_1)\ge-\frac{1}{6}$ for $z_1\in N_4^{3-}(y_1)$, which yields $f_4(y)\ge f_3(y)+\frac{1}{6}\ge\frac{1}{3}$.
%Let $z'\in L_4(y_1)$, we see $G-c_1-y_1$ has a path of length at most two with ends $y$ and $c_2$ because $G+z'c_1\in \mathcal{F}_G$, which implies $t_{c_2}^2\ge\frac{1}{3}$ and so $f_3(y_1)\ge f_2(y_1)+t_{c_2}^2+t_{c_3}^2\ge\frac{1}{2}$. But then $f_4(y)\ge f_3(y)+\frac{1}{6}\ge\frac{1}{3}$, a contradiction, as claimed.}
Note that $y_1c_1\notin E(G)$, otherwise $f_1(c_1)\ge g^*(c_1)\ge\frac{1}{3}n_3(c_1)+\frac{1}{6}n_4(c_1)$ because $y,y_1\subseteq N_3^+(c_1)$ and $n_5^-(c_1)=0$, which means $t_{c_1}^2\ge\frac{1}{3}$ and so $f_3(y)\ge f_2(y)+t_{c_1}^2\ge\frac{1}{3}$. Then
$G-y-c_2$ has a $2$-path $P(y_1c_1)$ as $G+zc_2\in\mathcal{F}_G$. Let $P(y_1c_1)=y_1d_1c_1$.
Then $d_1\in V_4$ as $g(y_1)=\frac{1}{6}$. By Ineq.~(\ref{g*}), $f_1(c_1)\ge g^*(c_1)\ge\frac{1}{6}n_3(c_1)+\frac{1}{6}(n_4(c_1)-1)$ since $y\in N_3^+(c_1)$ and $n_5^-(c_1)=0$. Clearly, $f_1(c_2)\ge0$, which means $c_2\in N_4^{1+}(c_1)$ and so $n_4^{1-}(c_1)\le n_4(c_1)-1$. Hence, $t_{c_1}^2\ge\frac{1}{6}$.
But then $f_3(y)\ge f_2(y)+t_{c_1}^2\ge\frac{1}{3}$, a contradiction.

Thus, $G$ has no $3$-path with ends $y$ and $y_1$. Clearly, $n_4^{3-}(y_1)=0$, else let $v\in N_4^{3-}(y_1)$, $G+zv\notin \mathcal{F}_G$ by Lemma~\ref{4-5-0}($a,b$).
Since $G+zx_1\in\mathcal{F}_G$ and $y_1\in N_3^1(y)$, we see $G$ has a $3$-path $P(y_1x_1)=y_1c_1c_2x_1$ such that $c_1\in V_4$ and $c_2\in V_3$.
If $n_5^-(c_1)=0$, then by Observation~\ref{ob3}(2), $ g^*(c_1)\ge\frac{1}{6}n_3(c_1)+\frac{1}{6}n_4(c_1)$ because $y_1\in N_3^+(c_1)$ and $n_3(c_1)\ge2$, which means $t_{c_1}^2\ge\frac{1}{6}$ and so $f_3(y_1)\ge f_2(y_1)+t_{c_1}^2\ge\frac{1}{6}$.  If $n_5^-(c_1)\ne0$, then let $w\in N_5^-(c_1)$. Then $G-x_1-c_1$ has a $2$-path $P(y_1c_2)=y_1c_3c_2$ such that $c_3\in V_4$ because $G+wx_1\in\mathcal{F}_G$ and $y_1\in N_3^1(y)$. Hence, for $i\in\{1,3\}$,  $g^*(c_i)\ge \frac{1}{12}n_3(c_i)+\frac{1}{6}n_4(c_i)+\frac{1}{3}$ because $y_1\in N_3^+(c_i)$, $n_3(c_i)\ge2$ and $c_2\in N_3^+(c_i)\cup N_3^{-1}(c_i)$. This implies $t_{c_i}^2\ge\frac{1}{12}$. Thus, $f_3(y_1)\ge f_2(y_1)+t_{c_1}^2+t_{c_3}^2 \ge\frac{1}{6}$.
But then in both cases $f_4(y)\ge f_3(y)+\frac{1}{6}\ge\frac{1}{3}$, a contradiction. This proves $(c)$.\medskip

Let $N_1(\alpha_1)=\{\alpha\}$. By the choice of $\alpha$, $y$ belongs to some $4$-cycle with vertices $y,z_1,z_2,z_3$ in order. By $(c)$ and Claim 3, $z_1,z_3\in V_4$.  Note that $n_5^-(z_i)=0$ for any $i\in\{1,3\}$, which means $f_1(z_i)\ge g^*(z_i)$. Since $G+zz_2\in\mathcal{F}_G$, we have $z_1z_3\in E(G)$ or $G-y-z_2$ has a $2$-path $P(z_1z_3)=z_1z_4z_3$.  We assert $g(z_i)\ge0$ for any $i\in\{1,3\}$. W.l.o.g., suppose $g(z_1)<0$. Then $z_2,z_4\in V_5$. By Ineq.~(\ref{g*}), for any $i\in\{2,4\}$, $g^*(z_i)\ge\frac{1}{6}n_4(z_i)$ because $n_4(z_i)\ge2$ and $n_4^+(z_i)+n_4^{-1}(z_i)\ge2$, which means $t_{z_i}^*\ge\frac{1}{6}$. Hence, $f_2(z_1)\ge g^*(z_1)+t_{z_2}^*+t_{z_4}^*\ge\frac{1}{3}$, which means $t_{z_1}^2\ge\frac{1}{3}$. But then $f_3(y)\ge f_2(y)+t_{z_1}^2\ge\frac{1}{3}$, a contradiction, as asserted. Then $g^*(y)\ge0$. Obviously, $t_{z_i}^2<\frac{1}{6}$ for some $i\in\{1,3\}$, else $f_3(y)\ge f_2(y)+t_{z_1}^2+t_{z_3}^2\ge\frac{1}{3}$. We further assert $N_5(z_1)\cap N_5(z_3)\ne\emptyset$. Suppose not. Then $z_1,z_3\in V_4^2$. By Ineq.~(\ref{g*}), for any $i\in\{1,3\}$, $g^*(z_i)\ge\frac{1}{6}n_3(z_i)+\frac{1}{6}n_4(z_i)$ because $n_3^+(z_i)+n_3^{-1}(z_i)+n_4^2(z_i)\ge3$, which yields $t_{z_i}^2\ge\frac{1}{6}$, a contradiction, as asserted. W.l.o.g., let $z_2\in N_5(z_1)\cap N_5(z_3)$. By Ineq.~(\ref{g*}), $g^*(z_2)\ge\frac{1}{3}n_4(z_2)$ because $\{z_1,z_3\}\subseteq N_4^+(z_2)$, which means $t_{z_2}^*\ge\frac{1}{3}$.  Note that $t_{z_i}^2<\frac{1}{6}$ for some $i\in\{1,3\}$. W.l.o.g., suppose $t_{z_1}^2<\frac{1}{6}$.  Then $z_1\in V_4^2$, otherwise $g^*(z_1)\ge 0$ and $f_1(z_1)\ge g^*(z_1)+t_{z_2}^*\ge \frac{1}{3}$, which implies $t_{z_1}^2\ge\frac{1}{6}$. Moreover, $z_3\in V_4^2$, otherwise $g^*(z_3)\ge\frac{1}{6}$ because $z_1\in N_4^2(y)$, which implies $f_2(z_{3})\ge g^*(z_{3})+t_{z_2}^*-\frac{1}{6}\ge \frac{1}{3}$ and so $f_3(y)\ge f_2(y)+t_{z_3}^2\ge\frac{1}{3}$. It is easy to see that if $z_4$ exists, then $z_4\in V_3\cup V_4$ otherwise $t_{z_4}^{*}\ge\frac{1}{3}$ and $t_{z_1}^2\ge\frac{1}{6}$. By Ineq.~(\ref{g*}), $g^*(z_1)\ge \frac{1}{6}n_3(z_1)+\frac{1}{6}(n_4(z_1)-2)$ because $n_3^+(z_1)+n_3^{-1}(z_1)+n_4^2(z_1)\ge2$ and $n_3(z_1)+n_4(z_1)\ge 2$.  Hence, $f_1(z_1)\ge g^*(z_1)+t_{z_2}^*\ge \frac{1}{6}n_3(z_1)+\frac{1}{6}n_4(z_1)$. But then $t_{z_1}^2\ge\frac{1}{6}$, a contradiction. This completes the proof of Lemma~\ref{4-0}. \qed

\noindent{\bf Proof of Lemma~\ref{yz16-40}:}
To prove ($a$), suppose $n_4^+(y)\ne0$. Let $z_0\in N_4^+(y)$. By $(\divideontimes)$ and Corollary~\ref{f2g}, $g^*(y)\le f_3(y)<0$ and $g^*(z)<0$.
By Lemmas~\ref{4-0} and \ref{4-5-0}$(a,b)$, $G$ has a $3$-path $P=zww_1z_1$ such that $d(v)=2$ for any $v\in V(P)$, and $yz_1\notin E(G)$, where $z_1\in V_4$, $w,w_1\in V_5$. By Observation~\ref{ob}(1),  $N_4^+(y)=N_4^1(y)=\{z_0\}$, which yields $f_2(y)= g^*(y)\ge-\frac{1}{6}$. Let $N_3(z_1)=\{y_1\}$, $N_4(z_0)=\{z_2\}$ and $N_2(y)=\{x\}$. Then $G$ has a $2$-path with vertices $x,x',y_1$ in order because  $G+wx\in\mathcal{F}_G$ and $g^*(y)<0$.
Then $y_1z_2\notin E(G)$, else $y_1x'xyz_0z_2$ is a $6$-cycle. This means $G$ has no $2$-path with ends $y_1$ and $z_0$.
Then $G$ has a $3$-path $P(yz_0)$ because $G+wz_0\in\mathcal{F}_G$. Let $P(yz_0)=ya_1a_2z_0$. Then $a_1\in N_4^-(y)$ and $a_2\in V_5$.  Note that $z_0\in N_4^+(a_2)$ and $n_4(a_2)\ge2$. By Observation~\ref{ob3}(2), $g^*(a_2)
\ge\frac{1}{6}n_4(a_2)+\frac{1}{6}n_5(a_2)$ and so $t_{a_2}^*\ge\frac{1}{6}$. By Lemma~\ref{4-53-0 }, $n_5^-(a_1)=0$. Hence, $f_2(a_1)\ge g^*(a_1)+t_{a_2}^*\ge\frac{1}{6}$ and so $t_{a_1}^2\ge\frac{1}{6}$. But then $f_3(y)\ge f_2(y)+t_{a_1}^2\ge0$, a contradiction. This proves $(a)$.

Now we shall prove ($b$). Suppose $n_5^-(z_1)\ne0$. Let $w_1\in N_5^-(z_1)$.
By $(\divideontimes)$ and $g(z)=\frac{1}{6}$, $-\frac{1}{6}\le f_2(y)\le f_3(y)<0$.
Then $G$ has a $4$-path with ends $y$ and $z_1$ containing $z$ as $G+yw_1\in\mathcal{F}_G$. Clearly, $G-z_1$ has no $3$-path with ends $y$ and $z$, otherwise $f_3(y)\ge0$ by similar analysis with the case $P(yz_0)$ exists in the proof of Lemma~\ref{yz16-40}($a$). Then $G-y$ has a $3$-path $P(zz_1)$.
Let $P(zz_1)=za_1a_2z_1$.  By Observation~\ref{ob}(2), $g(z_1)=\frac{1}{6}$. Hence, $a_1,a_2\in V_5$. Then $G-z_1-a_1$ has a path of length at most two with ends $z$ and $a_2$ because $G+w_1a_1\in\mathcal{F}_G$. Hence, $a_2\in N_5^2(a_1)$, which means $g^*(a_1)\ge\frac{1}{3}n_4(a_1)$ because $z\in N_4^+(a_1)$ and so $t_{a_1}^*\ge\frac{1}{3}$. By Lemma~\ref{55-0}, $n_5^-(z)=0$.
Hence, $f_1(z)\ge g^*(z)+t_{a_1}^*\ge\frac{1}{3}$ and so $t_{z}^2\ge\frac{1}{6}$. But then $f_3(y)\ge f_2(y)+t_{z}^2\ge0$, a contradiction.\qed

\end{document}